\documentclass[12pt]{amsart}
\usepackage[english]{babel}
\usepackage[all]{xy}
 \usepackage[latin1]{inputenc}
  \usepackage[T1]{fontenc}
  \usepackage{amsmath,amssymb}
  \usepackage{amsmath}
  \usepackage{amsthm}
\usepackage{amssymb}
\usepackage{mathrsfs}
\usepackage{enumerate}
\usepackage{graphicx}
\usepackage{tikz}\parindent=0.cm
\usetikzlibrary{shapes.geometric}
\usetikzlibrary{positioning}
\usetikzlibrary{arrows}
\usetikzlibrary{shapes.multipart}
\usepackage{physics}
\usepackage[notcite, final, notref]{showkeys}
\usepackage{dsfont}

\newtheorem{theorem}{Theorem}[section]

\newtheorem{lemma}[theorem]{Lemma}
\newtheorem{proposition}[theorem]{Proposition}

\newtheorem{corollary}[theorem]{Corollary}
\newtheorem{definition}[theorem]{Definition}

\theoremstyle{definition}
\newtheorem{remark}[theorem]{Remark}

\newtheorem{example}[theorem]{Example}

\newcommand{\sgn}{{\rm sign}\,}

\newcommand{\e}{\varepsilon}

\newcommand{\ct}{\circledast}

\usepackage{tikz}
\usetikzlibrary{arrows,shapes,positioning,shadows,trees,matrix}
\usepackage{xcolor}

\begin{document}

\title[Scaled global operator]{Scaled global operators and Fueter variables on non-zero scaled hypercomplex numbers}

\author[D. Alpay]{Daniel Alpay}
\address{(DA) Schmid College of Science and Technology \\
Chapman University\\
One University Drive
Orange, California 92866\\
USA}
\email{alpay@chapman.edu}

\author[I. Cho]{Ilwoo Cho}
\address{(IC) Department of Mathematics and Statistics \\
Saint Ambrose University \\
508 W. Locust St.
Davenport, IA 52803\\
USA}
\email{choilwoo@sau.edu}

\author[M. Vajiac]{Mihaela Vajiac}
\address{(MV) Schmid College of Science and Technology \\
Chapman University\\
One University Drive
Orange, California 92866\\
USA}
\email{mbvajiac@chapman.edu}

\thanks{D. Alpay thanks the Foster G. and Mary McGaw Professorship in
  Mathematical Sciences, which supported his research.}

\maketitle

\begin{abstract}
  In this paper we describe the rise of global operators in the scaled quaternionic case, an important extension from the quaternionic case to the family of scaled hypercomplex numbers $\mathbb{H}_t,\, t\in\mathbb{R}^*$, of which the $\mathbb{H}_{-1}=\mathbb{H}$ is the space of quaternions
  and $\mathbb{H}_{1}$ is the space of split quaternions. We also describe the scaled Fueter-type  variables associated to these operators, developing a coherent theory in this field. We use these types of variables to build different types of function spaces on  $\mathbb{H}_t$. Counterparts of the Hardy space and of the Arveson space are also introduced
  and studied in the present setting.
  The two different adjoints in the scaled hypercomplex numbers lead to two parallel cases in each
  instance. Finally we introduce and study the notion of rational function.
\end{abstract}

\tableofcontents

\mbox{}\\

\newpage
\section{Introduction and the Algebra of scaled Hypercomplex numbers}
\setcounter{equation}{0}
\label{Intro-Section}
In the case of $\mathbb H$, the skew-field of quaternions, generated on the real numbers by $1,i,j,k$ with the usual Cayley table, a global operator $G$ was defined in ~\cite{GLOBAL} as:
\begin{equation}
G_q=(x_1^2+x_2^2+x_3^2)\frac{\partial}{\partial x_0}+(ix_1+jx_2+kx_3)\left(x_1\frac{\partial}{\partial x_1}
    +x_2\frac{\partial}{\partial x_2}+x_3\frac{\partial}{\partial x_3}\right).
  \end{equation}
  
  In \cite{adv_prim}, the related operator, 
\begin{equation}
  V_q =\frac{\partial}{\partial x_0}-\frac{1}{ix_1+jx_2+kx_3}\left(x_1\frac{\partial}{\partial x_1}
    +x_2\frac{\partial}{\partial x_2}+x_3\frac{\partial}{\partial x_3}\right),
  \end{equation}
(which we still name "the global operator") was studied extensively by two of the co-authors (together with Kamal Diki). In particular new Fueter-type variables associated to $V_q$ were found and analyzed, and Fueter-type expansions were given for smooth functions in its kernel.\smallskip

In recent works~\cite{alpay_cho_3,MR4658650,MR4591396,MR4592128,alpay_cho_4}, the space of quaternions was extended to a family of spaces indexed
by the real numbers, and, for  $t\in\mathbb R$, these papers introduced and studied the scaled families of rings $\mathcal H^t_2$ of $\mathbb C^{2\times 2}$ matrices,
with elements  of the form:
\begin{equation}
  q_t=\begin{pmatrix}a&tb\\ \overline{b}&\overline{a}\end{pmatrix},\quad a,b\in\mathbb C,\, t\in\mathbb R.
\end{equation}
When $t=-1$ one obtains the matrix representation of the skew field of quaternions, while the case $t=1$ corresponds to the matrix representation of the ring of split-quaternions. In the latter case ($t=1$), hyperbolic numbers correspond to real $a$ and $b$.\smallskip

In the present paper we first extend the results of \cite{adv_prim} to the scaled family of rings mentioned above, interpolating from the space of quaternions to the one of split quaternions. Then we study counterparts of the classical Hardy space in the present setting. An important new point with respect to the previous work \cite{adv_prim} is the appearance of Krein spaces (as in \cite{alss_IJM}).\smallskip

We now review the definitions and main results of the scaled-hypercomplex algebras introduced in
\cite{alpay_cho_3,MR4658650,MR4591396,MR4592128,alpay_cho_4}. We start with the first definition of the space $\mathcal H_2^t$ as above, with two
equivalent characterizations to follow:
\begin{definition}
\label{matrix-def}
Let $t\in\mathbb R$. The space of scaled hypercomplex numbers is:
\begin{equation}\label{qt890}
\mathcal H_2^t= \left\{ \begin{pmatrix}a&tb\\ \overline{b}&\overline{a}\end{pmatrix}| \,a,b\in\mathbb C\right\}.
\end{equation}
This space is an algebra over the field of real numbers, with the usual addition and multiplication of matrices.
\end{definition}
In another characterization we will see that $\mathcal H_2^t$ becomes a real vector space.
\begin{remark}
We note that $\mathcal H_2^t$ is not closed under the usual matrix conjugation when $t\not=\pm 1$, since
\begin{equation}
\label{usual*}
  q_t^*=\begin{pmatrix}\overline{a}&b\\ t\overline{b}&a\end{pmatrix}\not\in\mathcal H_2^t
\end{equation}
if $t\not =\pm 1$ or if $b\not=0$. 
\end{remark}

One can define two natural adjoints (or conjugations) associated to \eqref{qt890}, and which leave $\mathcal H_2^t$ invariant, namely
\begin{equation}
  \label{q1}
  q_t^{\circledast}=\begin{pmatrix}\overline{a}&-tb\\ -\overline{b}&a\end{pmatrix}
\end{equation}
(that is, the adjugate matrix of $q_t$) and, see \eqref{second-adjoint}, 
\begin{equation}
  \label{q2}
  q_t^{[*]}=\begin{pmatrix}a&t\overline{b}\\ b&\overline{a}\end{pmatrix}.
  \end{equation}
See \cite{alpay_cho_3,MR4591396,MR4592128} and \cite{alpay_cho_4}.
  We note that these two adjoints are contravariant (see Lemma \ref{contra-variant} and Lemma~\ref{[*]adjoint})
  \[
(q_tp_t)^{\ct}=p_t^{\ct}q_t^{\ct}\quad {\rm and}\quad (q_tp_t)^{[*]}=p_t^{[*]}q_t^{[*]},\quad p_t,q_t\in\mathcal H_2^t,
    \]

are involutive, and commute:
  \[
\left(q_t^{\ct}\right)^{[*]}=\left(q_t^{[*]}\right)^{\ct}=\begin{pmatrix}\overline{a}&-t\overline{b}\\-b&a\end{pmatrix}.
    \]
It is therefore not possible to obtain another involutive contravariant from iterated compositions of these adjoints.\\
    
  It is easy to see that $\mathcal H_2^t$ is ring-isomorphic to the ring $\mathbb H_t$ defined below.
  \begin{definition} Let  $\mathbb H_t$ be the space of elements of the form:
\begin{equation}
  \label{new-def}
  q_t=a+bj_t,\quad a,b\in\mathbb C[i],
  \end{equation}
  where $i, j_t$ and $k_t=ij_t$ satisfy the generalized Cayley laws, induced by the matrix multiplication in Definition~\ref{matrix-def}:
  
\begin{center}
  \begin{tabular}{|c|c|c|c|c|}
    \hline
    $\nearrow$&$1$&$i$& $j_t$&$k_t$\\
    \hline
    $1$&$1$&$i$& $j_t$&$k_t$\\
    \hline
    $i$&$i$&$-1$ &$k_t$ &$-j_t$\\
    \hline
    $j_t$& $j_t$&$-k_t$&$t$&$-ti$ \\
    \hline
    $k_t$&$k_t$ &$j_t$ &$ti$ &$t$\\
    \hline
  \end{tabular}.
  \end{center}
 Writing $a=x_0+x_1 i$ and $b=x_2+x_3 i$ we have that the ring $\mathbb{H}_{t}$ is the vector space over $\mathbb R$ of elements of the form:
\begin{equation}\label{new-def-vs}
  q_t=x_0+x_1i+x_2j_{t}+x_3k_{t},\quad x_0,\;x_1,\;x_2,\;x_3\in\mathbb{R}.
\end{equation}
   \end{definition}
  
   We will often in this paper pass from one representation to the other and denote both of them by $q_t$.\smallskip
   
We now introduce a third equivalent definition, where scaled hypercomplex numbers are seen as pairs of complex numbers. For an arbitrarily fixed non-zero scale $t\in\mathbb{R}$, and $q_t=a+bj_t$ as in~\eqref{new-def}, we write the corresponding pair of complex numbers $\left(a,b\right)_t$ and we define the scaled multiplicative operation $\cdot_{t}$ on $\mathbb{C}^{2}$ induced by the multiplication of matrices  as:
\begin{equation}
\label{dot-t}
  \left(a_{1},b_{1}\right)_t\cdot_{t}\left(a_{2},b_{2}\right)_t\overset{\textrm{def}}{=}\left(a_{1}a_{2}+tb_{1}\overline{b_{2}},\:a_{1}b_{2}+b_{1}\overline{a_{2}}\right)_t.
\end{equation}

As proven in~\cite{MR4591396,MR4592128} the space of pairs $(a,b)_t$  forms
a unital ring with unity $\left(1,0\right)$, where $+$ is
the usual vector addition, and $\cdot_{t}$
is the operation \eqref{dot-t} above.

We will use these three definitions interchangeably, whenever each becomes more convenient to describe the associated differential theory, spaces of functions, domains of holomorphy and so forth. \smallskip

  Observe that, if
\[
a=x_0+x_1i \;\,\,\mathrm{and\;}\,\,b=x_2+x_3i
\]
are complex numbers with $i=\sqrt{-1}$ and $x_0,x_1,x_2,x_3\in\mathbb{R}$,
and $\left(a,b\right)_t \in\mathbb{H}_{t}$, with some abuse of notation, we have that the realization of
$\left(a,b\right)_t$ is identified with:

\begin{equation}
  \label{complex-pair}
q_t=a+bj_t=\left(a,b\right)_{t}=\left(\begin{array}{cc}
a & tb\\
\overline{b} & \overline{a}
\end{array}\right)=\left(\begin{array}{ccc}
x_0+x_1i &  & tx_2+tx_3i\\
\\
x_2-x_3 i &  & x_0-x_1i
\end{array}\right).
\end{equation}

We then identify the basis elements using: 

\begin{equation}
  \label{real-matrix-form}
q_t=x_0\left(\begin{array}{cc}
1 & 0\\
0 & 1
\end{array}\right)+
x_1\left(\begin{array}{cc}
i & 0\\
0 & -i
\end{array}\right)+
x_2\left(\begin{array}{cc}
0 & t\\
1 & 0
\end{array}\right)+
x_3\left(\begin{array}{cc}
0 & ti\\
-i & 0
\end{array}\right),
\end{equation}
and we obtain:
\[
1=\left(\begin{array}{cc}
1 & 0\\
0 & 1
\end{array}\right)=\left(1,0\right)_{t},\;\; i_t=\left(\begin{array}{cc}
i & 0\\
0 & -i
\end{array}\right)=\left(i,0\right)_{t},
\]
and
\[
j_t=\left(\begin{array}{cc}
0 & t\\
1 & 0
\end{array}\right)=\left(0,1\right)_{t},\;\;
k_t=\left(\begin{array}{cc}
0 & ti\\
-i & 0
\end{array}\right)=\left(0,i\right)_{t},
\]
in $\mathcal{H}_{2}^{t}$, which means that every element $q_t$ of $\mathcal{H}_{2}^{t}$ is expressed by
\[
q_t=x_0\left(1,0\right)_{t}+x_1\left(i,0\right)_{t}+x_2\left(0,1\right)_{t}+x_3\left(0,i\right)_{t},
\]
for some $x_l\in\mathbb{R},\,0\le l \le 3$, (as seen in~\cite{alpay_cho_3,MR4591396, MR4592128}).\\

The $\mathbb{R}$-basis $\mathbf{B}_{t}=\left\{ \left(1,0\right)_t,\:\left(i,0\right)_t,\:\left(0,1\right)_t,\:\left(0,i\right)_t\right\} $ satisfies the following properties in the Cayley table of $\mathbb{H}_{t}$, i.e.:

\begin{eqnarray*}
  i_t^{2}&=&\left(-1,0\right)_t=-1, \quad j_t^{2}=t\left(1,0\right)_t=t=k_t^{2},\\
  i_t\cdot_{t}j_t&=&k_t,\quad  j_t\cdot_{t}k_t=-ti_t,\\
  k_t\cdot_{t}i_t&=&j_t,\quad k_t\cdot_{t}j_t=ti\\
  j_t\cdot_{t}i_t&=&-k_t,\quad  i_t\cdot_{t}k_t=-j_t.
 \end{eqnarray*}

\begin{remark}
From now on we write $i$ rather than $i_t$, meaning that we view the ring $\mathbb H_t$
containing a fixed copy of the complex numbers.
\end{remark}

When $t=1$, that is, the split quaternion case, we have
\[
q_{1}^\ct=\begin{pmatrix}\overline{a}&-b\\ -\overline{b}&a\end{pmatrix}\quad{\rm while}\quad
q_{1}^{[*]}=\begin{pmatrix}{a}&\overline{b}\\ {b}&\overline{a}\end{pmatrix}.
\]
With $i_{1},j_{1},k_{1}$ corresponding to $t=1$, and using the notation \eqref{new-def},
one has $q_{1}=a-bj_{1}$
and
\begin{equation}
  q_1^{\ct}=\overline{a}-bj_1\quad{\rm and}\quad q_1^{[*]}=a+\overline{b}j_1.
\end{equation}
The latter case was the one considered in \cite{alss_IJM}.\smallskip

When $t=-1$, that is, for the quaternions, we have
\[
q_{-1}^{\ct}=\begin{pmatrix}\overline{a}&b\\ -\overline{b}&a\end{pmatrix}\quad {\rm and}\quad
q_{-1}^{[*]}=\begin{pmatrix}a&-\overline{b}\\ b&\overline{a}\end{pmatrix}.
\]
Here $q_1^{\ct}$ corresponds to the usual quaternion conjugate, and $q_1^{[*]}$ leads to new directions in analytic quaternionic theory.
With $i,j,k$ corresponding to $t=-1$, one has $q_{-1}=a+bj$ and 
\begin{equation}
  q_{-1}^{\ct}=\overline{a}-bj=\overline{a+bj}=\overline{q_{-1}}
\end{equation}
(since $bj$ has a zero real part) corresponding to classical quaternionic analysis, while
\begin{equation}
  q_{-1}^{[*]}={a}-\overline{b}j={a}+\overline{b}\overline{j}.
\end{equation}
As just mentioned a few lines above, the corresponding theory is new in quaternionic analysis.\\

The paper consists of seven sections and two appendices, of which this introduction (Section~\ref{Intro-Section}) is the first. In Section~\ref{Algebraic-Section+Theory-Functions} we introduce the algebra of hypercomplex numbers, their conjugates, and properties, asa well as an introduction to a theory of functions on these spaces.\\
In Section~\ref{scaled-global} we introduce the scaled global operator $V_t$ in this case and introduced a new set of Fueter-type variables associated to this operator. In the same section (Section~\ref{scaled-global}) we also build parallel theories of regular functions using these new variables (see Definition~\ref{mu-def}) and the previously introduced Fueter-variables associated to the scaled Fueter operator (see Definition~\ref{zeta_t_def}).\\
In Section~\ref{Hardy_spaces_cdast} we introduce and describes Hardy spaces with respect to the conjugate $\circledast$, define the Blaschke product in this case and discuss the interpolation problem in this context, while in Section~\ref{Hardy_q_t} we do the same for the $[*]$ conjugate.

In Section~\ref{Averson_Blaschke} we build the associated Averson space related to the new scaled Fueter variables associated to the scaled global operator, while in Section~\ref{Rational} we establish the premises of a theory of rational functions in this case. Both of these sections are further developed in works that the authors are currently in a preliminary stage of completion.

For clarity, we also write Appendix~\ref{append_krein} and Appendix~\ref{Pontryagin_Q_t}, to remind the reader of the concepts and properties of Krein and Pontryagin spaces, as they apply to our work.

\section{Conjugates, their properties, and a theory of functions}
\setcounter{equation}{0}
\label{Algebraic-Section+Theory-Functions}

In this section, as mentioned in the introduction, one can see that $\mathcal H_2^t$ is not closed under the usual matrix conjugation as in~\eqref{usual*} when $t\not=\pm 1$, unless $t\not =\pm 1$ or $q_t$ is a real scalar matrix.  However, as in~\eqref{q1} and~\eqref{q2}, one can define
two natural adjoints (or conjugations, or conjugate) associated to \eqref{qt890} which leave $\mathcal H_2^t$ invariant, namely
\begin{equation*}
  q_t^{\circledast}=\begin{pmatrix}\overline{a}&-tb\\ -\overline{b}&a\end{pmatrix}\qquad  q_t^{[*]}=\begin{pmatrix}a&t\overline{b}\\ b&\overline{a}\end{pmatrix}.
\end{equation*}
These conjugates lead to different theories, see \cite{alpay_cho_3,MR4591396} and \cite{alpay_cho_6} respectively. On the other hand, just as for the usual $*$-matrix conjugation, we have:
\begin{lemma}
The two new conjugates are contravariant with respect to the matrix product (or, equivalently the product $\cdot_t$ in the representation of the ring in terms of pairs).
\label{contra-variant}
\end{lemma}
\begin{proof}
For the first conjugate we have  with $q_t=\begin{pmatrix}a&tb\\ \overline{b}&\overline{a}\end{pmatrix}$ and $p_t=
  \begin{pmatrix}c&td\\ \overline{d}&\overline{c}\end{pmatrix}$,
  \[
    \begin{split}
      (q_tp_t)^{\ct}&=\left(\begin{pmatrix}a&tb\\ \overline{b}&\overline{a}\end{pmatrix}  \begin{pmatrix}c&td\\ \overline{d}&\overline{c}\end{pmatrix}\right)^{\ct}\\
      &\begin{pmatrix}ac+tb\overline{d}&t(ad+b\overline{c})\\ c\overline{b}+\overline{a}\overline{d}&t\overline{b}+\overline{a}\overline{c}\end{pmatrix}^{\ct}\\
      &=\begin{pmatrix}\overline{a}\overline{c}+t\overline{b}d &-(c\overline{b}+\overline{a}\overline{d})\\
        -(ad+b\overline{c})&ac+tb\overline{d}\end{pmatrix}\\
      &=\begin{pmatrix}\overline{c}&-td\\ -\overline{d}&{c}\end{pmatrix}\begin{pmatrix}\overline{a}&-tb\\ -\overline{b}&{a}\end{pmatrix}\\
      &=p_t^{\ct}q_t^{\ct}.      
      \end{split}
    \]
The second conjugate is considered in Lemma~\ref{[*]adjoint}.
\end{proof}

 \subsection{The regular matrix adjoint $*$}
 \label{regular_*}
$\mathcal H_2^t$ is not closed under the usual matrix conjugation when $t\not=\pm 1$,
\begin{equation*}
  q_t^*=\begin{pmatrix}\overline{a}&b\\ t\overline{b}&a\end{pmatrix}\not\in\mathcal H_2^t
\end{equation*}
if $t\not =\pm 1$ or if $q_t$ is not a real scalar matrix, however this conjugate is useful in defining the matrix Hilbert Schmidt norm and the operator norm. For $q_t$ as above, the Hilbert-Schmidt norm is equal to
\begin{equation}
  \|q_t\|_{HS}=\sqrt{{\rm Tr}\, q_tq_t^*}=\sqrt{2|a|^2+(1+t^2)|b|^2}.
  \label{HS}
  \end{equation}
  
  Let us define the operator norm $\|q_t\|_{op}$ as follows:
\begin{definition}
\label{norm-op-def}
 The {\em operator norm} of $q_t\in \mathcal H_t $ is denoted by $\|q_t\|_{op}$ and it is defined to be the square root of the largest eigenvalue of $q_tq_t^*$, using the regular adjoint $q_t^*$ as in~\eqref{usual*}.
\end{definition}

For ease of calculations, we now compute the operator norm of $q_t$:
\begin{proposition}
  For any  $q_t\in \mathcal H_t $ it holds that 
  \begin{equation}
    \label{norm-op}
\|q_t\|_{op}=\sqrt{\frac{|b|^2(1+t^2)+2|a|^2+\sqrt{|b|^4(1-t^2)^2+4|a|^2|b|^2(1+t)^2}}{2}}
    \end{equation}
  \end{proposition}

  \begin{proof}
    We have
    \[
      \begin{split}
        \det(q_tq_t^*-\lambda I_2)&=\det\begin{pmatrix}|a|^2+t^2|b|^2-\lambda&ab(1+t)\\
          \overline{a}\overline{b}(1+t)&|b|^2+|a|^2-\lambda\end{pmatrix}\\
        &\hspace{-1cm}=\lambda^2-\lambda(|b|^2+|a|^2+|a|^2+t^2|b|^2)+(|b|^2+|a|^2)(|a|^2+t^2|b|^2)-|a|^2|b|^2(1+t)^2\\
      \end{split}
    \]
    and the discriminant, say $\Delta$,  of the corresponding eigenvalue equation is:
    \[
      \begin{split}
        \Delta&=(|b|^2+|a|^2+|a|^2+t^2|b|^2))^2-4((|b|^2+|a|^2)(|a|^2+t^2|b|^2)-|a|^2|b|^2(1+t)^2)\\
        &=(|b|^2+|a|^2-|a|^2-t^2|b|^2))^2+4|a|^2|b|^2(1+t)^2\\
        &=|b|^4(1-t^2)^2+4|a|^2|b|^2(1+t)^2.
    \end{split}
  \]
  Thus
  \[
    \lambda=\frac{|b|^2(1+t^2)+2|a|^2\pm\sqrt{|b|^4(1-t^2)^2+4|a|^2|b|^2(1+t)^2}}{2}
    \]
    with largest eigenvalue corresponding to the $+$ sign, and hence the result.
  \end{proof}

  \begin{corollary}
    \label{opop}
    We have the following:
    \begin{equation}
      \|q_t\|_{op}=      \|q_t^{\ct}\|_{op}=      \|q_t^{[*]}\|_{op}
\end{equation}
and
 \begin{equation}
   \|q_t\|_{HS}=
   \|q_t^{\ct}\|_{HS}=      \|q_t^{[*]}\|_{HS}.
\end{equation}

    \end{corollary}
    \begin{proof}
This is because $\|q_t\|_{op}$ depends only on the absolute values of $a$ and $b$.
      \end{proof}
    
  \begin{remark}
    For $t=-1$ we get $\|q_{-1}\|=\sqrt{|a|^2+|b|^2}$ (as it should be, usually written as $|q|=\sqrt{|a|^2+|b|^2}$) while we get
\[
    \|q_{1}\|=|a|+|b|
\]
    for $t=1$. Note that the latter {\sl is not} a Hilbert space norm.
    \end{remark}
  
We recall that $\mathcal H_2^t$ endowed with the real-valued bilinear form
  \[
\langle p_t,q_t \rangle={\rm Tr} \, q_t^*p_t
    \]
  is a real four dimensional Hilbert space.\smallskip

This norm is cumbersome to deal with, so we define a similar norm that induces the same topology, a scaled    Euclidean-type norm.

\begin{definition}
We will now define the {\em scaled Euclidean norm} of $\mathbb H_t$ to be 
\begin{equation}
\label{scaled-E}
\|q_t\|_{E}=\sqrt{|a|^2+|t|\,|b|^2}=\sqrt{x_0^2+x_1^2+|t|\,(x_2^2+x_3^2)},
\end{equation}
and the scaled inner product associated $\|\cdot \|_{E}$ by $[\cdot,\,\cdot]_E$.
\end{definition}
For $q_t=x_0+x_1i+x_2j_{t}+x_3k_{t}$ and $p_t=y_0+y_1i+y_2j_{t}+y_3k_{t}$ the associated inner product is:
\begin{equation}
\label{ip_E}
[q_t,\,p_t]_E=x_0y_0+x_1y_1+|t|(x_2y_2+x_3y_3).
\end{equation}
\begin{remark}
This norm gives the same topology as the operator norm given in \eqref{norm-op} or the Hilbert-Schmidt norm given in \eqref{HS}.
\end{remark}

We conclude this subsection with a remark on matrices with entries in $\mathcal H_2^t$. Consider Let $G=(g_{uv})_{u,v=1}^N\in(\mathcal H_2^t)^{N\times N}$.
Identifying $\mathcal H_2^t$ with a subset of matrices of $\mathbb C^{2\times 2}$ one can identify $G$ with a matrix in $\mathbb C^{2N\times 2N}$.
Even if this new matrix is invertible in $\mathbb C^{2N\times 2N}$, it is not clear that it will have the structure corresponding to an element
of $(\mathcal H_2^t)^{N\times N}$. First a definition:
      \begin{definition}
        Let $G=(g_{uv})_{u,v=1}^N\in(\mathcal H_2^t)^{N\times N}$.
        We say that $G$ is invertible in $\mathcal H_2^t$ if there exist a matrix
        $H=(h_{uv})_{u,v=1}^N\in(\mathcal H_2^t)^{N\times N}$ such that $GH=HG=I_{(\mathcal H_2^t)^{N\times N}}$.
        The matrices $G_m=(G_{uv})_{u,v=1}^m$ are called the main minors of $G$.
        \end{definition}

      Before we consider this above invertibility question in greater details we mention a number of formulas, often presented for complex
      matrices, but in fact valid for matrices with entries in any algebra, and can be traced at least with the work of Schur
      \cite{schur, schur2}. More recent references include \cite{Dym_CBMS},
       \cite{henderson1981deriving}, \cite[p. 18]{MR1084815}, \cite[p. 113]{frazer}.

\begin{lemma}
Let $(A,B,C,D)\in(\mathcal H_2^t)^{n\times n}\times(\mathcal H_2^t)^{n\times m}\times(\mathcal H_2^t)^{m\times n}\times(\mathcal H_2^t)^{m\times m}$,
and assume that $A$ (resp. $D$) is invertible. Then:
\begin{equation}
\label{wiz}
\begin{pmatrix}A&B\\ C&D\end{pmatrix}=\begin{pmatrix}I_n&0\\ CA^{-1}&I_m\end{pmatrix}\begin{pmatrix}A&0\\ 0&D-CA^{-1}B \end{pmatrix}\begin{pmatrix}I_n&A^{-1}B\\ 0&I_m\end{pmatrix},
\end{equation}
and
\begin{equation}
\label{wiz1}
\begin{pmatrix}A&B\\ C&D\end{pmatrix}=\begin{pmatrix}I_n&BD^{-1}\\ 0&I_m\end{pmatrix}
\begin{pmatrix}A-BD^{-1}C&0\\ 0&D \end{pmatrix}\begin{pmatrix}I_n&0\\ D^{-1}C&I_m\end{pmatrix}
\end{equation}
respectively.
\end{lemma}
As a consequence, with
\begin{equation}
  A^\times=A-BD^{-1}C\quad{\rm and}\quad D^\square=D-CA^{-1}B,
\end{equation}
we have:
\begin{lemma}
Assume $A$ invertible. Then, $M$ is invertible if and only if $D^{\square}$ is
invertible, and it holds that:
\begin{eqnarray}
\nonumber
M^{-1}&=&\begin{pmatrix}I_n&-A^{-1}B\\ 0&I_m\end{pmatrix}\begin{pmatrix}A^{-1}&0\\ 0&(D^\square)^{-1} \end{pmatrix}\begin{pmatrix}I_n&0\\ -CA^{-1}&I_m\end{pmatrix}\\
  \nonumber
  & &\\
\nonumber
& &\\
\label{inv123}
&=&\begin{pmatrix}
 A^{-1}+A^{-1}B(D^\square)^{-1}CA^{-1}&-A^{-1}B(D^\square)^{-1}\\
-(D^\square)^{-1}CA^{-1}&(D^\square)^{-1}\end{pmatrix}.
\end{eqnarray}
Similarly, assuming $D$ invertible, the matrix $M$ is invertible if and only if $A^\times$ is invertible,
and we then have:
\begin{eqnarray}
\nonumber
  M^{-1}&=&\begin{pmatrix}I_n&0\\ -D^{-1}C&I_m\end{pmatrix}
                                            \begin{pmatrix}(A^\times)^{-1}&0\\ 0&D^{-1}
                \end{pmatrix}\begin{pmatrix}I_n&-BD^{-1}\\ 0&I_m\end{pmatrix}\\
\nonumber
& &\\
\nonumber
& &\\
\label{inv1234}
&=&\begin{pmatrix}(A^\times)^{-1}&-({A^\times})^{-1}BD^{-1}\\
-D^{-1}C({A^\times})^{-1}&D^{-1}+D^{-1}C({A^\times})^{-1}BD^{-1}\end{pmatrix}.
\end{eqnarray}
\end{lemma}

        As a prelimary example, consider the case $N=2$, and $G=\begin{pmatrix}g_{11}&g_{12}\\g_{21}&g_{22}\end{pmatrix}$.  Assume that
        $g_{11}$ is invertible, and define
        \[
g_{22}^{\square}=g_{11}-g_{12}g_{22}^{-1}g_{21}.
          \]
          Then it follows from the formula \eqref{inv123} that $G$ is invertible if $g_{22}^\square$ is invertible.
          Using the above formulas one can give {\sl sufficient conditions} in terms of $N$ element of $\mathcal H_2^t$
          for the matrix $G$ to be invertible in $(\mathcal H_2^t)^{N\times N}$.

\subsection{The adjoint ${\circledast}$}
\label{cdast}

As written in the introduction we define the new adjoint ${\circledast}$ as follows:
\begin{definition}
For $q_t=\begin{pmatrix}a&tb\\ \overline{b}&\overline{a}\end{pmatrix}, \,q_t\in\mathcal H_t$ we define 
$$q_t^{\circledast}=\begin{pmatrix}\overline{a}&-tb\\ -\overline{b}&a\end{pmatrix}.$$
\end{definition}
The bilinear real form associated to ${\circledast}$ is
\begin{equation}
  \label{form1}
  [q_t,p_t]_{\circledast}=a\overline{c}+\overline{a}c-t(b\overline{d}+\overline{b}d)
\end{equation}
with $p_t=\begin{pmatrix}c&td\\ \overline{d}&\overline{c}\end{pmatrix}$.

\begin{lemma}
The $\circledast$-adjoint satisfies
  \begin{equation}\label{qqast}
q_t  q_t^{\circledast}=q_t^{\circledast}q_t=(\det q_t) I_2=(|a|^2-t|b|^2)I_2.
\end{equation}
\end{lemma}
As we saw above, the associated bilinear real form is
\begin{equation}
  \label{form1}
  [q_t,p_t]_{\circledast}={\rm Tr}\, q_t^{\ct}p_t=a\overline{c}+\overline{a}c-t(b\overline{d}+\overline{b}d).
\end{equation}

This $\ct$ adjoint leads to a theory which is completely parallel to the study of slice hyperholomorphic functions.

We set
\[
\frac{q_t+q_t^{\circledast}}{2}={\rm Re}\, q_t
\]
We first note that
\begin{equation}
{\rm Re}\, q_t= ({\rm Re}\, a)I_2=x_0I_2,
  \end{equation}
where we use the notation of equation~\eqref{real-matrix-form}
\begin{lemma}
It holds that $q_t^\circledast=q_t$ if and only if $q$ is a scalar matrix, $q_t=xI_2$ for some $x\in\mathbb R$.
\end{lemma}

\begin{proof}
This is a direct computation.
\end{proof}

Formula \eqref{qqast} shows that $q_tq_t^{\circledast}$ need not be positive for $t>0$. More precisely:

\begin{theorem}
  \label{pont-ht}
 Let $t\in\mathbb R^*$. The space $\mathcal H_2^t$ endowed with the form
\begin{equation}
  \label{form-qsrar}
[p_t,q_t]_{\circledast}={\rm Tr}~(q_t^{\circledast}p_t)={\rm Tr}~(p_t^{\circledast}q_t)
\end{equation}
is a real Hilbert space of dimension four for $t<0$ and a real Pontryagin space of dimension four, with (real) index of negativity equal to two, for $t>0$.
\end{theorem}

\begin{proof}
With $p_t=\begin{pmatrix}c&td\\\overline{d}&\overline{c}\end{pmatrix}$ we have
\begin{equation}
  \label{tyu-12345}
{\rm Tr}~(q_t^{\circledast}p_t)={\rm Tr}~(p_t^{\circledast}q_t)=a\overline{c}+\overline{a}c-t(\overline{b}d+b\overline{d}).
\end{equation}
Furthermore,  

\begin{equation*}
  \begin{split}
    [p_t,q_t]_{\circledast}&={\rm Tr}~(q_t^{\circledast}p_t)\\
    &=a\overline{c}+\overline{a}c-t(\overline{b}d+b\overline{d})\\
    &=\frac{1}{2}\left\{(a+\overline{a})(c+\overline{c})-(a-\overline{a})(c-\overline{c})-t\left((b+\overline{b})(d+\overline{d})-(b-\overline{b})(d-\overline{d})\right)\right\}\\
    &=\frac{1}{2}\left\{(a+\overline{a})(c+\overline{c})+(a-\overline{a})\overline{(c-\overline{c})}-t\left((b+\overline{b})(d+\overline{d})+(b-\overline{b})\overline{(d-\overline{d})}\right)\right\},
  \end{split}
  \end{equation*}
  which shows that in the case of $t>0$ the form has two positive squares and two negative squares and $\mathcal H_2^t$ becomes a Pontryagin space with $ [p_t,q_t]_{\circledast}$ and it is a Hilbert space for $t<0$.
which ends the proof (see also \cite{alpay_cho_3} and \cite{alpay_cho_4}).
  \end{proof}

As a corollary we have the following result.

\begin{corollary}
  \label{help-facto} Let $\mathcal H^t_2$ be the scaled hypercomplex space for a fixed $t$, where we assume that $t\not=0$. If, moreover, for a fixed $p_t\in\mathcal H^t_2$ we have ${\rm Tr}\, q_t^{\ct}p_t=0$ for all $q_t\in\mathcal H^t_2$, then $p_t=0$.
\end{corollary}

\begin{proof}
  Let
  \[
  p_t=\begin{pmatrix}c&td\\ \overline{d}&\overline{c}\end{pmatrix}\quad{\rm and}\quad
  q_t=\begin{pmatrix}a&tb\\ \overline{b}&\overline{a}\end{pmatrix}.
  \]
  Then, as in~\eqref{form1},
  \[
{\rm Tr}\,q_t^{\ct}p_t=a\overline{c}+\overline{a}c-t(b\overline{d}+\overline{b}d).
\]
Assume ${\rm Tr}\,q_t^{\ct}p_t=0$ and first take $b=0$. The choice $a=1$ and $a=i$ lead respectively to ${\rm Re}\,c=0$ and ${\rm Im}\, c=0$ so that $c=0$. The case of $d$ is treated similarly.
  \end{proof}
  
The first adjoint $\ct$ leads to a theory which is completely parallel to the study of slice hyperholomorphic functions and the authors are preparing another work in this direction as well.

In order to properly define spaces of functions, we will write the connection between the scaled norm and the $\circledast$ inner product as follows. Let $J_{\circledast}$ be the matrix:

\begin{equation}
\label{J-circle}
J_{\circledast}=\left(\begin{array}{cccc}
1 & 0 & 0 & 0\\
0 & 1 & 0 & 0\\
0 & 0 & -\frac{t}{|t|} & 0\\
0 & 0 & 0 & -\frac{t}{|t|}
\end{array}\right)_{\mathbb R}=\left(\begin{array}{cc}
1 & 0 \\
0 & -\frac{t}{|t|}
\end{array}\right)_{\mathbb C}
\end{equation}
$J_{\ct}$ is a signature matrix (both self-adjoint and equal to its inverse) and induces a signature operator (which we still denote by $J_{\ct}$ on $\mathbb H_t$ such that 
  \begin{equation}
    [q_t,q_t]_{\ct}=[q_t,J_{\circledast}q_t]_E
  \end{equation}
  where $[\cdot,\cdot]_E$ is the inner product associated to the norm $\|\cdot\|_E$, as in~\eqref{ip_E}.\smallskip

   Note that, using the representation \eqref{new-def} we have
  \begin{equation}
    J_{\ct}q_t=\overline{a}-bj_t.
    \end{equation}

    \subsection{The adjoint $[*]$.}
    \label{bracket_*}
  
  We now turn to the definition of adjoint $q_t^{[*]}$, as in the introduction.
  
 \begin{definition} For a fixed scalar $t\in\mathbb{R}$, and $q_t=\left(a,b\right)_t \in\mathbb H_t$  we define $q_t^{[*]}$ by
\begin{equation}
\left(a,b\right)_t^{[*]}\overset{\mathrm{def}}{=}\left(a,\overline{b}\right)_t\;\;\mathrm{in\;\;}\mathbb{H}_{t}.
\label{second-adjoint}
\end{equation}
\end{definition}
Equivalently, if we take $q_t$ to be an element of $\mathcal{H}_{2}^{t}$ rather than $\mathbb{H}_{t}$ then $q_t^{[*]}$ becomes:
\begin{equation}
\left(\begin{array}{cc}
a & tb\\
\overline{b} & \overline{a}
\end{array}\right)^{[*]}=\left(\begin{array}{cc}
a & t\overline{b}\\
b & \overline{a}
\end{array}\right)
\end{equation}
on $\mathcal{H}_{2}^{t}$.

  The second adjoint $[*]$ satisfies the following identity in matrix form:
  \begin{equation}
    \label{J*J}
    q_t^{[*]}=Jq_t^*J\quad{\rm with }\quad J=\begin{pmatrix}0&1\\1&0\end{pmatrix}.
  \end{equation}
   Note that, even if, $q_t^*$ and $J$ are not in $\mathcal H_t$ for general t, we have that $q_t^{[*]}\in\mathcal H_t$
and hence it becomes a Pontryagin space adjoint on $\mathbb C^{2\times 2}$ (see Section \ref{append_krein} and Example \ref{Cnpont}; see also \cite{alpay_cho_3}).

We also note that, in general, $q_t$ is not normal for this adjoint,
  (i.e. it does not commute with $q_t^{[*]}$):
  \begin{equation}
    \label{notnormal}
    q_tq_t^{[*]}\not=q_t^{[*]}q_t.
  \end{equation}

  The adjoint $[*]$ is associated to the bilinear real form
\begin{equation}
  \label{form2}
  [q_t,p_t]_{[*]}=a{c}+\overline{ac}-t(b{d}+\overline{bd})
\end{equation}
(as compared with \eqref{form1}), and leads to a new theory, more ``non-commutative''. A way to see this is to compare the definitions
of the corresponding Blaschke factors later on in Section~\ref{Averson_Blaschke} (see formulas \eqref{bl-1} and \eqref{bl-2} respectively).

The second adjoint $[*]$ satisfies
  \begin{equation}
    \label{J*J}
    q_t^{[*]}=Jq_t^*J\quad{\rm with }\quad J=\begin{pmatrix}0&1\\1&0\end{pmatrix}.
  \end{equation}

    \begin{lemma}
    \label{[*]adjoint}
      It holds that, for every $p_t,q_t\in\mathbb H_t$
      \begin{eqnarray}
        (q_t^{[*]})^{[*]}&=&q_t\\
        (p_t\cdot_t q_t)^{[*]}&=&q_t^{[*]}\cdot_t p_t^{[*]}
        \end{eqnarray}
      \end{lemma}
\begin{proof}
      The proof follows from \eqref{J*J} and of the corresponding properties for the regular adjoint.
      We leave the details to the reader.
      \end{proof}
      We now make more precise the fact (see~\eqref{notnormal} above) that $q_t$ is not normal in general with respect to this adoint.

      \begin{proposition}
        Let $q_t=\begin{pmatrix}a&tb\\ \overline{b}&\overline{a}\end{pmatrix}\in\mathcal  H_2^t$. Then,
        \begin{equation}
          \begin{split}
            q_tq_t^{[*]}&=\begin{pmatrix} a^{2}+tb^{2}& 2{\rm Re}\,~a\overline{b}\\
2{\rm Re}~a\overline{b}&\overline{a}^{2}+t\overline{b}^{2}\end{pmatrix}\\
q_t^{[*]}q_t&=\begin{pmatrix} a^{2}+t\overline{b}^{2}& 2t{\rm Re}~ab\\
2 {\rm Re}~ab&\overline{a}^2+tb^2\end{pmatrix}.
 \end{split}
 \label{3.4}
 \end{equation}
\end{proposition}

\begin{proof}
The formulas \eqref{3.4} are shown by straightforward computations. See \cite{alpay_cho_6} for details.
\end{proof}

   \begin{lemma}
The form
      \[
   [q_t,p_t]_{[*]}=   {\rm Tr}\left( p_t^{[*]}q_t\right),\quad p,q\in\mathcal H_2^t
      \]
      is bilinear on $\mathcal H_2^t$ with respect to the real numbers, satisfies
      \[
        [q_t,p_t]_{[*]}=[p_t,q_t]_{[*]}
        \]
      and is not degenerate on $\mathcal H_2^t$.
      \end{lemma}
      \begin{proof}
        Set
        \[
          q_t=\begin{pmatrix}a&tb\\ \overline{b}&\overline{a}\end{pmatrix}\quad{\rm and}\quad          p_t=\begin{pmatrix}c&td\\ \overline{d}&\overline{c}\end{pmatrix},
          \]
and assume that $q\in\mathbb H_t$ is such that $[q_t,p_t]_{[*]}=0$ for all $p\in\mathbb H_t$. 
Since
\[
  q_t^{[*]}p_t=\begin{pmatrix}  ac+t\overline{b}\overline{d}&t(ad+\overline{b}\overline{c})\\
bc+\overline{a}\overline{d}&tbd+\overline{a}\overline{c}\end{pmatrix}
\]
we have
\[
  [q_t,p_t]_{[*]}={\rm Re}\,(ac+t\overline{b}\overline{d}).
\]
The choice $d=0$ and $c\in\mathbb R$ gives ${\rm Re}\,a=0$ and the choice $c\in\mathbb R$ gives ${\rm Im}\, a=0$. So $a=0$.
To see that $b=0$ take then $d$ real and then $d$ purely imaginary.
\end{proof}

By the above analysis, one can get the following result (see \cite{alpay_cho_6}).

\begin{theorem}
\label{bracket-Pontryagin}
If $t\not=0$, the pair $\left(\mathcal{H}^t_{2},\left[,\right]_{[*]}\right)$ is a Pontryagin space over $\mathbb{R}$.
\end{theorem}
\begin{proof}
We have that:
\begin{eqnarray*}[q_t,p_t]_{[*]} &=&2\, {\rm Re}\,(ac+t\overline{b}\overline{d})\\
&=&ac+\overline{a}\overline{c}+t(bd+\overline{b}\overline{d})\\
&=& \frac{1}{2}\left\{(a+\overline{a})(c+\overline{c})-(a-\overline{a})\overline{(c-\overline{c})}+t\left((b+\overline{b})(d+\overline{d})-(b-\overline{b})\overline{(d-\overline{d})}\right)\right\},
\end{eqnarray*}
which shows that the form has two positive squares and two negative squares and the result follows.
\end{proof}

We will now write the connection between the scaled norm and the $[*]$ inner product as follows. Let $J_{[*]}$ be the matrix:

\begin{equation}
\label{J-[*]}
J_{[*]}=\left(\begin{array}{cccc}
1 & 0 & 0 & 0\\
0 & -1 & 0 & 0\\
0 & 0 & \frac{t}{|t|} & 0\\
0 & 0 & 0 & -\frac{t}{|t|}
\end{array}\right)_{\mathbb R}
\end{equation}

$J_{[*]}$ is a signature matrix (both self-adjoint and equal to its inverse) and induces a signature operator (which we still denote by $J_{[*]}$ on $\mathbb H_t$ such that 
such that
 \begin{equation}
    [q_t,q_t]_{[*]}=[q_t,J_{[*]}q_t]_E.
  \end{equation}
  Note that, using the representation $q_t=a+b j_t$ we have
  \begin{equation}
    J_{[*]}q_t={a}+\overline{b}j_t.
    \end{equation}

\begin{remark}
As a reminder, we now have $\left(\mathcal{H}^t_{2},\left[,\right]_{[*]}\right)$ a Pontryagin space for all $t\neq 0$, while  $\left(\mathcal{H}^t_{2},\left[,\right]_{\circledast}\right)$ is a Pontryagin space for $t>0$ and a Hilbert space for $t<0$.
\end{remark}
      \begin{remark}
        Let $t=\pm 1$. It holds that $q_1^*=q_1^{[*]}$ if and only if $b=\overline{b}$ and $a=\overline{a}$, i.e. if and only if $q_1$ is an hyperbolic number.
        Similarly, let $t=-1$. Then $q_1^*=q_1^{[*]}$ if and only if $b=-\overline{b}$ and $a=\overline{a}$.
      \end{remark}


In a table form, the adjoints of $i$, $j_t$ and $k_t$ with respect to $\ct$ and $[*]$ are:
\begin{center}
      \begin{tabular}{|c|c|c|}
    \hline
    &$\ct$ &$[*]$\\
    \hline
  $1$&$1$&$1$\\
    \hline
    $i$ & $-i$ & $i$\\
    \hline
    $j_t$&$-j_t$&$j_t$\\
    \hline
    $k_t$&$-k_t$&$-k_t$\\
    \hline
  \end{tabular}
\end{center}
or, equivalently, in terms of matrices:
\begin{center}
  \begin{tabular}{|c|c|c|}
    \hline
    &$\ct$ &$[*]$\\
    \hline
  $\begin{pmatrix}1&0\\0&1\end{pmatrix}$&$\begin{pmatrix}1&0\\0&1\end{pmatrix}$&$\begin{pmatrix}1&0\\0&1\end{pmatrix}$\\
    \hline
    $\begin{pmatrix}i&0\\0&-i\end{pmatrix}$&$-\begin{pmatrix}i&0\\0&-i\end{pmatrix}$&$\begin{pmatrix}i&0\\0&-i\end{pmatrix}$\\
    \hline
    $\begin{pmatrix}0&t\\1&0\end{pmatrix}$&$-\begin{pmatrix}0&t\\1&0\end{pmatrix}$&$\begin{pmatrix}0&t\\1&0\end{pmatrix}$\\
    \hline
    $\begin{pmatrix}0&ti\\-i&0\end{pmatrix}$&$-\begin{pmatrix}0&ti\\-i&0\end{pmatrix}$&$-\begin{pmatrix}0&ti\\-i&0\end{pmatrix}$\\
    \hline 
  \end{tabular}.
\end{center}
  
  Throughout the paper we will make use of the three equivalent definitions of scaled hypercomplex numbers and their conjugates and use either the matrix, vector, or complex pair notation as needed. As a reminder, we include them here:
\begin{equation*}
q_t^*=\begin{pmatrix}\overline{a}&b\\ t\overline{b}&a\end{pmatrix} \qquad  q_t^{\circledast}=\begin{pmatrix}\overline{a}&-tb\\ -\overline{b}&a\end{pmatrix}\qquad  q_t^{[*]}=\begin{pmatrix}a&t\overline{b}\\ b&\overline{a}\end{pmatrix}.
\end{equation*}

  These conjugates will lead interesting function theories and it is worth mentioning here the example of $f(q_t)=\sqrt{1-q_t}$ which will have a familiar corresponding power series decomposition depending on the two adjoints $\circledast$ and $[*]$.

 In this paper, we do not develop the degenerate case $t=0$, as it is not relevant here. This case is currently being worked on in another project. 
  \subsection{Spaces of functions on $\mathcal H_2^t$}
\label{Function_spaces}

We start with considering the spectrum of a matrix given by the $\circledast$ adjoint. This is the only case that can be considered, since, for the $[*]$ adjoint we have that $q_t$ is not normal for this adjoint (see~\eqref{notnormal}) and that $q_t q_t^{[*]}$ is not scalar.

\begin{proposition}
 For $q_t\in\mathbb H_t$ and $A\in\mathbb H_t^{m\times m}$ we have that:
  \begin{equation}
\label{seconde-equation}
    \sum_{n=0}^\infty q_t^nA^n=(1-q_t^{\circledast}A)(q_tq_t^{\circledast}A^2-2({\rm Re}\, (q_t) )A+I_2)^{-1},
      \end{equation}
whenever the sum on the left is convergent.
\end{proposition}

       \begin{proof}
         We prove \eqref{seconde-equation};  We have
         \[
           \begin{split}
             \left(    \sum_{n=0}^\infty q_t^nA^n\right)\left(q_tq_t^{\circledast}A^2-2({\rm Re}\, (q_t^{\circledast}) )A+I_2\right)&=\\
             &\hspace{-3cm}=\sum_{n=0}^\infty q_t^n(q_tq_t^{\circledast})A^{n+2}-2\sum_{0}^\infty q_t^n ({\rm Re}\, (q_t)) A^{n+1}+ \sum_{n=0}^\infty q_t^nA^n\\
             &\hspace{-3cm}=\sum_{n=0}^\infty q_t^{n+1}q_t^{\circledast}A^{n+2}-\sum_{0}^\infty q_t^n (q_t+q_t^{\circledast}) A^{n+1}+ \sum_{n=0}^\infty q_t^nA^n\\
             &\hspace{-3cm}=\sum_{n=0}^\infty q_t^{n+1}q_t^{\circledast}A^{n+2}-\sum_{n=0}^\infty q_t^{n+1}A^{n+1}-\sum_{n=0}^\infty q_t^nq_t^{\circledast}A^{n+1}+ \sum_{n=0}^\infty q_t^nA^n\\
             &\hspace{-3cm}=\underbrace{\sum_{n=0}^\infty q_t^{n+1}q_t^{\circledast}A^{n+2}-\sum_{n=0}^\infty q_t^nq_t^{\circledast}A^{n+1}}_{-q_t^{\circledast}A}            \underbrace{+ \sum_{n=0}^\infty q_t^nA^n             -\sum_{n=0}^\infty q_t^{n+1}A^{n+1}}_{I_2}\\
             &\hspace{-3cm}=I_2-q_t^{\circledast}A.
           \end{split}
           \]
  \end{proof}


\begin{corollary} For any $p_t,q_t\in\mathbb H_2^t$ we have that:
  \begin{equation}
    \label{seconde-equation-1}
    \sum_{n=0}^\infty q_t^np_t^ n=(1-q_t^{\circledast}p_t)(q_tq_t^{\circledast}p_t^2-2({\rm Re}\, (q_t) )p_t+1)^{-1},
  \end{equation}
  whenever the sum on the left is convergent.
  \end{corollary}
\begin{proof}
Setting then $A=p_t$ in \eqref{seconde-equation} gives \eqref{seconde-equation-1}.
  \end{proof}

\begin{remark}
In the case where $t=1$, the quaternionic case, this spectrum first appeared in \cite{MR2496568,MR2661152,MR2752913}.
\end{remark}
We now introduce the Cauchy product of power series in $\mathcal H_2^t$.
In view of the non-commutativity there are various possible ways to define power series and their product; here we follow \cite{MR51:583}, and consider
power series of the form
\begin{equation}
\label{pwerseriesqt}
  f(q_t)=\sum_{n=0}^\infty q_t^n\alpha_n,
\end{equation}
where both the variable $q_t$ and the coefficients $\alpha_0,\alpha_1,\ldots$ are in $\mathcal H_2^t$.

\begin{definition}
 Let $f(q_t)=\sum_{n=0}^\infty q_t^n\alpha_n$ and $g(q_t)=\sum_{n=0}^\infty q_t^n\beta_n$, then the Cauchy-star product of $f$ and $g$ is:
 \begin{equation}
 \label{star-f-g}
 (f\star g)(q_t)=\sum_{n=0}^{\infty}\sum_{k=0}^n q_t^{n-k}\alpha_{n-k}\star\, q_t^k\beta_k
 \end{equation}
where:
\[
  q_t^l\alpha\star_tq_t^m\beta=q_t^{l+m}\alpha\cdot_t \beta, \quad l,m\in\mathbb N_0,\quad \alpha,\beta\in\mathcal H_2^t,
\]
is called the star-product.
\end{definition}

Throughout the paper we will call both of the above $\star$ products.
The $\star$-product does not respect point-evaluation, but as in the quaternionic setting we have:

\begin{lemma}
  \label{multi-3}
  In the above notation assume that $f(q_t)\not=0$. We then have:
  \begin{equation}
    (f\star g)(q_t)=f(q_t)g\left(f(q_t)^{-1}q_tf(q_t)\right).
    \label{multi}
    \end{equation}
  \end{lemma}

  \begin{proof}
    In the above notation we have:
    \[
      \begin{split}
        (f\star g)(q_t)&=\sum_{n=0}^\infty q_t^nf(q_t)\beta_n\\
        &=\sum_{n=0}^\infty f(q_t)f(q_t)^{-1}q_t^nf(q_t)\beta_n\\
        &=\sum_{n=0}^\infty f(q_t)\left(f(q_t)^{-1}q_tf(q_t)\right)^n\beta_n\\
        &=f(q_t)\left(\sum_{n=0}^\infty f(q_t)\left(f(q_t)^{-1}q_tf(q_t)\right)^n\beta_n\right)\\
&=        f(q_t)g\left(f(q_t)^{-1}q_tf(q_t)\right).
      \end{split}
      \]
    \end{proof}

The proof of the following important lemma is easy and will be omitted.

\begin{lemma}
  \label{extens}
  The power series \eqref{pwerseriesqt} is uniquely determined by the values $f(x)$, where $x$ is real and in an open interval of the form $(-\e,\e)$, for some $\e>0$.
  Furthermore the restriction to the real line of the $\star$-product reduces to the regular (non-commutative) product of power series of a real variable with matrix coefficients.
\end{lemma}

Thus these are series of matrices (with a certain structure), and the results
from series of matrices will hold, using the operator norm, and the claims can be proved by first considering the case $q_t\in\mathbb R$. In particular:

\begin{proposition}  
  \label{unique-123}
  Let $f(q_t)=\sum_{n=0}^\infty q_t^n\alpha_n$ and $g(q_t)=\sum_{n=0}^\infty q_t^n\beta_n$ be power series of the form
  \eqref{pwerseriesqt}, converging in a neighborhood $\|q_t\|_{op}<\epsilon$ of the zero
  matrix in $\mathcal H_2^t$. The star product $f\star g$ converges in a (possibly smaller) neighborhood of the origin.
  \end{proposition}

\begin{proof}
  Let $(\gamma_n)$ be the corresponding convolution sequence:
  \[
\gamma_n=\sum_{k=0}^n\alpha_k\beta_{n-k},\quad n=0,1,\ldots
\]
Restricting $q_t=x_0I_2$ where $x\in\mathbb R$, and resorting to the classical theory, we see that the series $\sum_{n=0}^\infty x_0^n\gamma_n$
converges for $|x_0|<\eta$ for some $\eta>0$. Thus
\[
\frac{1}{R}=\limsup_{n\rightarrow \infty}\|\gamma_n\|_{op}^{\frac{1}{n}}>0.
  \]
  For $\|q_t\|<R$ the series $\sum_{n=0}^\infty \|q_t\|_{op}^n\cdot\|\gamma_n\|_{op}$ converges and so does $f\star g$.
  \end{proof}

\begin{proposition}
  Assume $\alpha_0=1$, and write $f(q_t)=1-g(q_t)$ with $g(q_t)=\sum_{k=1}^\infty q_t^k\alpha_k$. Then, the inverse $f^{-\star}(q_t)$ is equal to a convergent power series of the form
  \eqref{pwerseriesqt} in a neighborhood  of the origin and moreover
  \[
f^{-\star}(q_t)=\sum_{n=0}^\infty g^{\star n}(q_t)
\]
where $g^{\star n}=g\star g\dots \star g$, $n-$times, and the convergence is in the topology of $\mathbb C^{2\times 2}$.
  \end{proposition}

\begin{proof}
  The result is true for $q_t=x\in\mathbb R$ and then follows since the restriction to any subinterval of the real line uniquely determines a power series of the form
  \eqref{pwerseriesqt} (see Lemma \ref{extens}).
  \end{proof}

  
The following lemma will play a key role in the study of rational functions.

\begin{lemma}
  \label{real-conv}
Let $\alpha=(\alpha_n)_{n\in\mathbb N_0}$ be a sequence of elements in $\mathcal H_2^t$ and let $\alpha^{\circledast}=(\alpha_n^{\circledast})_{n\in\mathbb N_0}$ denote the corresponding sequence of $\circledast$-adjoints. Then the sequence of self-convolution $\alpha\otimes \alpha^{\circledast}$ is real.
\end{lemma}

\begin{proof}
  We have
  \[
(  \alpha\otimes \alpha^{\circledast})_n  =\sum_{\substack{\ell,k\in\left\{0,\ldots, n\right\}\\ k+\ell=n}} \alpha_k\alpha_\ell^{\circledast}
\]
which is equal to its own $\circledast$-adjoint since
\begin{equation}
  \label{qqast1}
  \sum_{\substack{\ell,k\in\left\{0,\ldots, n\right\}\\ k+\ell=n}} \alpha_k\alpha_\ell^{\circledast}=\sum_{\substack{\ell,k\in\left\{0,\ldots, n\right\}\\ k+\ell=n}} \alpha_\ell\alpha_k^{\circledast}
  \end{equation}
\end{proof}

\begin{definition}
\label{circle-conj-f}
  Let $f(q_t)=\sum_{n=0}^\infty q_t^n\alpha_n$ be converging in a neighborhood $\|q_t\|<\epsilon$ of the zero element of $\mathcal H_2^t$. We define
  \begin{equation}
    f^{\circledast}(q_t)=\sum_{n=0}^\infty q_t^n\alpha_n^{\circledast}.
    \end{equation}
  \end{definition}

  We note that, for $q_t=x\in\mathbb R$ for which the power series converges,
  \[
    f^{\circledast}(x)=(f(x))^{\circledast}.
    \]

    \begin{proposition}
The power series $f(q_t)\star f^{\circledast}(q_t)$ has real coefficients and commute with every other converging power series.
      \end{proposition}

      \begin{proof}
This follows from Lemma  \ref{real-conv}.
        \end{proof}
      
  \begin{proposition}
    Assume that  $\alpha_0$ in~\eqref{pwerseriesqt} is invertible. Then,
    \begin{equation}
            (f(q_t))^{-\star}=f^{\circledast}(q_t)(f(q_t)\star f^{\circledast}(q_t))^{-1}= (f(q_t)\star f^{\circledast}(q_t))^{-1} f^{\circledast}(q_t).
    \end{equation}
    \end{proposition}

    \begin{proof} We observe that, since $f(q_t)\star f^{\circledast}(q_t)$ is real valued and $\star$ is associative:
      \[
        (f(q_t))^{-\star} \, (f(q_t)\star f^{\circledast}(q_t))=  (f(q_t))^{-\star} \star (f(q_t)\star f^{\circledast}(q_t))=
        (f(q_t)^{-\star}\star f((q_t)))\star f^{\circledast}(q_t)=f^{\circledast}(q_t).
        \]
      \end{proof}
In preparation of the study of rational functions we specialize the previous result to polynomial
    
\begin{corollary}
  \label{coro-2-19}
  Let $P(x)=\sum_{n=0}^Nx^n\alpha_n$ be a polynomial of the real variable $x$, with coefficients $\alpha_0,\ldots, \alpha_N \in \mathcal H_2^t$.
  Then the polynomial
  \begin{equation}
    \label{ppstar}
P(x)(P(x))^{\circledast}
\end{equation}
has real coefficients and
\begin{equation}
(P(q_t))^{-\star}=P^{\circledast}(q_t)(P(q_t)\star P^{\circledast}(q_t))^{-1}= (P(q_t)\star P^{\circledast}(q_t))^{-1} P^{\circledast}(q_t).
\label{quotient-123}
\end{equation}
\end{corollary}

We note that \eqref{quotient-123} can be rewritten as
\[
  (P(q_t))^{-\star}=\frac{P^{\circledast}(q_t)}{(P(q_t)\star P^{\circledast}(q_t))}.
\]
We present an illustration of the previous corollary, which pertains to Blaschke factors; see Section \ref{bl1} for the latter.

\begin{example}
  Let $\alpha\in\mathcal H_2^t$. Then,
  \[
    \begin{split}
      (1-q_t\alpha)^{-\star}&=      (1-q_t\alpha^{\ct})((1-q_t\alpha)\star(1-q_t\alpha^{\ct}))^{-1}\\
      &=\frac{1-q_t\alpha^{\ct}}{q_t^2(\alpha\alpha^{\ct})-2q_t{\rm Re}~\alpha+1}
    \end{split}
    \]
  \end{example}

We conclude this section with:
\begin{proposition}
\label{B-star-comm}
 We have that $q_t-\alpha$ and $(1-q_t\alpha)$ $\star$-commute. It then follows that $q_t-\alpha$ and $(1-q_t\alpha)^{-\star}$  also $\star$-commute.
 \end{proposition} 

 \begin{proof}
   The claim amounts to show that
   \[
     (1-q_t\alpha)\star(q_t-\alpha)=     (q_t-\alpha)\star(1-q_t\alpha),
   \]
   but both sides are equal to $-q_t^2\alpha+q_t(\alpha^2+1)-\alpha$.   
   
   The second part follows.
 \end{proof}

 In fact the above proposition is a special case of the following result, whose proof is left to the reader.

 \begin{proposition}
  Let $f(q_t)=\sum_{n=0}^\infty q_t^n\alpha_n$ and $g(q_t)=\sum_{n=0}^\infty q_t^n\beta_n$ be power series of the form
  \eqref{pwerseriesqt}, converging in a neighborhood $\|q_t\|_{op}<\epsilon$ of zero. 
  If their coefficients pairwise commute, i.e. $\alpha_n\beta_m=\beta_m\alpha_n,\,\,\forall n,m\in\mathbb N_0,$
then $f$ and $g$ will also $\star$-commute.
   \end{proposition}
   
   We conclude with a definition generalizing the quaternioinc setting:

   \begin{definition}
   \label{sph}
     Let $\alpha\in\mathbb H_t$. We denote by $[\alpha]$, and call sphere associated to $\alpha$, the set of points $\beta\in\mathbb H_t$ such that
\begin{eqnarray}
  {\rm Re}\,\alpha&=&{\rm Re}\, \beta\\
  \alpha\alpha^{\ct}&=&  \beta\beta^{\ct}
                       \end{eqnarray}
     \end{definition}


\section{The Scaled Global Operator for Scaled Hypercomplex Numbers}
\setcounter{equation}{0}
\label{scaled-global}

For the scaled variable $q_t$ we have its vector part:
\[
  \vec{q_t}=ix_1+j_t x_2+k_t x_3,
  \]
  and, similar to the global operator $G_q$ on the space of quaternions, we can define the global operator associated to the scaled quaternions $\mathbb H_t$ to be:

\begin{equation}
  G_t=\vec{q_t}\frac{\partial}{\partial x_0}-\left(x_1\frac{\partial}{\partial x_1}+x_2\frac{\partial}{\partial x_2}+x_3\frac{\partial}{\partial x_3}\right)
\end{equation}
and its associated, simplified version:
\begin{equation}
  V_t=\frac{\partial}{\partial x_0}-\frac{1}{ix_1+j_tx_2+k_t x_3}\left(x_1\frac{\partial}{\partial x_1}+x_2\frac{\partial}{\partial x_2}+x_3\frac{\partial}{\partial x_3}\right)
\end{equation}

\subsection{The Gleason problem associated to the $V_{t}$- operator and the $V_t$-Fueter variables}
\setcounter{equation}{0}

As in~\cite{adv_prim} we solve the Gleason problem associated to the new $V_{t}$- operator and define $\mu_{t,l}$, the Fueter-type variables which are in its kernel and which generate the series solutions to $V_t f=0$, where $f$ is real-analytic.\\
Assume $f$ smooth and in the kernel of $V_t$, which yields:
\[
  \frac{\partial f}{\partial x_0}=\frac{1}{ix_1+j_t x_2+k_tx_3}\left(x_1\frac{\partial f}{\partial x_1}+x_2\frac{\partial f}{\partial x_2}+x_3\frac{\partial f}{\partial x_3}\right).
\]
  
Varying $x$ by a real parameter $u$, we obtain:
\begin{equation}
  \begin{split}
    \frac{{\rm d}f(ux)}{{\rm d}u}&=x_0\frac{\partial f}{\partial x_0}(ux)+x_1\frac{\partial f}{\partial x_1}(ux)+x_2\frac{\partial f }{\partial x_2}(ux)+x_3\frac{\partial f}{\partial x_3}(ux)\\
    &=\frac{x_0}{\vec{q_t}}\left(x_1\frac{\partial f}{\partial x_1}(ux)+x_2\frac{\partial f}{\partial x_2}(ux)+x_3\frac{\partial f}{\partial x_3}(ux)\right)+\\
    &\hspace{5mm}+x_1\frac{\partial f}{\partial x_1}(ux)+x_2\frac{\partial f }{\partial x_2}(ux)+x_3\frac{\partial f}{\partial x_3}(ux)\\
    &=\mu_{1,t}(x)\frac{\partial f}{\partial x_1}(ux)+\mu_{2,t}(x)\frac{\partial f}{\partial x_2}(ux)+\mu_{3,t}(x)\frac{\partial f}{\partial x_3}(ux)
  \end{split}
\end{equation}
where:
\begin{equation}
\label{mu_t}
\mu_{l,t}(x)=x_l\left(1+\frac{x_0}{\vec{q_t}}\right),\quad l=1,2,3.
\end{equation}

\begin{definition}
\label{mu-def}
The functions $\mu_{l,t}$, $l=1,2,3$, are called the {\rm $V_{q_t}$-Fueter variables}.
\end{definition}

For $t=-1$, corresponding to the field of quaternions $\mathbb H_{-1}$, these variables were introduced in \cite{adv_prim}.\\

With $\alpha=(\alpha_1,\alpha_2,\alpha_3)\in\mathbb N_0^3$ we recall that
\[
|\alpha|=\alpha_1+\alpha_2+\alpha_3
  \]
and we use  the  multi-variable notation for $x=(x_1,x_2,x_3)\in\mathbb R^3$, and more generally for commuting variables:
\[
\begin{split}
  x^\alpha&=x_1^{\alpha_1}x_2^{\alpha_2}x_3^{\alpha_3}.
\end{split}
\]
We can now define:
\begin{definition}
\label{mu-alpha}
For the $V_t$-Fueter variables $\mu_u$ and $\alpha=(\alpha_1,\alpha_2,\alpha_3)\in\mathbb N_0^3$, we define the product:
\begin{equation}
  \mu_t^\alpha(x)=  \mu_{1,t}^{\alpha_1}(x)  \mu_{2,t}^{\alpha_2}(x)  \mu_{3,t}^{\alpha_3}(x),
\end{equation}
for every $x\in\mathbb{H}_t^*$.
\end{definition}

\begin{proposition}
  The $V_{t}$-Fueter variables pairwise commute and for $\alpha\in\mathbb N_0^3$ it holds that
  \begin{equation}
    \begin{split}
      \mu_{t}^\alpha&=\mu_{1,t}^{\alpha_1}\mu_{2,t}^{\alpha_2}\mu_{3,t}^{\alpha_3}\\
      &=x^\alpha\left(1+\frac{x_0}{\vec{q_t}}\right)^{|\alpha|}
      \end{split}
    \end{equation}
    and
    \begin{equation}
      \mu_{t}^\alpha\big|_{x_0=0}=x^\alpha
    \end{equation}
  \end{proposition}
 \begin{proof}
   Since here too the  $V_{t}$-Fueter variables $\mu_l$ commute, in the above definition we do not need to use the symmetric product and the proof
   goes along the lines of the quaternionic case which was presented in  \cite{adv_prim}. 
 \end{proof}

We can now prove that the $V_t$-Fueter products $\mu_t^\alpha$ are in the kernel of $V_t$ on $\mathbb{H}_t^*$:
\begin{theorem} 
\label{KerV_t}

It holds that $\mu_t^\alpha$ are in the kernel of the operator $V_t$ on any open domain $\Omega\subset \mathbb{H}_t^*$. Moreover, we have:
  \begin{equation}
V_t\mu_t^\alpha(x)=0,
    \end{equation}
    for every $x\in\Omega$.
  \end{theorem}
  
  \begin{proof} Following the argument in~\cite{adv_prim}, we divide the verification into a number of steps. \\

    STEP 1: {\sl It holds that
      \begin{eqnarray}
        \label{kobe_t}
\frac{\partial }{\partial x_1}\frac{1}{\vec{q_t}}=\frac{i}{\vec{q_t\,}^2}+\frac{2x_1}{\vec{q_t\,}^3},\nonumber\\
\frac{\partial }{\partial x_2}\frac{1}{\vec{q_t}}=\frac{j_t}{\vec{q_t\,}^2}+\frac{2x_2 t}{\vec{q_t\,}^3}, \\
\frac{\partial }{\partial x_3}\frac{1}{\vec{q_t}}=\frac{k_t}{\vec{q_t\,}^2}+\frac{2x_3 t}{\vec{q_t\,}^3}.\nonumber
\end{eqnarray}
}

Indeed, for the first case we have
\[
\frac{1}{\vec{q_t}}=\frac{\vec{q_t}}{\vec{q_t\,}^2}=\frac{-\vec{q_t}}{x_1^2+(x_2^2+x_3^2)t}.
  \]
  Hence
  \[
\begin{split}
  \frac{\partial }{\partial x_1}\frac{1}{\vec{q_t}}&=\frac{-i (x_1^2+(x_2^2+x_3^2)t)
    +2x_1\vec{q_t}}{(x_1^2+(x_2^2+x_3^2)t)^2}\\
  &=\frac{i\vec{q_t}^{\,2}}{\vec{q_t}^{\,4}}+\frac{2x_1\vec{q_t}}{\vec{q_t}^{\,3}}
\end{split}
\]
and the result follows. The other cases are similar and equations~\eqref{kobe_t} are proven.\smallskip


STEP 2: {\sl It now holds that
  \begin{eqnarray}
    \label{osaka123_t}
    \frac{\partial}{\partial x_1}\left(1+\frac{x_0}{\vec{q_t}}\right)^{|\alpha|}&=
   \sum_{\substack{t,s\in\mathbb N_0\\t+s=|\alpha|}}
   \left(1+\frac{x_0}{\vec{q_t}}\right)^t\left(\frac{x_0i}{\vec{q_t\,}^2}+\frac{2x_1x_0}{\vec{q_t\,}^3}\right)
   \left(1+\frac{x_0}{\vec{q_t}}\right)^s,\nonumber \\
    \frac{\partial}{\partial x_2}\left(1+\frac{x_0}{\vec{q_t}}\right)^{|\alpha|}&=
    \sum_{\substack{t,s\in\mathbb N_0\\t+s=|\alpha|}}
    \left(1+\frac{x_0}{\vec{q_t}}\right)^t\left(\frac{x_0j_t}{\vec{q_t\,}^2}+\frac{2t x_2x_0}{\vec{q_t\,}^3}\right)
   \left(1+\frac{x_0}{\vec{q_t}}\right)^s,\\
    \frac{\partial}{\partial x_3}\left(1+\frac{x_0}{\vec{q_t}}\right)^{|\alpha|}&=
    \sum_{\substack{t,s\in\mathbb N_0\\t+s=|\alpha|}}
    \left(1+\frac{x_0}{\vec{q_t}}\right)^t\left(\frac{x_0k_t}{\vec{q_t\,}^2}+\frac{2t x_3x_0}{\vec{q_t\,}^3}\right)
   \left(1+\frac{x_0}{\vec{q_t}}\right)^s.\nonumber
  \end{eqnarray}
  }
\eqref{osaka123_t} is a direct consequence of \eqref{kobe_t} and of the formula for the derivative of $f^n$ when $f$ is a matrix-valued function (and in particular a global quaternionic valued) of (say) a real variable $w$:
\begin{equation}
  \label{rtyu}
\frac{{\rm d} f^n}{{\rm d}w}=\sum_{\substack{t,s\in\mathbb N_0\\t+s=|\alpha|}}f^tf^\prime f^s.
\end{equation}

STEP 3: {\sl We have
  \begin{equation}
    \label{osaka12345_t}
    \frac{\partial}{\partial x_0}\left(1+\frac{x_0}{\vec{q_t}}\right)^{|\alpha|}=
    \frac{|\alpha|}{\vec{q_t}}    \left(1+\frac{x_0}{\vec{q_t}}\right)^{|\alpha|-1}.
  \end{equation}
}
This is because $1+\frac{x_0}{\vec{q_t}}$ commutes with its derivative with respect to $x_0$, and formula \eqref{rtyu}
reduces then to the classical formula.\\

STEP 4: {\sl We now calculate
  \[
\frac{1}{\vec{q_t}}
    \sum_{u=1}^3x_u\frac{\partial}{\partial x_u}\mu^\alpha\, .
    \]
}

We have:
\[
  \begin{split}
    &\frac{1}{\vec{q_t}}
    \sum_{u=1}^3x_u\frac{\partial}{\partial x_u}x^\alpha\left(1+\frac{x_0}{\vec{q_t}}\right)^{|\alpha|}=\frac{1}{\vec{q_t}}\cdot \\
   & \big(x_1\frac{\partial}{\partial x_1}x^\alpha\left(1+\frac{x_0}{\vec{q_t}}\right)^{|\alpha|}+x_2\frac{\partial}{\partial x_2}x^\alpha\left(1+\frac{x_0}{\vec{q_t}}\right)^{|\alpha|}+x_3\frac{\partial}{\partial x_3}x^\alpha\left(1+\frac{x_0}{\vec{q_t}}\right)^{|\alpha|}\big)\\
    &=
\frac{1}{\vec{q_t}}\left[x_1 \alpha_1 x^{\alpha-(1,0,0)}\left(1+\frac{x_0}{\vec{q_t}}\right)^{|\alpha|} + x_1 x^\alpha  \sum_{\substack{t,s\in\mathbb N_0\\t+s=|\alpha|}}
   \left(1+\frac{x_0}{\vec{q_t}}\right)^t\left(\frac{x_0i}{\vec{q_t\,}^2}+\frac{2x_1x_0}{\vec{q_t\,}^3}\right)
   \left(1+\frac{x_0}{\vec{q_t}}\right)^s \right] + \\
   &+\frac{1}{\vec{q_t}}\left[x_2 \alpha_2 x^{\alpha-(0,1,0)}\left(1+\frac{x_0}{\vec{q_t}}\right)^{|\alpha|} + x_2 x^\alpha  \sum_{\substack{t,s\in\mathbb N_0\\t+s=|\alpha|}}
   \left(1+\frac{x_0}{\vec{q_t}}\right)^t\left(\frac{x_0j_t}{\vec{q_t\,}^2}+\frac{2t x_2x_0}{\vec{q_t\,}^3}\right)
   \left(1+\frac{x_0}{\vec{q_t}}\right)^s \right] +\\
   &+\frac{1}{\vec{q_t}}\left[x_3 \alpha_3 x^{\alpha-(0,0,1)}\left(1+\frac{x_0}{\vec{q_t}}\right)^{|\alpha|} + x_3 x^\alpha  \sum_{\substack{t,s\in\mathbb N_0\\t+s=|\alpha|}}
   \left(1+\frac{x_0}{\vec{q_t}}\right)^t\left(\frac{x_0k_t}{\vec{q_t\,}^2}+\frac{2t x_3x_0}{\vec{q_t\,}^3}\right)
   \left(1+\frac{x_0}{\vec{q_t}}\right)^s \right] ,\\
&=\frac{1}{\vec{q_t}}\left[|\alpha|x^\alpha\left(1+\frac{x_0}{\vec{q_t}}\right)^{|\alpha|}\right]\\
&+\frac{1}{\vec{q_t}}\left[x^\alpha\sum_{\substack{t,s\in\mathbb N_0\\t+s=|\alpha|-1}}
\left(1+\frac{x_0}{\vec{q_t}}\right)^t\left(\frac{x_1x_0i+x_2x_0j_t+x_3x_0k_t}{\vec{q_t\,}^2}+\frac{2x_1^2x_0+2t x_2^2 x_0+2tx_3^2x_0}
  {\vec{q_t\,}^3}\right)
\left(1+\frac{x_0}{\vec{q_t}}\right)^s\right]\\
&=\frac{x^\alpha}{\vec{q_t}}\left[|\alpha|\left(1+\frac{x_0}{\vec{q_t}}\right)^{|\alpha|}
+\sum_{\substack{t,s\in\mathbb N_0\\t+s=|\alpha|-1}}
\left(1+\frac{x_0}{\vec{q_t}}\right)^t\left(\frac{x_0\vec{q_t}}{\vec{q_t\,}^2}-\frac{2\vec{q_t\,}^2x_0}
  {\vec{q_t\,}^3}\right)
\left(1+\frac{x_0}{\vec{q_t}}\right)^s\right]\\
&=\frac{x^\alpha}{\vec{q_t}}\left[|\alpha|\left(1+\frac{x_0}{\vec{q_t}}\right)^{|\alpha|}
-\sum_{\substack{t,s\in\mathbb N_0\\t+s=|\alpha|-1}}
\left(1+\frac{x_0}{\vec{q_t}}\right)^t\frac{x_0}{\vec{q_t}}
\left(1+\frac{x_0}{\vec{q_t}}\right)^s\right]\\
&=
\frac{x^\alpha}{\vec{q_t}}\left[|\alpha|\left(1+\frac{x_0}{\vec{q_t}}\right)^{|\alpha|}
-|\alpha|\left(1+\frac{x_0}{\vec{q_t}}\right)^{|\alpha|-1}\frac{x_0}{\vec{q_t}}\right].
    \end{split}
  \]

  STEP 5: {\rm We can now compute $V_t\mu_t^\alpha$.}\\

 Using \eqref{osaka12345_t} and the previous step we have:
  \[
\begin{split}
    V_t\mu_t^\alpha&=    x^\alpha\frac{|\alpha|}{\vec{q_t}}  \left(1+\frac{x_0}{\vec{q_t}}\right)^{|\alpha|-1}-
\frac{x^\alpha}{\vec{q_t}}\left[|\alpha|\left(1+\frac{x_0}{\vec{q_t}}\right)^{|\alpha|}
  -|\alpha|\left(1+\frac{x_0}{\vec{q_t}}\right)^{|\alpha|-1}\frac{x_0}{\vec{q_t}}\right]\\
&=\frac{|\alpha|x^\alpha}{\vec{q_t}}\left(1+\frac{x_0}{\vec{q_t}}\right)^{|\alpha|-1}
\left(  1-\left(1+\frac{x_0}{\vec{q_t}}\right)+\frac{x_0}{\vec{q_t}}\right)\\
&=0.
\end{split}
\]
This ends the proof.
 \end{proof}

 \begin{remark}
We observe that both functions $\mu_t^\alpha$ and $\zeta_t^\alpha$ coincide with $x^\alpha$ when $x_0=0$. It is important to note that these are two different extensions of the same real function $x^\alpha$ leading to two different regular function theories. In fact, $\mu_t^\alpha$ is the $V_t$-regular extension of $x^\alpha$ while $\zeta_t^\alpha$ gives the classical global Fueter extension. However, the classical global Fueter variables $\zeta_t^\alpha$ extend $x^\alpha$ to the whole space of quaternions while $\mu_t^\alpha$ extend $x^\alpha$ to domains of $\mathbb{H}_t^*$.     
    \end{remark}

 In their work paper~\cite{alpay_cho_4}
 two of the co-authors have introduced the extension of the Fueter variables corresponding to the Fueter operator to the scaled case,
 namely to scaled Fueter variables for the scaled Fueter operator:

\begin{definition}
\label{zeta_t_def}
The scaled Fueter $\zeta_t$ variables are:
\begin{equation}
\label{zeta_t}
\zeta_{1,t}=x_1-x_0 i, \qquad \zeta_{2,t}=x_2+\frac{\sgn(t) x_0}{\sqrt{|t|}}j_t, \qquad \zeta_{3,t}=x_3+\frac{\sgn(t) x_0}{\sqrt{|t|}}k_t.
\end{equation}
\end{definition}

In the case of $t=-1$, namely, in the case of quaternions, we have the following (see~\cite{adv_prim} for details)
\begin{proposition}
 Let $\Omega$ be an open domain in $\mathbb{H}^*$. For every $x\in\Omega\subset \mathbb{H}^*$ it holds that
  \begin{equation}
    |\mu^\alpha|^2=|x|^{2\alpha}\left(1+\frac{x_0^2}{x_1^2+x_2^2+x_3^2}\right)^{|\alpha|},
    \end{equation}
    and, in particular, we have:
    \begin{equation}
      \label{ineq987}
|\mu_u(x)|^2\leq |\zeta_u(x)|^2,\quad u=1,2,3.
      \end{equation}
  \end{proposition}

  \begin{proof}
This follows from the fact that $\frac{x_0}{\vec{q}}$ has no real part on $\Omega$ and:
    \begin{equation}
|\mu_u(x)|^2=x_u^2+\frac{x_u^2x_0^2}{x_1^2+x_2^2+x_3^2}\le x_u^2+x_0^2=|\zeta_u(x)|^2.
      \end{equation}
    \end{proof}
    In the case of the scaled    hypercomoplex numbers, using the scaled Euclidean norm we have a similar result:
    
  \begin{proposition}
 Let $\Omega_t$ be an open domain in $\mathbb{H}_t^*$. For every $x\in\Omega\subset \mathbb{H}_t^*$ it holds that
  \begin{equation}
    |\mu_t^\alpha|_E^2=|x|_E^{2\alpha}\left(1+\frac{x_0^2}{x_1^2+|t|(x_2^2+x_3^2)}\right)^{|\alpha|},
    \end{equation}
    and, in particular, we have:
    \begin{equation}
      \label{ineq987}
|\mu_{u,t}(x)|_E^2\leq |\zeta_{u,t}(x)|^2,\quad u=1,2,3.
      \end{equation}
  \end{proposition}

  \begin{proof}
This follows from the fact that $\frac{x_0}{\vec{q}}$ has no real part on $\Omega_t$ and it is easy to check that:
    \begin{equation*}
|\mu_{u,t}(x)|_E^2=x_u^2+\frac{x_u^2x_0^2}{x_1^2+|t|(x_2^2+x_3^2)}\le x_u^2+\frac{x_0^2}{|t|}=|\zeta_{u,t}(x)|_E^2,
      \end{equation*}
      for $u=2,3$ and that:
         \begin{equation*}
|\mu_{1,t}(x)|_E^2=x_1^2+\frac{x_1^2x_0^2}{x_1^2+{|t|}(x_2^2+x_3^2)}\le x_1^2+x_0^2=|\zeta_{1,t}(x)|_E^2,
      \end{equation*}
    \end{proof}

  \begin{proposition}
For any $n\in\mathbb{N}$, the function $q_t^n$ is in $\ker V_t$ and that, moreover ,
    \[
q_t^n=\sum_{|\alpha|=n}\mu_t^\alpha c_{\alpha,n}
\]
where, with $\alpha=(\alpha_1,\alpha_2,\alpha_3)$,
\begin{equation}
c_{\alpha,n}=\frac{n!}{\alpha!} i^{\times \alpha_1}\times\, j_t^{\times \alpha_2}\times\, k_t^{\times \alpha_3},
\label{cocosisisicocosi_t}
\end{equation}
where the symmetric product is taken among all the products of the units $\mathbf{e_u}$.
  \end{proposition}

  \begin{proof}
  In $\mathbb{H}_t^*$ we have 
    \[
  \begin{split}
    q_t^n&=(x_0+\vec{q_t})^n\\
    &=\left(1+\frac{x_0}{\vec{q_t}}\right)^n  \left(\vec{q_t}\,\right)^n\\
    &=\left(1+\frac{x_0}{\vec{q_t}}\right)^n  \left(\sum_{|\alpha|=n} x^\alpha c_{\alpha,n}\right)\\
      &=\sum_{|\alpha|=n}x^\alpha\left(1+\frac{x_0}{\vec{q_t}}\right)^n c_{\alpha,n},
    \end{split}
  \]
  for some $c_{\alpha,n}\in\mathbb H_t$ which can be expressed in term of symmetrized products as in
  \eqref{cocosisisicocosi_t}
    by known formulas.
\end{proof}

We note that \eqref{cocosisisicocosi_t} does not take into account  the Cayley table of multiplication for the quaternions.


\subsection{Left Regular Functions on $\mathbb{H}_{t}$ where $t\protect\neq0$
in $\mathbb{R}$}

In this section, we review some main results of \cite{alpay_cho_4},
providing motivations of our works. Let $t\in\mathbb{R}$ and $\mathbb{H}_{t}=\mathrm{span}_{\mathbb{R}}\left(\left\{ 1,i,j_{t},k_{t}\right\} \right)$,
the $t$-scaled hypercomplexes. We consider functions,
\[
f:\mathbb{H}_{t}\rightarrow\mathbb{H}_{t},
\]
in the $t$-scaled hypercomplex variable,
\[
w=x_{0}+x_{1}i+x_{2}j_{t}+x_{3}k_{t},\;\mathrm{with\;}x_{l}\in\mathbb{R},
\]
for all $l=0,1,2,3$. Especially, we are interested in $\mathbb{R}$-differentiable functions on an open connected set $\Omega$ of the $t$-hypercomplex
$\mathbb{R}$-space $\mathbb{H}_{t}$. For motivations, see e.g.,
\cite{alss_IJM}, \cite{as2}, \cite{MR0265618},
\cite{MR2822209} and \cite{bibref14}.
\begin{definition}
Let $t\in\mathbb{R}$ be fixed, and $U\subset\mathbb{H}_{t}$, an
open subset. Define a set,
\[
\mathcal{F}_{t.U}\overset{\textrm{def}}{=}\left\{ f:\mathbb{H}_{t}\rightarrow\mathbb{H}_{t}\left|f\textrm{ is }\mathbb{R}\textrm{-differentiable on }U\right.\right\} .
\]
\end{definition}

Remark that the openness on $\mathbb{H}_{t}$, here, is determined
by the semi-norm topology for $\mathbb{H}_{t}$ in terms of the semi-norm
$\left\Vert .\right\Vert _{t}$ of (2.1.10), for a fixed $t\in\mathbb{R}$.
We first, let $t\neq0$ in $\mathbb{R}$, and $\mathcal{F}_{t,U}$,
the family of $\mathbb{R}$-differentiable functions on an open subset
$U$ of the $t$-scaled hypercomplexes $\mathbb{H}_{t}$.
Define the
differential operators $\nabla_{t}$ and $\nabla_{t}^{\ct}$ on
$\mathcal{F}_{t,U}$ by
\begin{equation}
\label{nabla}
\nabla_{t}=\frac{\partial}{\partial x_{0}}+i\frac{\partial}{\partial x_{1}}-j_{t} \frac{\sgn(t)}{\sqrt{\left|t\right|}}\frac{\partial}{\partial x_{2}}-k_{t}\frac{\sgn(t)}{\sqrt{\left|t\right|}}\frac{\partial}{\partial x_{3}},
\end{equation}
and
\begin{equation}
\label{nabla_ct}
\nabla_{t}^{\ct}=\frac{\partial}{\partial x_{0}}-\frac{\partial}{\partial x_{1}}i+\frac{\sgn(t)}{\sqrt{\left|t\right|}}\frac{\partial}{\partial x_{2}} j_{t}+\frac{\sgn(t)}{\sqrt{\left|t\right|}}\frac{\partial}{\partial x_{3}} k_{t},
\end{equation}
where
\[
\sgn\left(t\right)=\left\{ \begin{array}{ccc}
1 &  & \mathrm{if\;}t>0\\
-1 &  & \mathrm{if\;}t<0,
\end{array}\right.
\]
for all $t\in\mathbb{R}\setminus\left\{ 0\right\} $. 

The definitions~\eqref{nabla} and~\eqref{nabla_ct} of the differential operators $\nabla_{t}$
and $\nabla_{t}^{\ct}$ on $\mathcal{F}_{t,U}$ generalizes those
on $\mathcal{F}_{-1}$ (the quaternionic case) and those on $\mathcal{F}_{1}$
(the split-quaternionic case), studied in \cite{MR2124899}, \cite{alss_IJM},
\cite{as2}, \cite{MR0265618}, \cite{MR2822209}
\begin{theorem}
If $\nabla_{t}$ and $\nabla_{t}^{\ct}$ are the differential
operators of~\eqref{nabla} and~\eqref{nabla_ct}, then

\medskip{}
\begin{equation}
  \nabla_{t}^{\ct}\nabla_{t}=\frac{\partial^{2}}{\partial x_{0}^{2}}+\frac{\partial^{2}}{\partial x_{1}^{2}}-\sgn\left(t\right)\frac{\partial^{2}}{\partial x_{2}^{2}}-\sgn\left(t\right)\frac{\partial^{2}}{\partial x_{3}^{2}}
\label{t-dAlembert_ct}
\end{equation}
\end{theorem}

\begin{proof}
See \cite{alpay_cho_4} for details.
\end{proof}
By \eqref{t-dAlembert_ct}, one obtains the following corollary.
\begin{corollary}
Let $\nabla_{t}$ and $\nabla_{t}^{\ct}$ be the differential
operators on $\mathcal{F}_{t}$, and $\varDelta_{t}$, the $t$-Laplacian.
Then
\[
\varDelta_{t}=\left\{ \begin{array}{ccc}
\varDelta_{-1}=\frac{\partial^{2}}{\partial x_{0}^{2}}+\frac{\partial^{2}}{\partial x_{1}^{2}}+\frac{\partial^{2}}{\partial x_{2}^{2}}+\frac{\partial^{2}}{\partial x_{3}^{2}} &  & \mathrm{if\;}t<0\\
\\
\varDelta_{1}=\frac{\partial^{2}}{\partial x_{0}^{2}}+\frac{\partial^{2}}{\partial x_{1}^{2}}-\frac{\partial^{2}}{\partial x_{2}^{2}}-\frac{\partial^{2}}{\partial x_{3}^{2}} &  & \mathrm{if\;}t>0.
\end{array}\right.
\]
\end{corollary}

\begin{proof}
It is proven by \eqref{t-dAlembert_ct}, whenever $t\neq0$ in $\mathbb{R}$.
\end{proof}
\begin{definition}
Let $\nabla_{t}$ be the operator~\ref{nabla} on $\mathcal{F}_{t,U}$,
and $f\in\mathcal{F}_{t,U}$. If 
\[
\nabla_{t}f=\frac{\partial f}{\partial x_{0}}+i\frac{\partial f}{\partial x_{1}}-j_{t}\frac{\sgn\left(t\right)}{\sqrt{\left|t\right|}}\frac{\partial f}{\partial x_{2}}-k_{t}\frac{\sgn\left(t\right)}{\sqrt{\left|t\right|}}\frac{\partial f}{\partial x_{3}}=0,
\]
then $f$ is said to be left $t$(-scaled)-regular on $U$. If
\[
f\nabla_{t}=\frac{\partial f}{\partial x_{0}}+\frac{\partial f}{\partial x_{1}}i-\frac{\sgn\left(t\right)}{\sqrt{\left|t\right|}}\frac{\partial f}{\partial x_{2}}j_{t}-\frac{\sgn\left(t\right)}{\sqrt{\left|t\right|}}\frac{\partial f}{\partial x_{3}}k_{t}=0,
\]
then $f$ is said to be right $t$(-scaled)-regular on $U$. If $f\in\mathcal{F}_{t,U}$
is both left and right $t$-regular, then it is called $t$(-scaled)-regular.

While, a function $f\in\mathcal{F}_{t,U}$ is $t$(-scaled)-harmonic,
if
\[
\varDelta_{t}f=\frac{\partial^{2}f}{\partial x_{0}^{2}}+\frac{\partial^{2}f}{\partial x_{1}^{2}}-\sgn\left(t\right)\frac{\partial^{2}f}{\partial x_{2}^{2}}-\sgn\left(t\right)\frac{\partial^{2}f}{\partial x_{3}^{2}}=0,
\]
where $\varDelta_{t}$
is the $t$-Laplacian \eqref{t-dAlembert_ct}.
\end{definition}

The following theorem illustrates the relation between $t$-regularity
and $t$-harmonicity.
\begin{theorem}
\label{regular-harmonic}
Let $f\in\mathcal{F}_{t,U}$, for $U\in\mathcal{T}_{t}$ in $\mathbb{H}_{t}$. Then if $f$ is left $t$-regular on $U$, it is $t$-harmonic on  $U$.
\end{theorem}

\begin{proof}
If $f$ is left $t$-regular in $\mathcal{F}_{t,U}$, then $\nabla_{t}f=0$
on $U$, and hence,
\[
\varDelta_{t}f=\nabla_{t}^{\ct}\left(\nabla_{t}f\right)=\nabla_{t}^{\ct}\left(0\right)=0,
\]
implying the $t$-harmonicity of $f$ on $U$.
\end{proof}
In Definition~\ref{zeta_t_def} we have introduced the entire $\mathbb{R}$-differentiable functions $\left\{ \zeta_{l}\right\} _{l=1}^{3}\subset\mathcal{F}_{t,\mathbb{H}_{t}}$, 
associated to a $\mathbb{H}_{t}$-variable $q_t=x_{0}+x_{1}i+x_{2}j_{t}+x_{3}k_{t}$,
with $x_l\in\mathbb{R},\, 0\le l\le 3$. In~\cite{alpay_cho_4} the following result was proven:

\begin{theorem}
If $\zeta_{l,t}$ are the functions~\eqref{zeta_t} in $\mathcal{F}_{t,\mathbb{H}_{t}}$,
then $\zeta_{l,t}$ are $t$-harmonic and $t$-regular on $\mathbb{H}_{t}$,
$\forall\, l=1,2,3.$
\end{theorem}

Now, define the so-called symmetrized product on $\mathbb{H}_{t}$.
One can do that since the $\mathbb{R}$-vector space $\mathbb{H}_{t}$
has a well-defined vector-multiplication $\left(\cdot\right)=\left(\cdot_{t}\right)$.
\begin{definition}
Let $h_{1},...,h_{N}\in\mathbb{H}_{t}$, for $N\in\mathbb{N}$. Then
the symmetrized product of $h_{1},...,h_{N}$ is defined by a new
hypercomplex number,

\medskip{}

\begin{equation}
\label{cross-product}
\overset{N}{\underset{n=1}{\times}}h_{n}\overset{\textrm{denote}}{=}h_{1}\times...\times h_{N}\overset{\textrm{def}}{=}\frac{1}{N!}\underset{\sigma\in S_{N}}{\sum}h_{\sigma(1)}h_{\sigma(2)}...h_{\sigma\left(N\right)},
\end{equation}

\medskip{}

\noindent where $S_{N}$ is the symmetric group over $\left\{ 1,...,N\right\} $.
\end{definition}

\begin{remark}Recall that we are considering
the cases where $t\neq0$ in $\mathbb{R}$, however, the symmetrized
product ($\times$) of~\eqref{cross-product} is well-defined for ``all'' $t\in\mathbb{R}$,
including the case where $t=0$.
\end{remark}

Let $f_{1},...,f_{N}:\mathbb{H}_{t}\rightarrow\mathbb{H}_{t}$ be
functions for $N\in\mathbb{N}$. By applying~\eqref{cross-product}, define a symmetrized-product
function of $f_{1},...,f_{N}$ by
\begin{equation}
\label{symm-prod}
\overset{N}{\underset{n=1}{\times}}f_{n}=\frac{1}{N!}\underset{\sigma\in S_{N}}{\sum}f_{\sigma(1)}f_{\sigma(2)}...f_{\sigma(N)},
\end{equation}
where
\begin{equation*}
\left(f_{\sigma(1)}f_{\sigma(2)}...f_{\sigma(N)}\right)\left(h\right)=f_{\sigma(1)}\left(h\right)...f_{\sigma\left(N\right)}\left(h\right),\;\forall h\in\mathbb{H}_{t},
\end{equation*}
for all $\sigma\in S_{N}$. So, one obtains
\[
\mathit{f_{j}^{(n)}}=\underset{n\textrm{-times}}{\underbrace{f_{j}\times f_{j}\times...\times f_{j}}}=f_{j}^{n},\;\;\forall n\in\mathbb{N},
\]
where we define $f^0$ to be:
\[
f^{(0)}=1,\;\mathrm{the\;constant\;\textrm{1-}function\;on\;}\mathbb{H}_{t},
\]
for all functions $f:\mathbb{H}_{t}\rightarrow\mathbb{H}_{t}$, where
$1=1+0i+0j_{t}+0k_{t}\in\mathbb{H}_{t}$.
\begin{definition}
Let $\mathbf{n}\overset{\textrm{denote}}{=}\left(n_{1},n_{2},n_{3}\right)\in\mathbb{N}_{0}^{3}$
be a triple of numbers in $\mathbb{N}_{0}=\mathbb{N}\cup\left\{ 0\right\} $,
and let $\left\{ \zeta_{l}\right\} _{l=1}^{3}$ be the $t$-harmonic
$t$-regular functions of~\eqref{zeta_t}. Define a new function $\zeta^{\mathbf{n}}\in\mathcal{F}_{1,\mathbb{H}_{t}}$
by
\begin{equation}
\label{zeta^n}
\zeta^{\mathbf{n}}\overset{\textrm{def}}{=}\frac{1}{\mathbf{n}!}\left(\zeta_{1}^{(n_{1})}\times\zeta_{2}^{(n_{2})}\times\zeta_{3}^{(n_{3})}\right),
\end{equation}
where
\begin{equation*}
\mathbf{n}!=\left(n_{1}!\right)\left(n_{2}!\right)\left(n_{3}!\right)\in\mathbb{N},
\end{equation*}
and
\[
\zeta_{l+1}^{(n_{l})}=\;\underset{n_{l}\textrm{-times}}{\underbrace{\zeta_{l+1}\times\zeta_{l+1}\times...\times\zeta_{l+1}}}=\zeta_{l+1}^{n_{l}},\;\forall l=1,2,3.
\]
\end{definition}

\begin{theorem}
\label{zeta-harmonic}
Let $\mathbf{n}=\left(n_{1},n_{2},n_{3}\right)\in\mathbb{N}_{0}^{3}$,
and let $\zeta^{\mathbf{n}}\in\mathcal{F}_{t,\mathbb{H}_{t}}$ be the
function~\eqref{zeta^n}. Then it is a $t$-harmonic $t$-regular
function, i.e.,
\begin{equation*}
  \zeta^{\mathbf{n}}\nabla_{t}=\nabla_{t}\zeta^{\mathbf{n}}=0,\;and\;\varDelta_{t}\zeta^{\mathbf{n}}=0,\;\mathrm{on\;}\mathbb{H}_{t}.
  \end{equation*}
\end{theorem}

\begin{proof}
Theorem~\ref{zeta-harmonic} follows from Theorem~\ref{regular-harmonic}. See
\cite{alpay_cho_4} for details.
\end{proof}
By Theorem~\ref{zeta-harmonic}, all functions $\left\{ \zeta^{\mathbf{n}}:\mathbf{n}\in\mathbb{N}_{0}^{3}\right\} $
of \eqref{zeta^n} are $t$-harmonic $t$-regular functions on $\mathbb{H}_{t}$. We will now see that they also generate the space of $t$-regular functions in this case.

%
%
%

In~\cite{alpay_cho_4} the authors have proven the following:
\begin{theorem}
Let $t\neq0$ in $\mathbb{R}$, and $f\in\mathcal{F}_{t,U}$, where
$0=0+0i+0j_{t}+0k_{t}\in U$ in $\mathbb{H}_{t}$. If $f$ is $\mathbb{R}$-analytic
on $U$, then
\[
f\textrm{ is left }t\textrm{-regular on }U,
\]
if and only if
\begin{equation}
  \label{zeta_series}
f=f\left(0\right)+\underset{\mathbf{n}\in\mathbb{N}^{3}}{\sum}\zeta^{\mathbf{n}}f_{\mathbf{n}},
\end{equation}
with
\[
f_{\mathbf{n}}=\frac{1}{\mathbf{n}!}\frac{\partial^{n_{1}+n_{2}+n_{3}}f}{\partial x_{1}^{n_{1}}\partial x_{2}^{n_{2}}\partial x_{3}^{n_{3}}}\left(0\right),\;\forall\mathbf{n}\in\mathbb{N}^{3}.
\]
\end{theorem}

\begin{proof}
See \cite{alpay_cho_4} for details. 
\end{proof}
Note that if $t<0$, then the $\mathbb{R}$-analyticity of a $t$-regular
function $f\in\mathcal{F}_{t,U}$ is automatically guaranteed by \eqref{zeta_series}.
Thus, if $t<0$ in $\mathbb{R}$ then, without the $\mathbb{R}$-analyticity
assumption of $f$, the above characterization~\eqref{zeta_series} holds. However,
if $t\geq0$, then the $\mathbb{R}$-analyticity assumption is needed
in \eqref{zeta_series}. See e.g., \cite{alpay_cho_4}.

\begin{remark} Please note that for the Fueter operator in the case of split quaternions other notations exist in literature and we write the equivalence of notation in this table:

\begin{center}
  \begin{tabular}{|c|c|c|c|}
    \hline
    &This paper&Notation in Paper \cite{alss_IJM} \\
    \hline
Variables    &$x_3$&$x_1$ \\
    \hline
    &$x_1$&$x_3$\\
    \hline
  Hypercomplex & $i$&$i$\\
  basis  &$j_1$&$k_2$\\
    &$k_1$&$k_1$\\
    \hline
    Differential &$\nabla_{1}$&$\nabla^+_{\mathbb R}$ \\
    operator &$\frac{\partial }{\partial x_{0}}+i\frac{\partial }{\partial x_{1}}-j_1\frac{\partial }{\partial x_{2}}-k_1\frac{\partial f}{\partial x_{3}}$&
$\frac{\partial }{\partial x_{0}}-k_1\frac{\partial }{\partial x_{1}}-k_2\frac{\partial }{\partial x_{2}}+i\frac{\partial f}{\partial x_{3}}$ \\
        \hline

                                Fueter variables            &                $\zeta_{1,1}(q_1)=x_1-x_0i$                    & $\zeta_1(x)=x_1+x_0k_1$ \\
(3.5)-(3.7) p.332   &      $\zeta_{2,1}(q_1)=x_2+x_0j_1$     &                          $\zeta_2(x)=x_2+x_0 k_2$ \\
 of \cite{alss_IJM}   &$\zeta_{3,1}(q_1)=x_3+x_0k_1$ & $\zeta_3(x)=x_3-x_0 i$ \\
    \hline
     \end{tabular}
\end{center}
\end{remark}

In general, other linear changes of variables can be considered for the setting of Fueter variables, for example in~\cite{MR2852190}.

%
 \section{Hardy space and interpolation for the first conjugate $\circledast$}
 \setcounter{equation}{0}
 \label{Hardy_spaces_cdast}
\subsection{Prelude}
\label{Hardy-prelude}
 The Hardy space $\mathbf H_2(\mathbb D)$ is the reproducing kernel Hilbert space of functions analytic in the open unit disk and with reproducing $\frac{1}{1-z\overline{w}}$. It is also the unique (up to a positive multiplicative factor in the inner product)
 Hilbert space of functions analytic in $\mathbb D$ such that $M_z^*=R_0$, where $M_z$ denotes the operator of multiplication by $z$ and $R_0$ is the backward-shift operator. Important issues pertaining to the classical theory include the Beurling-Lax theorem, Schur multipliers and related interpolation
 problems. In the present section we first study the counterparts of this space in the setting of the rings $\mathcal H_2^t$. For $t=-1$ this was done
 in \cite{MR2872478}. We mention that a difference occured there between the scalar and matrix-valued case in view of the lack of commutativity; see  \cite[p. 1767 (and in particular (62.38)]{acs-survey}).
 The counterpart of the Hardy space is now the following real Hilbert space:

\begin{definition}
We denote by $\mathbf H^t_2$ the space of power series of the form
\begin{equation}
  \label{gamma-n}
f(q_t)=\sum_{n=0}^\infty q_t^n\alpha_n
\end{equation}
where $\alpha_0,\alpha_1,\ldots\in\mathbb H_t$ and such that
\begin{equation}
\sum_{n=0}^\infty\|\alpha_n\|_{E}^2<\infty
\label{norm-12345}
\end{equation}
where the scaled norm \eqref{scaled-E} is used in \eqref{norm-12345}.
  \end{definition}

 Note that \eqref{norm-12345} is induced by the inner-product
  \[
[f,g]_{E}=2\sum_{n=0}^\infty {\rm Re}\,\left(a_n\overline{c_n}+|t|\,b_n\overline{d_n}\right)
    \]
    where 
    \[
\alpha_n=a_n+b_nj_t,\quad n=0,1,\ldots
\]
and with $g(q_t)=\sum_{n=0}^\infty q_t^n\gamma_n$,
    \[
\gamma_n=c_n+d_nj_t,\quad n=0,1,\ldots
\]
We also note that, under \eqref{norm-12345}, the power series  converges for $\|q_t\|<1$, where $\|\cdot\|$ is the operator norm.  The case $t=-1$, corresponding to the quaternionic setting, was first studied in \cite{MR2872478}. In that work and related ones, one views the Hardy space as a quaternionic
Hilbert space. This is not possible for general $t$ since $\mathbb H_{t}$ will not be a skew field in general.  We will endow $\mathbf H^t_2$ with two real symmetric forms, namely
 \begin{equation}
 \label{function-inner-1}
[f,g]_{\circledast}=\sum_{n=0}^\infty [g_n,f_n]_{\circledast},
\end{equation}
and
\begin{equation}
  \label{function-inner-2}
[f,g]_{[*]}=\sum_{n=0}^\infty [g_n,f_n]_{[*]}
\end{equation}
where $f(q_t)=\sum_{n=0}^\infty q_t^n f_n$ and $g(q_t)=\sum_{n=0}^\infty q_t^ng_n$ are elements of $\mathbf H_2^t$.
\begin{theorem}
\label{H_2-Krein}
  The space $\mathbf H_2^t$ endowed with either of the form \eqref{function-inner-1} or \eqref{function-inner-2} is a Krein space. In case of \eqref{function-inner-1} it is a
  Hilbert space for $t<0$. When using the form~\eqref{function-inner-1} we denote this space by  $\mathbf H_{2,\circledast}^t$, and by $\mathbf H_{2,[*]}^t$ when using~\eqref{function-inner-2}.
\end{theorem}

\begin{proof}
  The claim follows from the corresponding structure of the coefficient space $\mathbb H_t$ from Theorems~\ref{pont-ht} and~\ref{bracket-Pontryagin} endowed with the forms
   \begin{equation*}
[p_t,q_t]_{\circledast}={\rm Tr}~(q_t^{\circledast}p_t)=a\overline{c}+\overline{a}c-t(b\overline{d}+\overline{b}d)
\end{equation*}
and
\begin{equation*}
[p_t,q_t]_{[*]}={\rm Tr}~(q_t^{[*]}p_t)=a{c}+\overline{ac}-t(b{d}+\overline{bd}),
\end{equation*}
where
\begin{equation*}
      q_t=\begin{pmatrix}a&tb\\ \overline{b}&\overline{a}\end{pmatrix},\,\qquad p_t=\begin{pmatrix}c&td\\ \overline{d}&\overline{c}\end{pmatrix}.
    \end{equation*}
  \end{proof}

\begin{definition}
\label{Hardy-Krein}
We denote by $\mathbf H_{2,{\ct}}$ the Hardy space $\mathbf H_2^t$ endowed with the form~\eqref{function-inner-1} and by 
$\mathbf H_{2,[*]}^t$  the Hardy space $\mathbf H_2^t$ endowed with the form~\eqref{function-inner-2}, which, by Theorem~\ref{H_2-Krein}, become Hilbert and Krein spaces, respectively, depending on $\sgn t$.
 \end{definition} 
These Hilbert and  Krein spaces are studied in Subsections \ref{Hardy_q_t_1} and \ref{Hardy_q_t} respectively.

 In several complex variables, one considers in Schur analysis the Arveson
 space rather than the Hardy space of the ball. In the  last part of the present section we present the analog of the Arveson space for
 the $\mu_t$ variables. In the cases $t=1$ this was first done in \cite{MR2124899}. \smallskip

     \subsection{The Hardy space in one variable: Case $1$}
    \label{Hardy_q_t_1}
In this subsection, we consider the Hardy space $\mathbf H_{2,\circledast}^t$ when the adjoint $\circledast$ is considered,
\begin{equation}
  \label{not-same}
      q_t^{\circledast}=\begin{pmatrix}\overline{a}&-tb\\ -\overline{b}&a\end{pmatrix}.
    \end{equation}
    We recall that $\mathcal  H^t_2$ endowed with the inner product \eqref{form-qsrar}
    \[
[p_t,q_t]_{\circledast}={\rm Tr}~(q_t^{\circledast}p_t)
\]
is a real Pontryagin space of dimension $4$ and index of negativity $2$, in the case $t>0$ and it is a Hilbert space when $t<0$; see Lemma \ref{pont-ht}.

In what follows we will use the spaces $\mathcal H_2^t$ and $\mathbb H_t$ interchangeably, for example the matrix $I_2$ in the first becomes the number $1$ in the second.
\begin{theorem}
\label{ct-Krein}
$\mathbf H_{2,\circledast}^t$ is a reproducing kernel Krein right $\mathcal H_2^t$-module (or an $\mathbb H_t-$module) with reproducing kernel 
\begin{equation}
  \label{kern-1}
K(q_t,p_t)=\sum_{n=0}^\infty q_t^np_t^{\circledast n}
\end{equation}
when endowed with the form (we use the same notation as for $\mathcal H_2^t$; see \eqref{not-same})
\begin{equation}
  [f,g]_{\ct}=\sum_{n=0}^\infty \delta_n^{\circledast}\gamma_n
\end{equation}
with $g(q_t)=\sum_{n=0}^\infty q_t^n\delta_n$,  and associated inner product (indefinite when $t>0$ and definite when $t<0$)
\begin{equation}
  {\rm Tr}\,[f,g]_{\ct}=\sum_{n=0}^\infty {\rm Tr}\,\delta_n^{\circledast}\gamma_n
  \end{equation}
\end{theorem}

\begin{proof}
  The space of sequences $(\gamma_n)_{n\in\mathbb N_0}$ of elements of $\mathcal H_2^t$ such that \eqref{norm-12345} is in force is a real Hilbert space, and so is $\mathbf H_{2,\circledast}^t$ by transfer of structure. The norm \eqref{norm-12345} corresponds to the inner product
  \[
    \langle f,g\rangle=\sum_{n=0}^\infty {\rm Tr}~{\delta_n^*\gamma_n}
    \]
    where $g(q_t)=\sum_{n=0}^\infty q_t^n\delta_n$. To check the reproducing kernel property, we write, with $f$ as in \eqref{gamma-n} and the kernel defined by \eqref{kern-1}:
\[
  \begin{split}
{\rm Tr}\, [f(\cdot),K(\cdot,p_t)\alpha]_{\ct}&={\rm Tr}\,\left(\sum_{n=0}^\infty \alpha^{\ct}p_t^{n}\gamma_n\right)={\rm Tr}\, \alpha^{\ct}f(p_t).
  \end{split}
\]
\end{proof}

\begin{theorem}
  The operator $M_{q_t}$ of $star$-multiplication  on the left by $q_t$ is a continuous Krein space isometry, with adjoint the backward-shift operator
  \begin{equation}
    \label{M-*}
    ( M_{q_t}^*f)(q_t)=\sum_{n=1}^\infty q_t^{n-1}\alpha_n,\quad {with}\quad f(q_t)=\sum_{n=0}^\infty q_t^n\alpha_n.
  \end{equation}
  We also have that $M_{q_t}^*M_{q_t}f=f$ for $f\in\mathbf H_{2,\circledast}^t$ .
\end{theorem}

\begin{proof}

  We have with $f$ and $g$ as above,
  \[
    \begin{split}
      [M_{q_t}f,g]_{\ct}&={\rm Tr}\, \left(\sum_{n=0}^\infty \delta^{\ct}_{n+1}\gamma_n\right)\\
      &=[f,M_{q_t}^*g]_{\ct}.
    \end{split}
    \]
  \end{proof}


\subsection{Interpolation in the Complex Case}
\label{interpolation-prelude}
In the following subsections we begin a study of interpolation in the spaces $\mathbf H^t_{2,\ct}$ and $\mathbf H^t_{2,[*]}$. In the classical complex scalar setting, the most simple problem would be:\\

{\sl Given $N$ different points $a_1,\ldots, a_N$  in the open unit disk, describe the set of functions $f$
  in the Hardy space
  such that
  \begin{equation}
f(a_j)=0,\quad j=1,\ldots, N.
\end{equation}
}

It is a simple exercise, and a particular instance of Beurling's theorem, that the set of solutions of this problem consists of the functions of the form
$f(z)=b(z)g(z)$ where $b(z)=\prod_{n=1}^N\frac{z-a_n}{1-z\overline{a_n}}$ and $g$ runs through the Hardy space.
To sharpen the difference beetween the following subsection and the following two ones we present the following well-known result from $J$-theory
(see \cite[Theorem 1, p. 145, p. 146 line 5 and Theorem 2 p. 146]{pootapov}, \cite{siegel}, \cite{Dym_CBMS} and \cite[Lemma 4.1.7 p. 80]{MR1638044}
for more information).

 \begin{lemma}
   Let $E_1, E_2\in\mathbb C^{u\times v}$  be strict contractions. Then, the matrices
   \begin{equation}
     \label{formule-e1-e2}
(I_u-E_1E_1^*)^{-1/2}(E_2+E_1)(I_v+E_1^*E_2)^{-1}(I_v-E_1^*E_1)^{1/2}
\end{equation}
and
   \begin{equation}
     \label{formule-e1-e2-e3}
(I_u-E_1E_1^*)^{1/2}(I_v+E_2E_1^*)(E_2+E_1)(I_v-E_1^*E_1)^{-1/2}
\end{equation}
     are also strictly contractive.\\
     Assume that $u=v$ and $E_1E_1^*=cI_u$ for some real number $c>0$. Then, $E_1^*E_1=cI_u$ and
        \begin{equation}
     \label{formule-e1-e2-e4}
(E_2+E_1)(I_v+E_1^*E_2)^{-1}
\end{equation}
and
\begin{equation}
  (I_v+E_1^*E_2)^{-1}(E_2+E_1)
  \end{equation}
     are strictly contractive.
   \end{lemma}

When $u=1$ and replacing $E_1=-a\in\mathbb C^{1\times v}$ and $E_2=z\in\mathbb C^{1\times v}$ such that $zz^*<1$ we get the counterpart of the Blaschke product in the unit ball of $\mathbb C^v$, 
\[
  b_a(z)=\sqrt{1-aa^*}(z-a)(I_v-a^*z)^{-1}(I_v-a^*a)^{-1/2}
\]
which satisfies the formula (with $w\in\mathbb C^{1\times v}$ such that $ww^*<1$)
\begin{equation}
  \label{arv-123456789}
\frac{1-aa^*}{(1-za^*)(1-aw^*)}=\frac{1-b_a(z)b_a(w)^*}{1-zw^*}
  \end{equation}
see \cite{MR84i:32007}.\\

We remark that Blaschke factors, and more generally Blaschke products, are examples of rational functions taking unitary (or co-isometric in the case of $\mathbb C^d$) values on the boundary. 
\subsection{Blaschke factor and interpolation: Case $1$}
\label{bl1}
Since $\alpha\alpha^{\ct}$ is a scalar we can use \eqref{formule-e1-e2-e4} to define the $\circledast$-Blaschke factor.

\begin{definition}
Assume $\|\alpha\|<1$. We define the $\circledast$-Blaschke factor with a zero at $\alpha$ by:
\begin{equation}
 b_\alpha(q_t)= (q_t-\alpha)\star(1-q_t\alpha^{\circledast})^{-\star},
    \label{bl-1}
  \end{equation}
  where $S^{-\star}$ represents the $\star$ inverse of $S$.
\end{definition}

  \begin{lemma} Let $\alpha\in\mathcal H_2^t$ of norm striclty less that $1$.  Then,
    \begin{equation}
      \label{b-a-expansion}
      b_\alpha(q_t)=-\alpha+ \sum_{n=1}^\infty q_t^n(I_2-\alpha\alpha^{\ct}) \alpha^{\ct (n-1)}
\end{equation}
where the convergence in in $\mathbb C^{2\times 2}$.
    \end{lemma}

    \begin{proof}
      For $q_t=xI_2$ where $x$ is real and of absolute value less than $1$ we have
\[
      \begin{split}
        b_{\alpha}(q_t)&=(x-\alpha)(\sum_{n=0}^\infty x^n\alpha^{\ct n})\\
        &=\sum_{n=0}^\infty x^{n+1}\alpha^{\ct n}-\alpha-\sum_{n=1}^\infty x^n\alpha\alpha^{\ct n}\\
        &=-\alpha+\sum_{n=1}^\infty x^n(\alpha^{\ct (n-1)}-\alpha\alpha^{\ct n})\\
        &=-\alpha+        \sum_{n=1}^\infty x^n(I_2-\alpha\alpha^\ct)\alpha^{\ct (n-1)}\\
                &=-\alpha+     (\det\alpha)   \sum_{n=1}^\infty x^n\alpha^{\ct (n-1)}
              \end{split}
    \]
    and hence the result by Proposition \ref{unique-123}.
      \end{proof}

As in the slice-hyperholomorphic setting we have the following new phenomenon:
    
  \begin{lemma} Let $\alpha$ of norm strictly less than $1$. It holds that:
    \label{commu}
  \begin{equation}
( b_\alpha\star b_{\alpha^\circledast})(q_t)= \frac{q_t^2-2q_t{\rm Re \alpha}+\det \alpha}{q_t^2\det \alpha-2q_t{\rm Re}~ \alpha +1}
    \label{bl-1}
\end{equation}
\end{lemma}

\begin{proof}
We can write
  
  \[
    \begin{split}
      ( b_\alpha\star b_{\alpha^\circledast})(q_t)&=(q_t-\alpha)\star\underbrace{(1-q_t\alpha^{\circledast})^{-\star}\star(q_t-\alpha^{\circledast})}_{{\rm \star-commute;\,\,see \,\, Proposition\,\, \ref{B-star-comm}}}\star(1-q_t\alpha^{\circledast})^{-\star}\\
      &=      (q_t-\alpha)\star(q_t-\alpha^{\circledast})\star (1-q_t\alpha^{\circledast})^{-\star}\star(1-q_t\alpha)^{-\star}\\
      &=(q_t^2-2q_t{\rm Re}~\alpha)+\alpha\alpha^{\ct})\star (1-q_t\alpha^{\circledast})^{-\star}\star(1-q_t\alpha)^{-\star}\\
      &=(q_t^2-2q_t{\rm Re}~\alpha)+\alpha\alpha^{\ct})\star( (1-q_t\alpha)\star(1-q_t\alpha^{\ct}))^{-\star}\\
      &=(q_t^2-2q_t{\rm Re}~\alpha)+\alpha\alpha^{\ct})\star(q_t^2\det\alpha-2q_t{\rm Re}\, \alpha+1)^{-\star},
    \end{split}
    \]
and hence the result since the elements of the quotients have real coefficients, and hence the $star$-product can be replaced by the usual product.
  \end{proof}

In view of the following theorem note that there is a technical difficulty:  a densely defined isometry in a Krein (as opposed to Hilbert) space need not be continuous, and so it is not enough to check isometry on a dense set.
    
\begin{theorem}
  \label{thmbabounded}
The operator $M_{b_{\alpha}}$ of of $\star$-multiplication by $b_\alpha$ on the left is a continuous isometry from the Krein space $\mathbf H_{2,\circledast}$ into itself.
\end{theorem}

\begin{proof} We divide the proof into steps.\\

  STEP 1: {\sl $M_{b_a}$is an isometry on a dense set.}\smallskip

  Indeed, since
  \[
  (M_{b_\alpha}q_t^u\gamma)(M_{b_\alpha}q_t^v\delta)^{\ct}=\delta^{\ct}
  [M_{b_\alpha}q_t^u\gamma,M_{b_\alpha}q_t^v]_{\ct}\gamma
  \]
  we can assume $\gamma=\delta=I_2$. We first assume $u=v$. We have
  \[
  \begin{split}
    [M_{b_\alpha}q_t^u,M_{b_\alpha}q_t^u]_{\ct}&=\\
    &\hspace{-2.5cm}=[-q_t^u\alpha+\sum_{n=1}^\infty q_t^{n+u}(I_2-\alpha\alpha^{\ct})\alpha^{\ct(n-1)},
      -q_t^u\alpha+\sum_{m=1}^\infty q_t^{m+u}(I_2-\alpha\alpha^{\ct})\alpha^{\ct(m-1)}]_{\ct}\\
      &\hspace{-2.5cm}=\alpha^{\ct}\alpha+\sum_{n=1}^\infty\alpha^{n-1}(I_2-\alpha\alpha^{\ct})^2\alpha^{\ct(n-1)}\\
      &\hspace{-2.5cm}=\alpha^{\ct}\alpha+(I_2-\alpha\alpha^{\ct})^2(I_2-\alpha\alpha^{\ct})^{-1}\\
      &\hspace{-2.5cm}=I_2,
    \end{split}
  \]
  where we have used that $\alpha\alpha^{\ct}=\alpha^{\ct}\alpha$ to compute
  \[
\sum_{n=1}^\infty\alpha^{n-1}\alpha^{\ct(n-1)}=(I_2-\alpha\alpha^{\ct})^{-1}
    \]
for $\|\alpha\|<1$ and the fact that $[q_t^{n+u},q_t^{m+u}]_{\ct}=0$ if $m\not =n$.\\

We now assume $u\not=v$.  We now have $[q_t^{n+u}q_t^{m+v}]_{\ct}=0$ unless $n+u=m+v$. Without loss of generality we will assume $u>v$,
so that  $m=n+(u-v)>n$ and write:
  \[
  \begin{split}
    [M_{b_\alpha}q_t^u,M_{b_\alpha}q_t^v]_{\ct}&=\\
    &\hspace{-2.5cm}=[-q_t^u\alpha+\sum_{n=1}^\infty q_t^{n+u}(I_2-\alpha\alpha^{\ct})\alpha^{\ct(n-1)}    ,    -q_t^v\alpha    +\sum_{m=1}^\infty q_t^{m+v}(I_2-\alpha\alpha^{\ct})\alpha^{\ct(m-1)}]_{\ct}\\
    &\hspace{-2.5cm}=\underbrace{[-q_t^u\alpha,   \sum_{m=1}^\infty q_t^{m+v}(I_2-\alpha\alpha^{\ct})\alpha^{\ct(m-1)}]_{\ct}}_{{\rm only} \, m=u-v\,\, {\rm contributes}}+
    \underbrace{[\sum_{n=1}^\infty q_t^{n+u}(I_2-\alpha\alpha^{\ct})\alpha^{\ct(n-1)},-q_t^v\alpha]_{\ct}}_{=0\,\, {\rm since}\,\, n+u>v}\\
    &\hspace{-2.0cm}+[\sum_{n=1}^\infty q_t^{n+u}(I_2-\alpha\alpha^{\ct})\alpha^{\ct(n-1)},   \sum_{n=1}^\infty q_t^{n+u}(I_2-\alpha\alpha^{\ct})\alpha^{\ct(n+u-v-1)}]_{\ct}\\
    &\hspace{-2.5cm}=-(I_2-\alpha\alpha^{\ct})\alpha^{u-v}+\alpha^{u-v}(I_2-\alpha\alpha^{\ct})^2 [\sum_{n=1}^\infty q_t^{n+u}\alpha^{\ct(n-1)},   \sum_{n=1}^\infty q_t^{n+u}\alpha^{\ct(n-1)}]_{\ct}\\
    &\hspace{-2.5cm}=-(I_2-\alpha\alpha^{\ct})\alpha^{u-v}+(I_2-\alpha\alpha^{\ct})^2\alpha^{u-v}(I_2-\alpha\alpha^{\ct})^{-1}\\
    &\hspace{-2.5cm}=0.
    \end{split}
  \]

  STEP 2: {\sl $M_{b_a}$ is continuous.}\smallskip
  We note that $M_{b_a}=M_{q_t-\alpha}\circ M_{(1-q_t\alpha^{\ct})^{-1}}$.
The result then follows from the fact that $M_{q_t}$ is continuous of norm equal to $1$ and so
  \begin{equation}
    \|M_{(1-q_t\alpha^{\ct})^{-1}}\|\le\frac{1}{1-\|\alpha\|}<\infty
    \label{mbounded}
    \end{equation}
    It follows that $M_{b_a}$ has a continuous extension to $\mathbf H_2^t$, which is an isometry since the form \eqref{function-inner-1} is continuous with respect to the inner product of $\mathbf H_2^t$.
\end{proof}

\begin{theorem}
  A function $f\in\mathbf H_{2,\ct}$ vanishes at $\alpha$ if and only if it can be written as $f=b_\alpha\star g$ where $g\in\mathbf H_{2,\ct}$
  and, moreover
  \[
    [f,f]_{\ct}=[g,g]_{\ct}
  \]
  \label{inter-pol-1}
  \end{theorem}

  \begin{proof}
    We proceed in a number of steps; the arguments are an adaptation of the classical arguments to the present setting. We denote by $\mathfrak M$ the real linear space spanned by the functions $(1-q_t\alpha)^{-\star}b$, where $b$ varies in $\mathcal H_2^t$.\\

    STEP 1: {\sl For any $g\in\mathbf H_2^t$ and $\alpha$ such that $\|\alpha\|<1$ (operator norm), we have
      \begin{equation}
        \label{ope-zero}
(b_\alpha\star g)(\alpha)=0.
        \end{equation}
      }\\
      Indeed, with $h(q_t)=(I_2-q_t\alpha)^{-\star}\star g(q_t)\stackrel{{\rm ref}}{=}\sum_{n=0}^\infty q_tt^nh_n$ we can write
      \[
        \begin{split}
          b_\alpha(q_t)\star g(q_t)&=(q_t-\alpha)\star h(q_t)\\
          &=\sum_{n=0}^\infty((q_t-\alpha)\star q_t^n)h_n\\
          &=\sum_{n=0}^\infty (q_t^{n+1}-q_t^n\alpha)h_n
        \end{split}
      \]
      which vanishes for $q_t=\alpha$.\\

      STEP 2: {\sl We have $\mathfrak M\subset \mathbf H_2^t\ominus b_a\star \mathbf H_2^t$.}\smallskip

    It suffices to verify that $\mathfrak M$ and $b_a\star \mathbf H_2^t$ are orthogonal.  For $f\in\mathbf H_2^t$ we have
\[
  [f(q_t),(1-q_t\alpha)^{-\star}b]_{\ct}={\rm Tr}\, b^{\ct}f(\alpha) = [ f(\alpha), b]_{\circledast}
\]
and so the claim follows from Step 1.\\

    STEP 3: {\sl We have $\mathfrak M\supset \mathbf H_2^t\ominus b_a\star \mathbf H_2^t$.}\small skip

    We follow the proof of Proposition 3.10 in \cite{alpay2024schur}. Let $g(q_t)=\sum_{n=0}^\infty q_t^ng_n\in \mathbf H_2^t\ominus B_a\star \mathbf H_2^t$. Then for every $k\in\mathbb N_0$,
    \[
      [g,b_a\star q_t^k]_{\ct}=0.
    \]
    Using \eqref{b-a-expansion} we can write
    \[
            [g,b_a\star q_t^k]_{\ct}=-\alpha^{\ct}g_k+\sum_{j=1}^\infty (I_2-\alpha\alpha^{\ct})\alpha^{j-1}g_{k+j}.
          \]
          Hence, 
          \begin{equation}
            \label{equa-1234}
            \begin{split}
              \alpha^{\ct}g_k              &=\sum_{j=1}^\infty (I_2-\alpha\alpha^{\ct})\alpha^{j-1}g_{k+j}\\
              &=(I_2-\alpha\alpha^{\ct})g_{k+1}+(I_2-\alpha\alpha^{\ct})\alpha g_{k+2}+(I_2-\alpha\alpha^{{\ct}2})\alpha^2 g_{k+3}+\cdots, \quad k=0,1,\ldots
              \end{split}
            \end{equation}
            We multiply on the left these equalities by $\alpha$ to obtain
            \[
\alpha         \alpha^{\ct}g_k     =\alpha(I_2-\alpha\alpha^{\ct})g_{k+1}+\alpha(I_2-\alpha\alpha^{\ct})\alpha g_{k+2}+\alpha(I_2-\alpha\alpha^{\ct})\alpha^2 g_{k+3}+\cdots, \quad k=0,1,\ldots
            \]
           For $k+1$ instead of $k$ this identity becomes
                         \[
                           \alpha         \alpha^{\ct}g_{k+1}=\alpha(I_2-\alpha\alpha^{\ct})g_{k+2}+\alpha(I_2-\alpha\alpha^{\ct})\alpha g_{k+3}+\alpha(I_2-\alpha\alpha^{\ct})\alpha g_{k+4}+\cdots, \quad k=0,1,\ldots
            \]
            Subtracting from \eqref{equa-1234}
            \[
               \alpha^{\ct}g_k       =(I_2-\alpha\alpha^{\ct})g_{k+1}+(I_2-\alpha\alpha^{\ct})\alpha g_{k+2}+(I_2-\alpha\alpha^{{\ct}})^2\alpha^2 g_{k+3}+c\dots, \quad k=0,1,\ldots
\]
            we obtain
            \[
\alpha\alpha^{\ct}g_{k+1}-\alpha^{\ct}g_k=-(I_2-\alpha\alpha^{\ct})g_{k+1}, \quad k=0,1,\ldots
\]
and so
\[
  g_{k+1}=\alpha^{\ct}g_k, \quad k=0,1,\ldots
\]
Thus $g_{k}=\alpha^{\ct k}g_0$, and
\[
g(q_t)=\sum_{k=0}^\infty q_t\alpha^{\ct k}g_0=(I_2-q_t\alpha^{\ct})^{-\star}g_0\in\mathfrak M.
\]
To conclude the proof, we write for $f\in\mathbf H_2^t$,
\[
  f(q_t)=\underbrace{f(q_t)-(1-q_t\alpha^{\ct})^{-\star}(1-\alpha\alpha^{\ct})f(\alpha)}_{\mbox{{\rm vanishes at $\alpha$}}}+\underbrace{(1-q_t\alpha^{\ct})^{-\star}(1-\alpha\alpha^{\ct})f(\alpha)}_{\in \mathfrak M}
  \]
We have thus $f(q_t)-(1-q_t\alpha^{\ct})^{-\star}(1-\alpha\alpha^{\ct})f(\alpha)\in b_a\star \mathbf H_2^t$, and the result follows.
      \end{proof}

      \begin{remark}
        \label{unfortunate}
        Unfortunately the previous result cannot be used iteratively because of the non-commutativity. If $\alpha_1$ and $\alpha_2\in\mathbb H_t$ are of small neough norm, then we indeed have $f(q_t)=b_{\alpha_1}(q_1)\star g(q_t)$ when $f(\alpha_1)=0.$ 
        On the other hand the condition
        $f(\alpha_2)=0$, i.e.
\[
  (b_{\alpha_1}(q_t)\star g(q_t))(\alpha_2)=0
  \]
  does not translate in a way conducive to computations, even when using Lemma \ref{multi-3}.
  Indeed, let $\alpha_1$ and $\alpha_2$ be in $\mathbb H_t$ and small enough. Assume that $f\in\mathbf H^t_{2,\ct}$ is such that $f(\alpha_1)=0$. Then
  \[
    f(q_t)=b_{\alpha_1}(q_t)\star g(q_t)
  \]
  where $g\in\mathbf H^t_{2,\ct}$. We now require $f(\alpha_2)=0$. Assuming $b_{\alpha_1}(\alpha_2)$ invertible, we can write
  \[
f(\alpha_2)=\left(b_{\alpha_1}(\alpha_2)\right)g(\alpha_3)
\]
where
\[
  \alpha_3=\left(b_{\alpha_1}(\alpha_2)\right)^{-1}\alpha_2\left(b_{\alpha_1}(\alpha_2)\right)g(\alpha_3)\not=\alpha_2
\]
in general. A special case occcurs when $\alpha_2=\alpha_1^{\ct}$. Then (see Proposition \ref{B-star-comm}),
\[
  \begin{split}
    b_{\alpha_1}(q_t)b_{\alpha_1^{\ct}}(q_t)&=(q_t-\alpha)\star(1-q_t\alpha^{\ct})^{-star}\star(q-\alpha^{ct})\star(1-q_t\alpha)^{-\star}\\
    &=(q_t-\alpha)\star(q-\alpha^{ct})\star(1-q_t\alpha^{\ct})^{-star}    \star(1-q_t\alpha)^{-\star}\\
    &=(q_t^2-2q_t ({\rm Re}\, \alpha)+\alpha\alpha^{\ct})(q_t(\alpha\alpha^\ct)-2q_t({\rm Re}\,\alpha)+1)^{-1}
        \end{split}
\]
and so we have interpolation on the sphere associated to $\alpha$ (see Definition~\ref{sph} for the latter).
\end{remark}

  In general one needs to proceed differently, in a global way.
  The following extension of Theorem \ref{inter-pol-1} can be seen as a finite dimensional version of the Beurling-Lax theorem in
  the present setting. First some notations. We define a matrix ${\mathbf G_{\ct}}\in\mathbf H_t^{N\times N}$ via
  \begin{equation}
    \label{opopG}
  {\mathbf G_{\ct}}=\begin{pmatrix}
  \sum_{n=0}^\infty \alpha_1^n\alpha_1^{{\ct} n}&\sum_{n=0}^\infty \alpha_1^n\alpha_2^{{\ct} n}&\cdots &\sum_{n=0}^\infty \alpha_1^n\alpha_N^{{\ct} n}\\
  \sum_{n=0}^\infty \alpha_2^n\alpha_1^{{\ct} n}&\sum_{n=0}^\infty \alpha_2^n\alpha_2^{{\ct} n}&\cdots &\sum_{n=0}^\infty \alpha_2^n\alpha_N^{{\ct} n}\\
 \vdots &\vdots && \vdots\\
  \sum_{n=0}^\infty \alpha_N^n\alpha_1^{{\ct } n}&\sum_{n=0}^\infty \alpha_N^n\alpha_2^{{\ct} n}&\cdots& \sum_{n=0}^\infty \alpha_N^n\alpha_N^{{\ct} n}\end{pmatrix}.
\end{equation}
where $\alpha_1,\alpha_2,\dots,\alpha_N\in\mathbb H_t$  are the interpolation points, assumed to be in the operator unit ball of $\mathbf H_t$, i.e.
such that
\begin{equation}
  \|\alpha_j\|_{op}<1,\quad j=1,\ldots, N.
\label{unit-call}
\end{equation}
Then, see Corollary \ref{opop}, the $\ct$-adjoints also are in the operator unit ball, and thus the power series in \eqref{opopG} converge.
\begin{theorem}
  Assume ${\mathbf G_{\ct}}$ invertible in $\mathbb H_{t}^{N\times N}$ and define
  \begin{equation}
     \label{theta_N}
     \Theta(q_t)=1-(1-q_t)\star \begin{pmatrix}(1-q_t\alpha_1^{\ct})^{-\star}&\cdots&(1-q_t\alpha_N^{\ct})^{-\star}\end{pmatrix}{\mathbf G_{\ct}}^{-1}
     \begin{pmatrix}(1-\alpha_1)^{-1}\\ \vdots \\ (1-\alpha_N)^{-1}\end{pmatrix}.
   \end{equation}
 A function $f\in\mathbf H_{2,\ct}$ vanishes at $\alpha_1,\alpha_2,\dots,\alpha_N$ if it can be written as $f=b_\alpha\star g$ where $g\in\mathbf H_{2,\ct}$
  and, moreover
  \[
    [f,f]_{\ct}=[g,g]_{\ct}
  \]
  \label{inter-pol-N}
  \end{theorem}

  Before presenting the proof we need some preliminary results. We first note that     
      \eqref{ope-zero} has the following generalization:

      \begin{lemma}
        \label{ope-zero-1}
Assume that $f\star g$ makes sense and that $f(\alpha)=0$. Then, $(f\star g)(\alpha)=0$.
\end{lemma}

\begin{proof} Let $f(q_t)=\sum_{n=0}^\infty q_t^nf_n$ and $g(q_t)=\sum_{m=0}^\infty q_t^mg_m$ where the coefficients $f_n$ and $g_m$ beong to $\mathbb H_t$. The claim follows from
the formula
\[
(  f\star g)q_t=\sum_{m=0}^\infty q_t^mf(q_t)g_m.
\]
\end{proof}

\begin{proposition}
  \label{power-exp}
Assume that the points $\alpha_1,\ldots, \alpha_N$ are in the open unit operator ball of $\mathbb H_t$, and let
  \begin{eqnarray*}
           A&=&{\rm diag}(\alpha_1^{\ct},\ldots, \alpha_N^{\ct})\\
                C&=&\begin{pmatrix}1&1&\cdots &1\end{pmatrix}.
                                                \end{eqnarray*}
The matrix  ${\mathbf G_{\ct}}$ is the unique solution of the Stein equation
  \begin{equation}
    {\mathbf G_{\ct}}-A^{\ct}{\mathbf G_{\ct}} A=C^{\ct}C.
    \label{stein-eq}
  \end{equation}
  and is given by the formula
  \begin{equation}
    \label{formula-G}
    {\mathbf G_{\ct}}=\sum_{n=0}^\infty A^{\ct n}C^\ct CA^n.
    \end{equation}
\end{proposition}

\begin{proof}
  From the hypothesis on the interpolation points we have $\|A\|_{op}<1$ and the series on the right hand side of \eqref{formula-G} converges.
  On the other hand, iterating \eqref{stein-eq} $(M-1)$-times we get
  \[
{\mathbf G_{\ct}}=\sum_{n=1}^{M} A^{\ct n}C^\ct CA^n +A^{\ct(M+1)}{\mathbf G_{\ct}} A^{(M+1)}
\]
The result follows by letting $M\rightarrow\infty$.
  \end{proof}

      \begin{proof}[Proof of Theorem \ref{inter-pol-N}]         We will proceed in steps and follow arguments from \cite{alpay2024schur}.\\

   STEP 1: {\sl It holds that $\Theta(\alpha_j)=0, \quad j=1,\ldots N$}\smallskip

   Without loss of generality we consider the case $j=1$. 
   We have that
   \[
     \begin{pmatrix}(1-q_t\alpha_1^{\ct})^{-\star}&\cdots&(1-q_t\alpha_N^{\ct})^{-\star}\end{pmatrix}(\alpha_1)
     \]
     is the first row of the matrix $G$ and that
     \[
       (1-q_t)\star \begin{pmatrix}(1-q_t\alpha_1^{\ct})^{-\star}&\cdots&(1-q_t\alpha_N^{\ct})^{-\star}\end{pmatrix}=(1-\alpha_1)\begin{pmatrix}g_{11}&g_{12}&\ldots
         &g_{1N}\end{pmatrix}.
       \]
       Thus
       \[
         \Theta(\alpha_1)=1-(1-\alpha_1)\begin{pmatrix}g_{11}&g_{12}&\ldots         &g_{1N}\end{pmatrix}{\mathbf G_{\ct}}^{-1}
          \begin{pmatrix}(1-\alpha_1)^{-1}\\ \vdots \\ (1-\alpha_N)^{-1}\end{pmatrix}=0
         \]
     
         STEP 2:  {\sl It holds that
                           \begin{eqnarray}
\label{theta-!!!}
                             \Theta(q_t)&=&D+q_tC\star(I_{\mathbb H_1^N}-q_tA)^{-\star}B
                             \end{eqnarray}
where (note that $A$ and $C$ were defined earlier in Proposition \ref{power-exp})
         
         \begin{eqnarray}
           A&=&{\rm diag}(\alpha_1^{\ct},\ldots, \alpha_N^{\ct})
           \label{t-1}\\
                B&=&(I_{\mathbb H_t^N}-A){\mathbf G_{\ct}}^{-1}(I_{\mathbb H_t^N}-A^{\ct})^{-1}C^{\ct}.\\
                     C&=&\begin{pmatrix}1&1&\cdots &1\end{pmatrix}\\
           D&=&1-C{\mathbf G_{\ct}}^{-1}(I_{\mathbb H_t^N}-A^{\ct})^{-1}C^{\ct}.
                \label{t-4}
         \end{eqnarray}
}

The computations spresented in the next step  are used later, see Example \ref{ex-1}, but we need at this stage equality \eqref{t-900} in Step 4.\\

STEP 3: {\sl The following three equalities hold:}

\begin{eqnarray}
  A^\ct{\mathbf G_{\ct}} A+C^\ct C&=&{\mathbf G_{\ct}},
                              \label{t-9}\\
  \label{t-890}
  B^{\ct}{\mathbf G_{\ct}} A+ D^\ct C&=&0,\\
  B^{\ct}{\mathbf G_{\ct}} B +D^\ct D&=&I.
                                 \label{t-900}
    \end{eqnarray}

 Indeed, the first equality is the Stein equation \eqref{stein-eq}. We now check \eqref{t-890} and write
    \[
      \begin{split}
        B^{\ct}{\mathbf G_{\ct}} A+ D^\ct C&=C(I_N-A)^{-1}{\mathbf G_{\ct}}^{-1}(I_N-A^\ct)
        {\mathbf G_{\ct}} A+(1-C(I-A)^{-1}{\mathbf G_{\ct}}^{-1}C^\ct)C \\
        &=C(I_N-A)^{-1}{\mathbf G_{\ct}}^{-1}\left\{(I_N-A^\ct)
          {\mathbf G_{\ct}} A+{\mathbf G_{\ct}}(I_N-A)-C^\ct C\right\}\\
        &=C(I_N-A)^{-1}{\mathbf G_{\ct}}^{-1}\left\{-A^\ct{\mathbf G_{\ct}} A+{\mathbf G_{\ct}}-C^\ct C\right\}\\
        &=0
      \end{split}
      \]
      thanks to the Stein equation \eqref{stein-eq}. Finally to check \eqref{t-900} we write
      \[
        \begin{split}
          B^{\ct}{\mathbf G_{\ct}} B +D^\ct D&=C(I_N-A)^{-1}{\mathbf G_{\ct}}^{-1}
          (I_N-A^\ct){\mathbf G_{\ct}} (I_{(\mathbf H_2^t)^N}-A){\mathbf G_{\ct}}^{-1}(I_{(\mathbf H_2^t)^N}-A^{\ct})^{-1}C^{\ct}+\\
            &\hspace{5mm}+(1-C(I-A)^{-1}{\mathbf G_{\ct}}^{-1}C^\ct)(1-C{\mathbf G_{\ct}}^{-1}(I_{(\mathbf H_2^t)^N}-A^{\ct})^{-1}C^{\ct})\\
            &=1+C(I_N-A)^{-1}{\mathbf G_{\ct}}^{-1}\left\{ (I_N-A^\ct){\mathbf G_{\ct}} (I_{(\mathbf H_2^t)^N}-A)+  \right.\\
            &\hspace{5mm}\left.  -(I_N-A^\ct){\mathbf G_{\ct}}-{\mathbf G_{\ct}}(I_N-A)+C^\ct C \right\}{\mathbf G_{\ct}}^{-1}(I_{(\mathbf H_2^t)^N}-A^{\ct})^{-1}C^{\ct}\\
            &=1+C(I_N-A)^{-1}{\mathbf G_{\ct}}^{-1}\left\{-{\mathbf G_{\ct}}+ A^\ct{\mathbf G_{\ct}} A+C^\ct C\right\}\times\\
          &\hspace{5mm}\times{\mathbf G_{\ct}}^{-1}(I_{(\mathbf H_2^t)^N}-A^{\ct})^{-1}C^{\ct}\\
          &=0,
        \end{split}
        \]
      here too thanks to the Stein equation \eqref{stein-eq}.\\

This follows from \eqref{mbounded} in the second step of Theorem \ref{thmbabounded} applied to each of the $\alpha_j$.\\

STEP 4: {\sl  The operator $M_\Theta$ of $\star$-multiplication on the left by $\Theta$ is an isometry from $\mathbf H_{2,\ct}$ into itself.}\smallskip

The fact that $M_\Theta$ is bounded follows from \eqref{mbounded} in the second step of Theorem \ref{thmbabounded} applied to each of the $\alpha_j$.
We now prove that $M_\Theta$ is isometric on polynomials. The continuity of $M_\Theta$ implies then the isometry property on the whole Hardy space.
We   first remark that $\Theta=\sum_{n=0}^\infty q_t^n\theta_n$ with
         \begin{equation}
           \theta_n= \begin{cases}\,D,\quad\hspace{10mm} n=0,\\
             \,\, CA^{n-1}B,\,\, n=1,2,\ldots
           \end{cases}
         \end{equation}
         So
         \begin{equation}
           \label{u-u}
           \begin{split}
             \theta_0^\ct\theta_0+\sum_{n=1}^\infty \theta_n^\ct\theta_n&=D^\ct D+B^\ct\left(\sum_{n=0}^\infty A^{\ct n}C^\ct CA^n\right)B\\
             &=D^\ct B+B^\ct {\mathbf G_{\ct}}B\\
             &=1,
           \end{split}
         \end{equation}
         where we have used \eqref{t-900} and the formula \eqref{formula-G} for ${\mathbf G_{\ct}}$, and for $k>0$,
         \begin{equation}
\label{u-u}
           \begin{split}
             \theta_0^\ct\theta_k+\sum_{n=1}^\infty \theta_n^\ct\theta_{n+k}&=D^\ct CA^{k-1}B+B^\ct\left(\sum_{n=0}^\infty A^{\ct n}C^\ct CA^n\right)A^{k-1}B\\
             &=(D^\ct C+B^\ct {\mathbf G_{\ct}})A^{k-1}B\\
             &=0,
           \end{split}
         \end{equation}
         where we have used \eqref{t-890} and here too the formula \eqref{formula-G} for ${\mathbf G_{\ct}}$.\\

         We follow STEP 1 in the proof of Theorem \ref{thmbabounded} and compute $[M_\Theta q_t^u\gamma,M_\Theta q_tv\delta]_{\ct}$ for $\gamma=\delta=1$.
         We may, and shall, assume that $u\ge v$ and set $k=u-v$. We have
\[
         \begin{split}
           [M_\Theta q_t^u,M_\Theta q_t^v]&=[M_\Theta 1,M_\Theta q_t^v]_{\ct}\\
           &=[\sum_{n=0}^\infty q_t^n\theta_n,\sum_{n=0}^\infty q_t^{n+k}\theta_n]_{\ct}\\
           &= \theta_0^\ct\theta_k+\sum_{n=1}^\infty \theta_n^\ct\theta_{n+k}\\
           &=\delta_{k,0}
           \end{split}
         \]
         by the previous step, and hence the result by linearity and continuity.\\

         STEP 5: {\sl Let $\mathfrak M$ be the real linear span of the functions $(1-q_t\alpha_j)^{-\star}b$, where $j=1,\ldots, N$ and $b\in\mathbb H_t$.
           Then
           \begin{equation}
M_\Theta \mathbf H_{2,\ct}^t\subset             \mathfrak M^{\perp_{\ct}},
             \end{equation}
             where $\perp_\ct$ is the Krein orthogonal with respect to $[\cdot,\cdot]_{\ct}$}\smallskip

           This follows from Step 1.
           
         \end{proof}

         The above result is a finite dimensional Beurling-Lax type theorem in the present setting.
 
 \section{Hardy space and interpolation for the second conjugate $[*]$}
 \setcounter{equation}{0}
 \label{Hardy_q_t}
\subsection{The Hardy space in one variable: Case $2$}
 We recall from Subsection~\ref{Hardy-prelude} that the Hardy space $\mathbf H_2(\mathbb D)$ is the reproducing kernel Hilbert space of functions analytic in the open unit disk and with reproducing $\frac{1}{1-z\overline{w}}$. 
We also recall that the counterpart of the Hardy space is now the following real Hilbert space $\mathbf H^t_2$ the space of power series of the form
\begin{equation*}
f(q_t)=\sum_{n=0}^\infty q_t^n\alpha_n
\end{equation*}
as in~\eqref{norm-12345}.

 We now focus on $\mathbf H^t_2$ with the real $[*]$ symmetric form~\eqref{function-inner-2}, namely 
\begin{equation*}
[f,g]_{[*]}=\sum_{n=0}^\infty [g_n,f_n]_{[*]}
\end{equation*}
where $g(q_t)=\sum_{n=0}^\infty q_t^ng_n$ is another element of $\mathbf H_2^t$ and we recall that Theorem~\ref{H_2-Krein} for this case states:
\begin{theorem}
  The space $\mathbf H_2^t$ endowed with the form \eqref{function-inner-2}  is a Krein space. This space was denote by  $\mathbf H_{2,[*]}^t$ in Definition~\ref{Hardy-Krein}.
\end{theorem}


    \setcounter{equation}{0}
\begin{proposition}
  The space $\mathbf H_{2,[*]}^t$ endowed with the real symmetric form \eqref{function-inner-2} is a reproducing kernel Krein space right $\mathcal H_2^t$-module with reproducing kernel
  with reproducing kernel
    \begin{equation}
k(q_t,p_t)=\sum_{n=0}^\infty q_t^n(p_t^n)^{[*]}.
\end{equation}
\end{proposition}

\begin{proof}
  We write
  \[
    \begin{split}
      [f(\cdot),k(\cdot,p_t)b]_{[*]}&=[\sum_{n=0}^\infty q_t^nf_n,\sum_{n=0}^\infty q_t^n(p_t^n)^{[*]}b]_{[*]}\\
      &=\sum_{n=0}^\infty[f_n,(p_t^n)^{[*]}b]_{[*]}\\
      &=\sum_{n=0}^\infty{\rm Tr}\,(b^{[*]}p_t^nf_n)\\
      &={\rm Tr}\, b^{[*]}f(p_t))\\
      &=[f(p_t),b]_{[*]}.
    \end{split}
    \]
\end{proof}

We see that $\mathbf H_{2,[*]}^t$ is an $\mathbb H_t$-module. For $t=-1$ the second adjoint gives a new theory for quaternions.

\subsection{Star product and Blaschke factor}
We remind the reader of the precursor to this theory, described in Subsection~\ref{interpolation-prelude}.
Since $\alpha\alpha^{[*]}$ is not (in general) a scalar matrix, one cannot define Blaschke factors using the formula
$(q_t-\alpha)\star (I_2-q_t\alpha^{[*]})^{-\star}$, but one could think of defining a Blaschke factor here by 

\begin{equation}
  \label{bl-2}
C_{\alpha}(q_t)=(1-\alpha\alpha^{[*]})^{-1/2}\star (q_t-\alpha)\star(1-q_t\alpha^{[*]})^{-\star}(1-\alpha^{[*]}\alpha)^{1/2};
\end{equation}
See \eqref{formule-e1-e2}.
But this definition does not allow for a counterpart of Theorem \ref{inter-pol-1}. The appropriate definition is
\[
B_\alpha(q_t)=(q_t-\alpha)\star(1-q_t\Gamma_\alpha \alpha^{[*]}\Gamma_\alpha^{-1})^{-\star}L_{\alpha}^{1/2}
\]
which is introduced (for regular adjoint) in \cite{alpay2024schur} in a different context, and takes its inspiration
from the formulas for a Blaschke factor in \cite[Section 4]{MR93b:47027}.
For $t\not=-1$ $\Gamma_\alpha$ is not self-adjoint, and need not be invertible not have a squareroot. But we have:
  \begin{lemma}
    $\Gamma_\alpha$ is $[*]$ self-adjoint. When
  \[
\sum_{n=1}^\infty\|\alpha\|_{op}^{2n}=\frac{\|\alpha\|_{op}^2}{1-\|\alpha\|_{op}^2}<1
\]
$\Gamma_\alpha$ is invertible and there exists an element of $\mathbb H_t$, denoted $\Gamma_\alpha^{1/2}$, which is $[*]$ self-adjoint and
whose square is $\Gamma_\alpha$.
  \end{lemma}

  \begin{proof}
    We write
 \[
  \Gamma_\alpha=1+\underbrace{\sum_{n=1}^\infty \alpha^n\alpha^{[*]n}}_{\e}
  \]
  For $\|\alpha\|_{op}<1$ we have
  \begin{equation}
    \label{786}
\|\e\|_{op}<\frac{\|\alpha\|_{op}^2}{1-\|\alpha\|_{op}^2}
  \end{equation}
  we can define $\Gamma_\alpha^{1/2}$ via
  \[
  \Gamma_\alpha^{1/2}=1+\sum_{n=1}^\infty\gamma_n\alpha^n,
  \]
where the numbers $\gamma_n$ are defined via  
\begin{equation}
  \label{sqrt-power}
\sqrt{1+t}=1+\sum_{n=1}^\infty \gamma_nt^n,\quad t\in\mathbb C\,\,{\rm such \,\, that}\,\,\,|t|<1.
\end{equation}
  \end{proof}

  Still following \cite{alpay2024schur} and \cite{MR93b:47027} we define
\[
  L_\alpha=\Gamma_\alpha-\Gamma_\alpha \alpha^{[*]}\Gamma_\alpha^{-1}\alpha\Gamma_\alpha
  \]
  for $\alpha$ such that $\Gamma_\alpha^{-1}$ exists.
  
\begin{lemma}  
  $L_\alpha$ satisfies $L_\alpha=L_\alpha^{[*]}$ and the set of $\alpha$ for which $L_\alpha$ is invertible and for which there exists $K_\alpha$  satisfying $K_\alpha=K_a^{[*]}$ and $K_\alpha^2=L_\alpha$ is open and not empty.
Furthermore we have
\begin{equation}
L_\alpha^{-1}=\alpha^{[*]}\alpha+\Gamma_\alpha^{-1}
\label{567}
\end{equation}
\label{lemma-567}
\end{lemma}

\begin{proof}
  Let as above $\Gamma_\alpha=1+\e$ with $\alpha$ be such that $\|\e\|_{op}<1$. Then, the series
  \[
\sum_{n=0}^\infty (-1)^n\e^n
    \]
converges in $\mathbb H_t$ and (since the operator norm is submultiplicative)
\[
\begin{split}
  \|\Gamma_\alpha^{-1}\|_{op}&\le\frac{1}{1-\|\e\|_{op}}\\
  &\le \frac{1}{1-\frac{\|\alpha\|_{op}^2}{1-\|\alpha\|^2_{op}}}\\
    &=\frac{1-\|\alpha\|^2_{op}}{1-2\|\alpha\|^2_{op}}<1
  \end{split}
\]
for $\alpha$ small enough. From
\[
L_\alpha=1+\underbrace{\Gamma_\alpha-1-\Gamma_\alpha \alpha^{[*]}\Gamma_\alpha^{-1}\alpha\Gamma_\alpha}_{\e_1}
\]
and using \eqref{786}, we have
\[
  \begin{split}
    \|\e_1\|_{op}&\le\frac{\|\alpha\|_{op}^2}{1-\|\alpha\|_{op}^2}+\|\Gamma_\alpha\|_{op}^2\|\Gamma_\alpha^{-1}\|_{op}\|\alpha\|^2_{op}\\
    &\le \frac{\|\alpha\|_{op}^2}{1-\|\alpha\|_{op}^2}+\left(1+\frac{\|\alpha\|^2_{op}}{1-\|\alpha\|^2_{op}}\right)^2\frac{1}{1-\frac{\|\alpha\|_{op}^2}{
        1-\|\alpha\|_{op}^2}} \|\alpha\|^2_{op}\\
    &\le  \frac{\|\alpha\|_{op}^2}{1-\|\alpha\|_{op}^2}+\left(\frac{1}{1-\|\alpha\|^2_{op}}\right)^2\frac{1-\|\alpha\|_{op}^2}{1-2
        \|\alpha\|_{op}^2}\|\alpha\|^2_{op}\\
    &<1
\end{split}
  \]
for $\alpha$ small enough. For such $\alpha$, $L_\alpha$ is also invertible and has a squareroot which is $[*]$ symmetric. The argument is the same
as in the preceding lemma. We now check \eqref{567}. We have
\[
\begin{split}
  (\alpha^{[*]}\alpha+\Gamma_{\alpha}^{-1})L_{\alpha}&={\alpha}^{[*]}{\alpha}(\Gamma_\alpha-\Gamma_\alpha \alpha^{[*]}\Gamma_\alpha^{-1}\alpha\Gamma_\alpha)+1-\alpha^{[*]}\Gamma_\alpha^{-1}\alpha\Gamma_\alpha\\
  &=1+\alpha^{[*]}\alpha\Gamma_\alpha+\alpha^{[*]}(1-\Gamma_\alpha)\Gamma_\alpha^{-1}\alpha\Gamma_\alpha-\alpha^{[*]}\Gamma_\alpha^{-1}\alpha\Gamma_\alpha\\
  &=1.
\end{split}
\]

To conclude the proof we note that $M_{B_\alpha}$ is continuous for $\alpha$ small enough, as is seen from the formula giving $B_\alpha$, and so the
isometry property extends to all of $H_{2,[*]}^t$ (recall that ina Krein space, a densely defined unitary map need not be continuous).
\end{proof}

\begin{theorem}
The operator of $\star$-multiplication by $B_\alpha$ on the left is isometric from the Krein space $\mathbf H_{2,[*]}^t$ into itself.
\end{theorem}

\begin{proof}
  It is enough to check on monomials and the proofs of Steps 1, 2, and 3 are similar to the $\circledast$ case. We refer to the proof of Theorem~\ref{thmbabounded} for more details.\\

  STEP 1: {\sl It holds that}

  \begin{equation}
    \label{aveiro-2024}
      B_\alpha(q_t)=-\alpha K_\alpha+\sum_{n=1}^\infty q_t^n\alpha^{[*](n-1)}\Gamma_\alpha^{-1}K_\alpha
    \end{equation}
    We set $q_t=x\in\mathbb R$ (small enough) in the definition of $B_\alpha$ and write:
    \[
      \begin{split}
        B_\alpha(x)&=(x-\alpha)\left(K_\alpha+\sum_{n=0}^\infty x^n\Gamma_\alpha\alpha^{[*]n}\Gamma_\alpha^{-1}K_\alpha\right)\\
      &=-\alpha K_\alpha+\sum_{n=0}^\infty x^{n+1}\Gamma_\alpha\alpha^{[*]n}\Gamma_\alpha^{-1}K_\alpha-
      \underbrace{\sum_{n=0}^\infty x^{n+1}\alpha\Gamma_\alpha\alpha^{[*]n}\Gamma_\alpha^{-1}K_\alpha}_{n=0\,\,{\mbox{corresponds to}}\,\, xK_\alpha}+xK_\alpha\\
      &=-\alpha K_\alpha+\sum_{n=1}^\infty x^{n}\Gamma_\alpha\alpha^{[*](n-1)}\Gamma_\alpha^{-1}K_\alpha-
      \sum_{n=1}^\infty x^{n}(\alpha\Gamma_\alpha\alpha^{[*]})\alpha^{[*](n-1)}\Gamma_\alpha^{-1}K_\alpha\\
      &=-\alpha K_\alpha+\sum_{n=1}^\infty x^{n}\alpha^{[*](n-1)}\Gamma_\alpha^{-1}K_\alpha
      \end{split}
      \]
      since $\Gamma_\alpha$ satisfies the Stein equation
      \begin{equation}
        \label{stein-4}
        \Gamma_\alpha-\alpha\Gamma_\alpha\alpha^{[*]}=1.
        \end{equation}
    
    STEP 2: {\sl We have}
      \[
        [B_\alpha\star q_t^mc,B_\alpha\star q_t^md]_{[*]}=[c,d]_{[*]}
      \]
      We have
      \[ [B_\alpha\star q_t^mc,B_\alpha\star q_t^md]_{[*]}=[\Delta c,d]_{[*]} ,\]
      where
      \[
        \begin{split}
          \Delta&=K_\alpha\alpha^{[*]}\alpha K_\alpha+K_\alpha\Gamma_\alpha^{-1}
          \left(\underbrace{\sum_{n=1}^\infty\alpha^{(n-1)}\alpha^{[*](n-1)}}_{\Gamma_\alpha}\right)\Gamma_{\alpha}^{-1}
          K_\alpha\\
          &=K_\alpha\alpha^{[*]}\alpha K_\alpha+K_\alpha\Gamma_\alpha^{-1}\Gamma_\alpha\Gamma_{\alpha}^{-1}          K_\alpha\\
          &=K_\alpha\left(\alpha^{[*]}\alpha +\Gamma_\alpha^{-1}\right)K_\alpha\\
          &=K_\alpha L_\alpha^{-1}K_\alpha \quad ({\rm by}\,\,\eqref{567}\\
          &=1
          \end{split}
        \]
        where we have used Lemmma \ref{lemma-567}. Note that $K_\alpha L_\alpha^{-1}K_\alpha=1$ and $K_\alpha^2=L_\alpha$ are indeed equivalent since both
        $K_\alpha$ and $L_\alpha$ are invertible.\\
        
          STEP 3: {\sl We have}
      \[
        [B_\alpha\star q_t^uc,B_\alpha\star q_t^vd]_{[*]}=0,\quad u\not=v.
        \]

  Without loss of generality we will assume $u>v$,      In a way similar to the proof of Step 2 we have with $h=u-v$:
  \[
    \begin{split}
      [B_\alpha\star q_t^{(u-v)}c,B_\alpha\star d]_{[*]}&=\\
      &\hspace{-3cm}=
      [-q_t^h\alpha K_\alpha c +\sum_{n=1}^\infty q_t^{h+n}\alpha^{[*](n-1)}\Gamma_\alpha^{-1}K_\alpha c,\sum_{n=h}^\infty q_t^n\alpha^{[*](n-1)}\Gamma_\alpha^{-1}
      K_\alpha d]_{[*]}\\
&\hspace{-3cm}=[-q_t^h\alpha K_\alpha c+\sum_{n=1}^\infty q_t^{h+n}\alpha^{[*](n-1)}\Gamma_\alpha^{-1}K_\alpha c, 
\alpha^{[*](h-1)}\Gamma_\alpha^{-1}K_\alpha d+\sum_{n=0}^\infty q_t^{h+n}\alpha^{[*](h+n-1)}\Gamma_\alpha^{-1}K_\alpha d]_{[*]}\\     
      &\hspace{-3cm}=      [\Delta_hc,d]_{[*]},
\end{split}
        \]
        where 
        \[
\begin{split}
  \Delta_h&= -K_\alpha\Gamma_\alpha^{-1}\alpha^h K_\alpha+\sum_{n=0}^\infty K_\alpha\Gamma_\alpha^{-1} \underbrace{\alpha^{h+n}\alpha^{[*]n}}_{\alpha^h\Gamma_\alpha}\Gamma^{-1}_\alpha K_\alpha\\
  &= -K_\alpha\Gamma_\alpha^{-1}\alpha^h K_\alpha+\sum_{n=0}^\infty K_\alpha\Gamma_\alpha^{-1}\alpha^{h}\Gamma_\alpha\Gamma^{-1}_\alpha K_\alpha\\
  &=0.
          \end{split}
\]
      \end{proof}

\begin{theorem}
$f\in\mathbf H^t_{2,[*]}$ vanishes at $\alpha$ if and only if $f=B_\alpha\star g$ where $g\in\mathbf H^t_{2,[*]}$ satisfies $[f,f]_{[*]}=[g,g]_{[*]}$
  \end{theorem}

  \begin{proof}
    One direction is clear, thanks to \eqref{ope-zero} (see Step 1 in the proof of Theorem \ref{inter-pol-1}. To prove the converse
    we proceed in a number of steps, and first define the operator
    \begin{equation}
      M^aq_t^nb=q_t^nab,\quad a,b\in\mathbb H_t.
    \end{equation}

    STEP 1: {\sl $M^a$ is bounded in $\mathbf H_{2,[*]}^t$ and it holds that $(M^a)^{[*]}q_t^nb=q_t^na^{[*]}b$}\smallskip

    $M^a$ is defined everywhere and closed (as is seen using the reproducing kernel property), and so is bounded in the Krein space
    $\mathbf H_{2,[*]}^t$. To compute the adjoint we write
    \[
      \begin{split}
        [M_a^{[*]}(q_t^nb),q_t^mc]_{[*]}&=[q_t^nb,M^a(q_t^mc)]_{[*]}\\
        &=[q_t^nb,q_t^mac)]_{[*]}\\
        &=\delta_{m,n} {\rm Tr} ((ac)^{[*]}b)\\
&=        \delta_{m,n} {\rm Tr} (c^{[*]}a^{[*]}b)\\
&=\delta_{m,n}[q_t^na^{[*]}b,q_t^mb]_{[*]}.
\end{split}
\]
        
STEP 2: {\sl The operator $M_{q_t}-M^a$ has closed range and}
  \[
    H^t_{2,[*]}={\rm ran}\,(M_{q_t}-M^a)\stackrel{[+]}{[*]}{\rm span}\left\{(1-q_ta^{[*]})^{-\star}b, b\in\mathbb H_1\right\}
  \]

  The closed range follows from the reproducing kernel property and the previous step. Let now $f$ orthogonal in the $[*]$ inner product to the range.
  Then
  \[
    (    M_{q_t}^{[*]}-(M^a)^{[*]})f=0
  \]
  With $f(q_t)=\sum_{n=0}^\infty q_t^nf_n$ and since, as is readily verified
  \[
    M_{q_t}^{[*]}q_t^nb\begin{cases}\,\,\, 0,\quad\hspace{5.6mm} n=0,\\
     \, q_t^{n-1}b,\quad n=1,2,\ldots
    \end{cases}
  \]
  we get
  \[
    f_n=a^{[*]}f_{n+1},\quad n=0,1,\ldots
  \]
  so that
  \[
    f(q_t)=(1-q_ta^{[*]})^{-\star}f_0
\]

STEP 3: {\sl We conclude the proof.}\smallskip

Let $f\in\mathbf H_{[2,[*]}^t$ be such that $f(a)=0$. By the previous step, $f\in{\rm ran}(M_{q_t}-M^a)$ so that
\[
  f(q_t)=(q_t-a)\star g(q_t)=  B_\alpha(q_t)\star   (1-q_t\Gamma_\alpha \alpha^{[*]}\Gamma_\alpha^{-1})^{\star}L_{\alpha}^{-1/2}
\]
which concludes the proof.

\end{proof}

The same conclusions as in Remark \ref{unfortunate} hold. One can define the couterpart of the matrix ${\mathbf G_{\ct}}$
and Problem \ref{inter-pol-N}, namely
  \begin{equation}
    \label{opopG}
  {\mathbf G_{[*]}}=\begin{pmatrix}
  \sum_{n=0}^\infty \alpha_1^n\alpha_1^{{{[*]}} n}&\sum_{n=0}^\infty \alpha_1^n\alpha_2^{{{[*]}} n}&\cdots &\sum_{n=0}^\infty \alpha_1^n\alpha_N^{{{[*]}} n}\\
  \sum_{n=0}^\infty \alpha_2^n\alpha_1^{{{[*]}} n}&\sum_{n=0}^\infty \alpha_2^n\alpha_2^{{{[*]}} n}&\cdots &\sum_{n=0}^\infty \alpha_2^n\alpha_N^{{{[*]}} n}\\
 \vdots &\vdots && \vdots\\
  \sum_{n=0}^\infty \alpha_N^n\alpha_1^{{{[*]} } n}&\sum_{n=0}^\infty \alpha_N^n\alpha_2^{{{[*]}} n}&\cdots& \sum_{n=0}^\infty \alpha_N^n\alpha_N^{{{[*]}} n}\end{pmatrix}.
\end{equation}
where $\alpha_1,\alpha_2,\dots,\alpha_N\in\mathbb H_t$  are the interpolation points, assumed to be in the operator unit ball of $\mathbf H_t$, i.e.
such that
\begin{equation}
  \|\alpha_j\|_{op}<1,\quad j=1,\ldots, N.
\label{unit-call}
\end{equation}


 \section{Arveson space, Blaschke factors and related topics}
 \setcounter{equation}{0}
\label{Averson_Blaschke}


\subsection{Reproducing kernel Krein spaces of $V_t$-Fueter series}
\label{Krein_q_t}

 We are interested in studying power series of the form
 \begin{equation}
   \label{rkk}
  f(x)= \sum_{\alpha\in\mathbb N_0^3}\mu_t^\alpha f_\alpha
 \end{equation}
 where the $f_\alpha\in\mathbb H_t$. Since the $\mu_t$ are not defined at $x=0\in\mathbb R^4$ one first needs to make precise where such power series
 will be considered. For estimates we will need a submultiplicative norm, and take the operator norm \eqref{norm-op}. We define
\begin{equation}
  O_{r,\rho}=\left\{x\in\mathbb R^4\,;\, r<|x_1^2+t(x_2^2+x_3^2)|\quad{\rm and}\quad  |x_u|<\rho,\,\, u=0,1,2,3\right\} 
  \end{equation}

  \begin{lemma}
In $O_{r,R,\rho}$ we have
\begin{equation}
\|\mu_t^\alpha(x)\|_{op}\le M_{r,R,\rho}^{|\alpha|}
\end{equation}
with
\begin{equation}
  M_{r,R,\rho}=\rho\left(1+\frac{3\rho^2}{r}\right)
  \end{equation}
\end{lemma}
\begin{proof} In a similar way to our work~\cite{adv_prim}, we have
  \[
    \begin{split}
      \|\mu_t^{\alpha}(x)\|_{op}&=|x_1|^{\alpha_1}|x_2|^{\alpha_2}|x_3|^{\alpha_3}\left\|\left(1+\frac{x_0}{\vec{q_t}}\right)^{|\alpha|}\right\|_{op}\\
        &\le \rho^\alpha\left\|1+\frac{x_0}{\vec{q_t}}\right\|_{op}^{|\alpha|}\\
        &=\rho^\alpha\left\|1-\frac{x_0\vec{q_t}}{x_1^2+t(x_2^2+x_3^2)}\right\|^{|\alpha|}_{op}\\
        &\le\rho^\alpha\left(1+\left\|\frac{x_0\vec{q_t}}{x_1^2+t(x_2^2+x_3^2)}\right\|_{op}\right)^{|\alpha|}\\
        &\le \rho^{|\alpha|}\left(1+\frac{3\rho^2}{r}\right)^{|\alpha|}
          \end{split}
    \]
\end{proof}

One could also consider
\begin{equation}
  O_{r,\rho,R}=\left\{x\in\mathbb R^4\,;\, r<|x_1^2+t(x_2^2+x_3^2)|<R\quad{\rm and}\quad  |x_u|<\rho,\,\, u=0,1,2,3\right\} 
  \end{equation}

The following is the counterpart of \cite[Proposition 5.1]{adv_prim}.
  
  \begin{proposition}
    Let $(c_\alpha)_{\alpha\in\mathbb N_0^3}$ be a family of positive fnumbers such that
    \[
\sum_{\substack{\alpha\in\mathbb N_0^3\\ c_\alpha>0}}\frac{M_{r_1,R_1,\rho_1}^{|\alpha|}}{c_\alpha}<\infty
\]
for all $r_1,\rho_1$ such that
\[
 0< r<r_1\quad and\quad 0<\rho_1<\rho.
\]
Then the function
 \begin{equation}
k_{\bf c}(x,y)=      \sum_{\substack{\alpha\in\mathbb N_0^3\\ c_\alpha\not=0}}\frac{\mu_t^\alpha(x)(\mu_t^\alpha(y))^{{\ct}}}{c_\alpha}
  \end{equation}
  converges in $O_{r,\rho}$, and is the reproducing kernel of the reproducing kernel Krein space of power series of the form \eqref{rkk} endowed with the symmetric form $[\cdot,\cdot]_{\ct}$, namely
    \begin{equation}
    \mathfrak H(k_{\bf c})
    =\left\{f=\sum_{\substack{\alpha\in\mathbb N_0^3\\ c_\alpha\not=0}}\mu_t^\alpha f_\alpha, \,\,\, f_\alpha\in\mathcal H_2^t\,\,|\,\,
      \sum_{\substack{\alpha\in\mathbb N_0^3\\ c_\alpha\not=0}}
      c_\alpha\|f_\alpha\|_{op}^2<\infty\right\}
    \end{equation}
\end{proposition}

\begin{proof}
We leave the proof to the reader, as it is a straightforward generalization of our work~\cite{adv_prim}.
\end{proof}

    \begin{proposition}
      Elements of $ \mathfrak H(k_{\bf c})$ are $V_{t}$-regular in
      \begin{equation}
O(k_{\bf c})=\left\{x\in\mathbb R^4\setminus\left\{(0,0,0,0)\right\}\, |\, \sum_{\substack{\alpha\in\mathbb N_0^3\\ c_\alpha\not=0}}\frac{\|\mu_t^\alpha\|_{op}^2}{c_\alpha}<\infty\right\}
        \end{equation}
      \end{proposition}

\begin{proof}
See our work~\cite{adv_prim} for the $t=-1$ case, this follows from Theorem~\ref{KerV_t}.
\end{proof}
  Similar considerations hold for power series of the form
         \begin{equation}
      \sum_{\substack{\alpha\in\mathbb N_0^3\\ c_\alpha\not=0}}\frac{\mu_t^\alpha(x)(\mu_t^\alpha(y))^{[*]}}{c_\alpha}
  \end{equation}

  The case $c_\alpha=\alpha!$ will correspond to the counterpart of the Arveson space.
  \subsection{Blaschke factor ($\mu_t$ variables)}
  We now denote the vector $\mu_t$ of $V_t$-Fueter variables by:
  \begin{equation}
\mu_t=\begin{pmatrix}\mu_{1,t}&\mu_{2,t}&\mu_{3,t}\end{pmatrix},
\end{equation}
where $\mu_{l,t}$ are defined in Equation~\eqref{mu_t}.\\
We first remark that for $a\in\mathbb  R^4$ such that
\begin{equation}
  \label{newcond}
\sum_{j=1}^3\|\mu_{j,t}(a)\|<1
\end{equation}
one can define the squareroots
\[
\sqrt{1-\mu_t(a)\mu_t(a)^\ct}\quad {\rm and}\quad (I_{\mathbb H_t^3}-\mu_t(a)^\ct\mu_t(a))^{-1/2},
\]
where 
\begin{equation*}
\mu_t(a)^\ct = \begin{pmatrix}\mu_{1,t}^{\ct}&\mu_{2,t}^{\ct}&\mu_{3,t}^{\ct}\end{pmatrix}^T,
\end{equation*}
via the power series \eqref{sqrt-power}. Since
\[
\mu_t(a)(\mu_t(a)^{\ct}\mu_t(a))^n=(\mu_t(a)\mu_t(a)^{\ct})^n\mu_t(a),\quad n=1,2,\ldots
\]
we furthermore have
\begin{equation}
  \label{567890}
  (1-\mu_t(a)\mu_t(a)^{\ct})^{1/2}\mu_t(a)(1_{\mathbb H_t^3}-\mu(a)^{\ct}\mu_t(a))^{-1/2}=\mu_t(a).
\end{equation}

\begin{definition}
\label{ba-1}
  Assume that $a\in\mathbb R^4$ satifies \eqref{newcond}. We define
  \begin{equation*}
    b_a(q_t)=\sqrt{1-\mu_t(a)\mu_t(a)^\ct}\star(1-\mu_t(q_t)\mu_t(a)^\ct)^{-\star} \star(\mu_t(q_t)-\mu_t(a))\cdot(I_{\mathbb H_t^3}-\mu_t(a)^\ct\mu_t(a))^{-1/2}
    \end{equation*}
  \end{definition}
  
Similarly, for $a\in\mathbb  R^4$ such that
\begin{equation}
  \label{newcond}
\sum_{j=1}^3\|\mu_{j,t}(a)\|<1
\end{equation}
one can define the square roots
\[
\sqrt{1-\mu_t(a)\mu_t(a)^{[*]}}\quad {\rm and}\quad (I_{\mathbb H_t^3}-\mu_t(a)^{[*]}\mu_t(a))^{-1/2},
\]
where 
\begin{equation*}
\mu_t(a)^{[*]} = \begin{pmatrix}\mu_{1,t}^{{[*]}}&\mu_{2,t}^{{[*]}}&\mu_{3,t}^{[*]}\end{pmatrix}^T,
\end{equation*}
via the power series \eqref{sqrt-power}. Since
\[
\mu_t(a)(\mu_t(a)^{[*]}\mu_t(a))^n=(\mu_t(a)\mu_t(a)^{[*]})^n\mu_t(a),\quad n=1,2,\ldots
\]
we furthermore have
\begin{equation}
  \label{567890}
  (1-\mu_t(a)\mu_t(a)^{[*]})^{1/2}\mu_t(a)(1_{\mathbb H_t^3}-\mu(a)^{[*]}\mu_t(a))^{-1/2}=\mu_t(a).
\end{equation}

  \begin{definition}
        \label{ba-2}
      Assume that $a\in\mathbb R^4$ satifies \eqref{newcond}. We define
    \begin{equation*}
    \beta_a(q_t)=\sqrt{1-\mu_t(a)\mu_t(a)^{[*]}}\star(1-\mu_t(q_t)\mu_t(a)^{[*]})^{-\star} \star(\mu_t(q_t)-\mu_t(a))\cdot(I_{\mathbb H_t^3}-\mu_t(a)^{[*]}\mu_t(a))^{-1/2}
    \end{equation*}
  \end{definition}
   
   \begin{remark} It is worth noting that this analysis can be written equivalently, for the corresponding Fueter $\zeta_t$ variables, and we have the same type of results as in the previous section, but with
  \[
\zeta_t(q_t)=\begin{pmatrix}\zeta_{1,t}&\zeta_{2,t}&\zeta_{3,t}\end{pmatrix}
    \]
    \end{remark}
   
    instead of $\mu_{l,t}$.

      \section{Rational functions}
      \setcounter{equation}{0}
      \label{Rational}
In this section we discuss rational functions in the case of the $q_t,\mu_{u,t}$ and $\zeta_{u,t}$ variables. We begin with a prologue, where the main features of the classical case are reviewed. 
\subsection{Prologue}
The easiest way to define a matrix-valued rational function in one or several complex variables is certainly to characterize it by having entries being quotient of polynomials (in the corresponding number of variables). Another approach was developped in the theory of linear systems, Schur analysis
and operator theory, and is called {\sl the state space method}; see e.g. \cite{MR2663312}. There the emphasis is on the notion of transfer function.
These various concepts make sense and a corresponding theory of linear system can be developped, when the complex numbers are replaced by a
(possibly non-commutative) ring; see for instance \cite{MR839186,MR0452844,SSR}. This provides one motivation for the present section, where
we  consider the counterpart of rational functions in the setting of the $\mathbb H_t$ algebras, both for power series and Fueter-like series. We conclude and remark that rational functions with respect to the $\star$ product will be defined independently of the adjoint, but some of their
properties will be fundamentally different depending on the chosen adjoint.\smallskip

We begin with the following definition, taken from \cite[Definition 2.1 p. 60]{MR2275397}, but of course of a much earlier origin, restricted to
$\mathbb H_t$.

\begin{definition}
  A rational $(\mathbb H_t^{n\times m}$-valued function of $d$ real variables $x_1,\ldots, x_d$, regular at the origin of
  $\mathbb R^d$ is an element of the minimal subring of the  $(\mathbb H_t^{n\times m}$-valued function real analytic in a neighborhood of the
  origin which contains the constant functions the real variables  and is closed under inversion.
  \end{definition} 

  \begin{theorem} (see \cite[Theorem 2.2 p 61]{MR2275397})
    A $\mathbb H_t^{n\times m}$-valued function real analytic in a neighbourhood of the origin of $\mathbb R^N$
    is rational if and only if it can be written as
    \begin{equation}
      \label{real-3}
f(x_1,\ldots, x_d)=D+C(I_{\mathbb H_t^N}-x_1A_1-x_2A_2-\cdots-x_dA_d)^{-1}(x_1B_1+x_2B_2+\cdots +x_dB_d),
\end{equation}
where $N\in\mathbb N$, $D\in \mathbb H_t^{n\times m}$, $ C\in\mathbb H_t^{n\times N}$, and $A_j\in \mathbb H_t^{N\times N},
B_j\in  \mathbb H_t^{N\times m}$, $j=1,\ldots, d$.
  \end{theorem}
  When $N=1$ one rewrites usually \eqref{real-3} as
  \begin{equation}
    \label{bgk-1}
f(x)=D+xC(I_{\mathbb H_t^N}-xA)^{-1}B.
    \end{equation}

  Expressions \eqref{real-3} and \eqref{bgk-1} are called realizations of the corresponding function $f$.When $d>1$, they correspond to transfer
  functions of a linear system of the Fornasini-Marchesini type; see \cite{MR54:2342,MR80c:93028} for the latter.\\  

  In terms of power series expansion \eqref{real-3} can be rewritten as, in multi-index notation,
  \begin{equation}
    f(x)=\sum_{k\in\mathbb N_0^d} x^kf_k
\end{equation}
where
\begin{equation}
f_k=\frac{(|k|-1)^!}{k!}C\begin{pmatrix}k_1A^{k_1-e_1}&k_2A^{k-e_2}&\cdots &k_qA^{k-e_d}\end{pmatrix}B,
  \end{equation}
  where $e_k$, $k=1,\ldots, d$, is the multi-index with all entries equal to $0$, besides the $k$-th one, equal to $1$. For $d=1$ we have
  
      \begin{equation}
f_k=CA^{k-1}B,\quad k=1,2,\ldots
\end{equation}

  In this section, we take $N=1$ and $N=3$, and wish to replace in \eqref{real-3} $x_1=x$ by $q_t$ (for $N=1$) and by the Fueter variables
  $\mu$ and $\zeta$ for $N=3$, and use the $\star$ product rather than the pointwise product. We give equivalent definitions in terms of
  backward-shift operators and related tools.

  \subsection{Rational function: power series in $q_t$}

      \begin{theorem}
        Let $f(q_t)=\sum_{k=0}^\infty q_t^kf_k$ be a converging power series with coefficients $f_n\in\mathbb  H_t^{n\times m}$. Then, the following are equivalent:
        \begin{enumerate}
        \item There exist $N\in\mathbb N_0$ and matrices $(A,B,C)\in\mathbb H_t^{N\times N}\times \mathbb H_t^{N\times m}\times \mathbb H_{t}^{n\times N}$ such that, with $D=f(0)$,
          \begin{equation}
            \label{realstar}
            f(q_t)=D+q_t\star C\star (I_N-q_t\star A)^{-\star}\star B
            \end{equation}

          \item With $A,B,C$ as above
            \begin{equation}
              \label{cab}
              f_k=CA^{k-1}B,\quad k=1,2,\ldots
\end{equation}
\item The real linear span $\mathfrak M(f)$ of the functions $R_0^kfb$, when $n=0,1,\ldots$ and $b$ runs through $\mathbb H_{t}^{N\times 1}$ is a
  finite dimensional real vector space.

\item $f$ is the $\star$ quotient of two matrix-valued polynomials 
  \[
    f=P_1\star P_2^{-\star}=P_3^{-\star}\star P_4
  \]
  such that the denominator polynomials $P_2$ and $P_3$ do not vanish at the origin.
  
  \item The restriction of $f$ to the real line is a rational function with no pole at the origin  and
    with coefficients being block matrices with entries in $\mathbb H_t$. 
\end{enumerate}
      \end{theorem}

      \begin{proof} We proceed into a number of steps.\smallskip \smallskip 

        STEP 1: {\sl $(1)$, $(2)$ and $(5)$ are equivalent.}\smallskip
        
        Assume that $(1)$ holds. Letting $q_t=x\in\mathbb R$ we get the power series expansion of $f(x)$ and obtain \eqref{cab}. Assume $(2)$, then for
        real $x$ near the origin, and with $f(0)=D$,
        \begin{equation}
          \label{rest}
          f(x)= f(0)+\sum_{k=1}^\infty f_kx^k=D+xC(I_N-xA)^{-1}B
        \end{equation}
        and this proves $(5)$. Assume now $(5)$. By the classical theory it admits a realization of the form \eqref{rest},
        which determines $f(q_t)$ in a unique way via the power expansion in power of $x$.\\

        STEP 2: {\sl $(1)$ and $(3)$ are equivalent.}\smallskip

        Assume that $(1)$ holds.  Then
        \[
          R_0^kfb=C\star(I_{\mathbb H_t^N} -q_t\star A)^{-\star} A^{k-1}B,\quad k=1,2,\ldots
        \]
        and so the span in question has dimension at most the span (as a real vector space) of the vectors of the form
        \begin{equation}
          \label{vad}
          A^{k-1}Bb
        \end{equation}
        which is finite, and so $(3)$ holds. Assume now $(3)$. Then we take $v_1,\ldots, v_N$ a basis of $\mathfrak M(f)$, and let
        \[
M(q_t)=\begin{pmatrix}v_1(q_t)&v_2(q_t)&\cdots &v_N(q_t)\end{pmatrix}
\]
Since $\mathfrak M(f)$ is $R_0$-invariant, there is a matrix $A\in\mathbb R^{N\times N}$ such that
\[
  R_0M=MA
\]
from which we get
\[
  M(q_t)=q_t\star M(q_t)A+M(0)
\]
and so\footnote{Since this specific $A$ has real entries we could write $M(0)(I_{\mathbb H_t^N}-q_tA)^{-1}$} $M(q_t)=M(0)\star(I_{\mathbb H_t^N}-q_tA)^{-\star}$.
To conclude we note that for every $u\in\mathbb H_{t}^{N}$ there exists a $v\in\mathbb H_t^{N}$ such that
\[
  R_0fu=M(q_t)v
\]
and so
\[
R_0f=M(q_t)B
\]
for some matrix in $\in\mathbb H_t^{N\times m}$ , and hence $(1)$ holds.\\

STEP 3: {\sl $(3)$ and $(5)$ are equivalent.}\smallskip

Assume $(3)$. Then the restriction of $f$ to the real axis is a rational function of a real variable, with coefficients matrices with entries in
$\mathbb H_t$. It admits thus a realization of the form \eqref{rest}, and so is the restriction of a uniquely determined function of the form \eqref{realstar}.
The converse is proved by remarking that $(I_{\mathbb H_t^N}-xA)^{-1}$ is the quotient of two polynomials.

      \end{proof}
      
      \begin{definition}
The power series $f(q_t)=\sum_{k=0}^\infty q_t^kf_k$ is called rational if it satisfies either of the conditions in the previous theorem.
        \end{definition}

        Using the classical properties of realizations one proves:

        \begin{theorem}
 The sum and product of two rational functions of compatible sizes are rational. Assuming its value at the origin invertible the inverse of a rational function
 is invertible.
          \end{theorem}

          \begin{proof}
            For the inverse, we need to have $n=m$ and $D$ invertible. Then (see e.g. \cite{bgk1}), with $f(x)$ as in \eqref{bgk-1}, we have $\det f(x)\not\equiv 0$ and
            \[
              f(x)^{-1}=D^{-1}-xD^{-1}C(I_{\mathbb H_t^N}-xA^\times)^{-1}BD^{-1},\quad{\rm where}\quad A^\times=A-BD^{-1}C.
            \]
            It follows that the function
            \[
D^{-1}-q_tDC^{-1}\star(I_{\mathbb H_t^N}-q_tA^{\times})^{-\star}BD^{-1}
\]
is the $\star$-inverse of $f(q_t)$. Similarly, still relying on the formulas from \cite{bgk1}, we have (with matrices of compatible sizes)
\[
  (D_1+q_tC_1\star(I_{\mathbb H_t^{N_1}}-q_tA_1)^{-\star}B_1)   \star (D_2+q_tC_2\star(I_{\mathbb H_t^{N_2}}-q_tA_2)^{-\star}B_2)=D+q_tC\star(I_{\mathbb H_t^N}-q_tA)^{-\star}B
\]
with $N=N_1+N_2$ and
\[
  \begin{split}
    A&=\begin{pmatrix}A_1&B_1C_2\\0&A_2\end{pmatrix}\\
    B&=\begin{pmatrix}B_1D_2\\ B_2\end{pmatrix}\\
    C&=\begin{pmatrix}C_1&D_1C_2\end{pmatrix}\\
    D&=D_1D_2.
  \end{split}
  \]
            \end{proof}

          \begin{theorem}
A function is rational if and only if it can be obtained after a finite number of $\star$ operations involving polynomials
            \end{theorem}

          \begin{proof}
The result is true when restricted to the real line and follows by replacing $x$ by $q_t$ and the pointwise product by the $\star$ product.
            \end{proof}

            As we noted above, the notion of rational function is the same for the two adjoints $\circledast$ and $[*]$, but some of their properties are
            completely different, in view of the different properties of the two adjoints.  In the present discussion, in the case of $\ct$ one can go one step
            further:
            
          \begin{theorem}
The power series $f(q_t)=\sum_{k=0}^\infty q_t^kf_k$ is rational if and only if it is the quotient of a matrix-polynomial and of a real scalar polynomial
            \end{theorem}

            \begin{proof}
              Assume first that $f$ is $\mathbb H_t^{a\times b}$-valued and rational. So each entry $f_{ij}$, $i=1,\ldots a$ and $j=1,\ldots b$ is rational. Using Corollary \ref{coro-2-19} write
              \[
                \begin{split}
                  f_{ij}(q_t)&=q_{ij}(q_t)\star(p_{ij}(q_t))^{-\star}\\
                  &=\frac{(q_{ij}\star p_{ij}^{\ct}(q_t))(q_t)}{(p_{ij}\star p_{ij}^\ct)(q_t)},
                \end{split}               
              \]
              The result follows by taking as denominator $P(q_t)=\prod_{i,j}(p_{ij}\star p_{ij}^\ct)(q_t)$
              Each entry $f_{ij}$ can be rewritten as
              \[
                f_{ij}(q_t)=\frac{q_{ij}(q_t)\star\underbrace{ \prod_{\substack{k=1,\ldots, a\\
                    \ell=1,\ldots b\\
                    k\not =i\\
                    \ell\not=j}}{(p_{k\ell}\star p_{k\ell}^\ct)(q_t)}}_{{\rm real\,\, polynomial}}}{P(q_t)}
                \]
\end{proof}

This characterisation is not possible for $[*]$ since $q_tq_t^{[*]}$ is not, in general, a scalar matrix.

\begin{example}
  \label{ex-1}
  The function $\Theta_N$ defined by \eqref{theta_N} is rational with a realization \eqref{theta-!!!} given by \eqref{t-1}-\eqref{t-4}.
  Furthermore, this realization is $\ct$-unitary with respect to ${\mathbf G_{\ct}}$ in the sense that
  \begin{equation}
    \label{unitary}
   \begin{pmatrix}A&B\\C&D\end{pmatrix}^{\ct}\begin{pmatrix}{\mathbf G_{\ct}}&0\\     0&I_{\mathbb H_t^N}\end{pmatrix}\begin{pmatrix}A&B\\C&D\end{pmatrix}=\begin{pmatrix}{\mathbf G_{\ct}}&0\\     0&I_{\mathbb H_t^N}\end{pmatrix}
  \end{equation}
\end{example}

This follows Step 2 in the proof of  Theorem \ref{inter-pol-N}, since one has to verify equalities \eqref{t-890}-\eqref{t-900}.\\

We note that \eqref{unitary} is equivalent to
\[
   \begin{pmatrix}A&B\\C&D\end{pmatrix}^{-1}\begin{pmatrix}{\mathbf G_{\ct}}^{-1}&0\\     0&I_{\mathbb H_t^N}\end{pmatrix}\left(\begin{pmatrix}A&B\\C&D\end{pmatrix}^{\ct}\right)^{-1}=\begin{pmatrix}{\mathbf G_{\ct}}^{-1}&0\\     0&I_{\mathbb H_t^N}\end{pmatrix}
  \]
and so is equivalent to

\begin{equation}
    \label{unitary}
   \begin{pmatrix}A&B\\C&D\end{pmatrix}\begin{pmatrix}{\mathbf G_{\ct}}^{-1}&0\\     0&I_{\mathbb H_t^N}\end{pmatrix}\begin{pmatrix}A&B\\C&D\end{pmatrix}^{\ct}=\begin{pmatrix}{\mathbf G_{\ct}}^{-1}&0\\     0&I_{\mathbb H_t^N}\end{pmatrix},
 \end{equation}
  In particular the matrix ${\mathbf G_{\ct}}^{-1}$ satisfies the Stein equation
\begin{equation}
  \label{stein-2}
{\mathbf G_{\ct}}^{-1}-A{\mathbf G_{\ct}}^{-1}A^{\ct}=BB^{\ct}.
    \end{equation}

    The following result hints to the fact that the theory of structured rational functions can be extended to the present setting; as was shown by a counterexample in \cite[(62.38) p. 1767]{acs-survey} one cannot directly force conditions on the boundary, but one has to consider conditions in
    terms of multiplication
    operators, or symmetric on the real line replacing the metric boundary conditions by operator-type conditions.

\begin{proposition}
  In the notation of \eqref{unitary} define $\Theta$ by \eqref{theta-!!!},
  \[
    \Theta(q_t)=D+q_tC\star(I_{\mathbb H_t^N}-q_tA)^{-\star}B.
  \]
  Then the operator of $\star$ multiplication by $\Theta$ on the right is an isometry from the Hardy space $\mathbf H_{2,\ct}$ into itself.
\end{proposition}

\begin{proof}
The  proof is the same as the proof of Step 5 in the proof of Theorem \ref{inter-pol-N}, where only equality \eqref{unitary} is used, and not the specific form of $A,B,C,D$.
  \end{proof}

  Similar considerations hold for the second adjoint. Similar claims hold for the Blaschke factor associated to $[*]$, as we now check, and for
  the Blaschke factors associated to the Fueter variables; see next subsection for the latter.
  
  \begin{proposition}
    The formula
    \begin{eqnarray}
      \label{real1}
      A&=&a^{[*]}\\
      B&=&\Gamma_\alpha^{-1}K_\alpha\\
      C&=&1\\
      D&=&-\alpha K_\alpha
           \label{real4}
    \end{eqnarray}
    provide a realization of $B_\alpha$, which satisfies moreover
    \begin{eqnarray}
      \label{uni-0001}
        \begin{pmatrix}A&B\\C&D\end{pmatrix}^{[*]}\begin{pmatrix}\Gamma_\alpha&0\\ 0&1\end{pmatrix}
        \begin{pmatrix}A&B\\C&D\end{pmatrix}&=&\begin{pmatrix}\Gamma_\alpha&0\\ 0&1\end{pmatrix}
      \\
      \label{uni-0002}
\begin{pmatrix}A&B\\C&D\end{pmatrix}\begin{pmatrix}\Gamma_\alpha^{-1}&0\\ 0&1\end{pmatrix}\begin{pmatrix}A&B\\C&D\end{pmatrix}^{[*]}&=&
        \begin{pmatrix}\Gamma_\alpha^{-1}&0\\ 0&1\end{pmatrix}
  \end{eqnarray}
  \end{proposition}

  \begin{proof}
    From \eqref{aveiro-2024}
    we have
    \[
    \begin{split}
      B_\alpha(q_t)&=-\alpha K_\alpha+\sum_{n=1}^\infty q_t^n\alpha^{[*](n-1)}\Gamma_\alpha^{-1}K_\alpha\\
      &=-\alpha K_\alpha+q_t\star (1-q_t\alpha^{[*]})^{-\star}\Gamma_\alpha^{-1}K_\alpha,
    \end{split}
    \]
    and hence we get the realization \eqref{real1}-\eqref{real4}. To check \eqref{uni-0001} amounts to check the following three equalities
    \begin{eqnarray}
      \alpha\Gamma_\alpha\alpha^{[*]}+1&=&\Gamma_\alpha,\\
      \alpha K_\alpha-\alpha K_\alpha&=&0,\\
      K_\alpha\Gamma_\alpha^{-1}K_\alpha+     K_\alpha a^{[*]}aK_\alpha&=&1.
    \end{eqnarray}
    The first one is the Stein equation \eqref{stein-4} and the second one is a tautology. The third one can be rewritten as
    \[
\Gamma_\alpha^{-1}+a^{[*]}a=K_\alpha^{-2}=L_\alpha^{-1}
\]
and is just \eqref{567}. To get \eqref{uni-0002} we rewrite \eqref{uni-0001} as
\[
          \begin{pmatrix}A&B\\C&D\end{pmatrix}^{[*]}\begin{pmatrix}\Gamma_\alpha&0\\ 0&1\end{pmatrix}
          \begin{pmatrix}A&B\\C&D\end{pmatrix}\begin{pmatrix}\Gamma_\alpha^{-1}&0\\ 0&1\end{pmatrix}=I_{\mathbb H_t^4},
          \]
          so that
          \[
            \begin{pmatrix}A&B\\C&D\end{pmatrix}
            \begin{pmatrix}\Gamma_\alpha^{-1}&0\\ 0&1\end{pmatrix}\begin{pmatrix}A&B\\C&D\end{pmatrix}^{[*]}
            \begin{pmatrix}\Gamma_\alpha&0\\ 0&1\end{pmatrix}
            =I_{\mathbb H_t^4},
          \]
          and the result follows.
        \end{proof}

        For the above computations in the setting of complex matrices see the arXiv preprint \cite[Theorem 3.12]{alpay2024schur}.
\subsection{Rational functions for Fueter variables}
We will consider this topic in greater details in a future publication but mention the following:
For the $\mu_t$ variables, the $\star$-product is defined by
\begin{equation}
  \mu_t^\alpha u\star\mu_t^\beta v=\mu_t^{\alpha+\beta}uv,\quad \alpha,\beta\in\mathbb N_0^3,\,\ u,v\in\mathbb H_t.
\end{equation}
The corresponding rational functions regular at the origin are the expressions of the form
  \begin{equation}
    \label{formula-09876}
  f(x)=D+C\star(I_{\mathbb H_t^N}-\sum_{k=1}^3\mu_{k,t}A_k)^{-\star}\star(\sum_{k=1}^3\mu_{k,t}B_k).
\end{equation}

With
\[
  \mu_t=\begin{pmatrix}\mu_{1,t}&\mu_{2,t}&\mu_{3,t}\end{pmatrix}
\]
and
\[
A=\begin{pmatrix}A_1\\ A_2\\  A_3\end{pmatrix},\quad B=\begin{pmatrix}B_1\\ B_2\\ B_3\end{pmatrix},
\]
we rewrite \eqref{real-3} as
\begin{equation}
f(x)=D+C\star(I_{\mathbb H_t^N}-(\mu_t \otimes I_N) A)^{-1})^{-\star}(\mu_t\otimes I_N)B
  \end{equation} 
which we write as
  \begin{equation}
f(x)=D+C\star(I_{\mathbb  H_t^N}-\mu_t A)^{-1})^{-\star}\mu_t B
  \end{equation}
  Note that this function is not defined at the origin of $\mathbb R^4$.\\
  
 For the $\zeta_t$ variables,  the $\star$-product is defined by
\begin{equation}
  \zeta_t^\alpha u\star\zeta_t^\beta v=\zeta_t^{\alpha+\beta}uv,\quad \alpha,\beta\in\mathbb N_0^3,\,\ u,v\in\mathbb H_t.
  \end{equation}
  \begin{equation}
    \label{formula-09876}
  f(x)=D+C\star\left(I_{\mathbb H_t^N}-\sum_{k=1}^3\zeta_{k,t}A_k\right)^{-\star}\star\left(\sum_{k=1}^3\zeta_{k,t}B_k\right).
\end{equation}

The formulas are the same as for the previous subsection, and were developped in \cite{MR2124899} for the case $t=-1$.\smallskip

In both cases the corresponding Blaschke factors \eqref{ba-1} and \eqref{ba-2} are rational. More precisely, and using \eqref{567890}, we can rewrite $b_a(q_t)$ as (see also \cite[(4.4) p. 12]{akap1})
\begin{equation}
  \label{ba-1-1}
  \begin{split}
    b_a(x)&=-\mu_t(a)+\\
    &\hspace{5mm}
    +(1-\mu_t(a)\mu_t(a)^{\ct})^{1/2} \star(1-\mu_t(x)\mu_t(a)^{\ct})^{-\star}\star \mu_t(a)(I_{\mathbb H_t^3}-\mu(a)_t^{\ct}\mu_t(a))^{1/2}.
    \end{split}
\end{equation}
The realization is equal to
\begin{equation}
  \begin{split}
    A&=\mu_t(a)^{\ct}\\
    B&=(I_{\mathbb H_t^3}-\mu_t(a)^{\ct}\mu_t(a))^{1/2}\\
    C&=(1-\mu_t(a)\mu_t(a)^{\ct})^{1/2}\\
    D&=-\mu_t(a),
  \end{split}
\end{equation}
where, as earlier, the squareroots are defined by the power expansion of $\sqrt{1-t}$ for $a\in\mathbb R^4$ such that
$\|\mu_t(a)\mu_t(a)^{\ct}\|_{op}<1$.
We repeat the argument: with $\sqrt{1-t}=1+\sum_{n=1}^\infty c_nt^n$, where $c_1,c_2,\ldots \in\mathbb R$ (see also \eqref{sqrt-power}), the  series
\[
  1+\sum_{n=1}^\infty c_n( \mu_t(a)\mu_t(a)^{\ct})^n\quad {\rm and}\quad I_{\mathbb H_t^3}+\sum_{n=1}^\infty c_n\mu_t(a)^{\ct}(\mu_t(a)\mu_t(a)^{\ct})^{n-1}
\mu_t(a)
\]
converge since the operator norm is sub-multiplicative.\smallskip

  The realization is ${\ct}$-unitary in the sense that
  \[
    \begin{pmatrix}A&B\\C&D\end{pmatrix}\begin{pmatrix}A&B\\C&D\end{pmatrix}^{\ct}=\begin{pmatrix}A&B\\C&D\end{pmatrix}^{\ct}
    \begin{pmatrix}A&B\\C&D\end{pmatrix}=I_{\mathbb H_t^{4\times 4}}.
  \]
  Equality \eqref{567890} is used to prove this equality.\\
  
  \begin{remark}
    Besides the case just encountered, there are three other cases for the Fueter variables, namely the classical Fueter variables, for the
    two adjoints, and the new Fueter variables, for the $[*]$ adjoint.
    \end{remark}
  
  
   
\appendix

   \section{Pontryagin and Krein spaces}
   \label{append_krein}
   \setcounter{equation}{0}
   Having in view the needs for the present work we briefly review the main relevant definitions and results on Pontryagin spaces and Krein spaces.
   For further information we refer in particular to \cite{azih,bognar,MR92m:47068,ikl}. Let therefore $(\mathcal V,[\cdot,\cdot])$ be a vector space
   on $\mathbb K$ where $\mathbb K=\mathbb R$ or $\mathbb K=\mathbb C$, endowed
   with a symmetric form (when $\mathbb K=\mathbb R$) or an Hermitian form (when $\mathbb K=\mathbb R$).
The discussion below is done for the complex setting and we choose the convention
\begin{equation}
  \label{hermi-form}
[av,bw]_{\mathcal V}=\overline{b}[v,w]_{\mathcal V}a,\quad v,w\in\mathcal V,\quad {\rm and}\quad a,b\in\mathbb C,
     \end{equation}
for the anti-linear variable. The real case is directly adapted. One can also consider the case of quaternions, but we prefer to include it in the
following section, pertaining to modules.

\begin{definition}
  The space $(\mathcal V,[\cdot,\cdot]_{\mathcal V})$ is called a Krein space if it can be written in the direct and orthogonal sum
  \begin{equation}
    \label{vwv}
    \mathcal V=   \mathcal V_+\stackrel{\cdot}{+}   \mathcal V_-
  \end{equation}
  where
  \begin{enumerate}
  \item The spaces $(   \mathcal V_+,[\cdot,\cdot]_{\mathcal V})$ and $(   \mathcal V_-,-[\cdot,\cdot]_{\mathcal V})$ are Hilbert spaces.

  \item
The sum \eqref{vwv} is direct, meaning that $   \mathcal V_+\cap    \mathcal V_-=\left\{0\right\}$.
\item
The sum  \eqref{vwv} is orthogonal, meaning that
\begin{equation}
[v_+,v_-]_{\mathcal V}=0,\quad \forall v_+\in\mathcal V_+ \quad and\quad v_-\in   \mathcal V_-.
\end{equation}
    \end{enumerate}
    
  \end{definition}

  \eqref{vwv} is called a {\it fundamental decomposition} and is not unique, but in the case where one of the component reduces to the zero vector
  space. Every fundamental decomposition generates a Hilbert space structure on $\mathcal V$  via
  \begin{equation}
\langle v,w\rangle=[v_+,w_+]_{\mathcal V}-[v_-,w_-]_{\mathcal V},
\end{equation}
where $v=v_++v_-$ and $w=w_++w_-$ are the corresponding fundamental decompositions of $v,w\in\mathcal V$.

\begin{proposition} (see \cite[p. 102 and Corollary IV.6.3 p. 92]{bognar})
  The space $(\mathcal V,\langle\cdot,\cdot\rangle)$ is a Hilbert space, and norms arising from orthogonal decompositions are equivalent,
  and define therefore the same topology on $\mathcal V$.
\end{proposition}

In particular, Riesz  theorem for the representation of bounded functions still holds in Krein spaces and the notion of bounded point evaluation in a Krein space of
functions leads to the notion of reproducing kernel Krein space.\smallskip

When $\mathcal V_-$ is finite dimensional the space $\mathcal V$ is called a Pontryagin space, and ${\rm dim}\, \mathcal V_-$ is its index of
negativity.

\begin{example}
  Let $J\in\mathbb C^{n\times n}$ satisfy $J=J^*=J^{-1}$ (a {\sl signature matrix}), and define on $\mathbb C^n$
  \begin{equation}
    [z,w]=w^*Jz,\quad z,w\in\mathbb C^n.
  \end{equation}
  Then, $(\mathbb C^n, [\cdot,\cdot]_J)$ is a finite dimensional Pontryagin space with negativity index equal to the number of negative eigenvalues of $J$.
  The adjoint $A^{[*]}$ of the matrix $A\in\mathbb C^{n\times n}$ with respect to this adjoint is
  \begin{equation}
    A^{[*]}=JA^*J.
  \end{equation}
  \label{Cnpont}
  \end{example}

One should be aware that there may be more than one given Krein (as opposed to Hilbert or Pontryagin) space of functions with a given reproducing kernel; this was proved in Laurent Schwartz paper \cite[\S13]{schwartz} and another counterexample can be found in \cite[\S4]{a2}.\\

In our work Pontryagin, and Hilbert spaces appear in a natural way as described in Section~\ref{Algebraic-Section+Theory-Functions} and Krein spaces appear in Section~\ref{Hardy_spaces_cdast} and beyond, see, for example, Theorem~\ref{ct-Krein}.
  


   \section{Pontryagin and Krein right $\mathcal H_2^t$-modules}
   \label{Pontryagin_Q_t}
   We assume now that $\mathcal V$ is a right $\mathcal H_2^t$-module, and now replace \eqref{hermi-form} by
   \begin{equation}
  \label{hermi-form-1}
  [vq,wp]_{\mathcal V}
  =p^{\circledast}[v,w]_{\mathcal V}q,\quad v,w\in\mathcal H,\quad {\rm and}\quad p,q\in\mathcal H_2^t,
     \end{equation}

     or

     \begin{equation}
  \label{hermi-form-2}
[vq,wp]_{\mathcal V}=p^{[*]}[v,w]_{\mathcal V}q,\quad v,w\in\mathcal H,\quad {\rm and}\quad p,q\in\mathcal H_2^t,
     \end{equation}
depending on the chosen adjoint.

\begin{definition}
  Let $\mathcal K$ be a vector space which is also a right $\mathcal H_2^t$-module, and let $[\cdot,\cdot]_{\mathcal K}$ be a
  $\mathcal H_2^t$-valued form which satisfies \eqref{hermi-form-1} or \eqref{hermi-form-2}. Then, $\mathcal K$ is called a $\mathcal H_2^t$ module if
  it is moreover a real Krein space when endowed with the symmetric form
\[
\langle v,w\rangle={\rm Tr}\, \underbrace{[v,w]_{\mathcal K}}_{\in\mathcal H^t_2}.
\]
\end{definition}

\begin{definition}
  In the notation of the previous definition, assume that $\mathcal K$ is a space of $\mathcal H_2^t$-valued
  functions defined on a set $\Omega$. Then, $\mathcal K$ is a
  reproducing kernel Krein space if there exists a function $K(z,w)$ defined for $z,w\in\Omega$ and with the following property
  \begin{equation}
    \label{point-eval-2}
    {\rm Tr}\, \underbrace{[f(\cdot), K(\cdot, w)p]_{\mathcal K}}_{\in\mathcal H^t_2}
    ={\rm Tr}\,p^{\circledast}f(w)     ,\quad w\in \Omega,\quad p\in\mathcal H_2^t,
   \end{equation}
   or
   \begin{equation}
{\rm Tr}\, \underbrace{[f(\cdot), K(\cdot, w)p]_{\mathcal K}}_{\in\mathcal H^t_2} ={\rm Tr}\, p^{[*]}f(w),\quad w\in \Omega,\quad p\in\mathcal H_2^t,
       \end{equation}
   depending on the chosen adjoint.
     \end{definition}

\begin{proposition}
  In the previous notation, it holds that
  \begin{equation}
    \label{point-eval}
       [f(\cdot), K(\cdot, w)b]_{\mathcal K}=b^{\circledast}f(w),\quad w\in \Omega,\quad b\in\mathcal H_2^t
\end{equation}
     or
\begin{equation}
       [f(\cdot), K(\cdot, w)b]_{\mathcal K}
       =b^{[*]}f(w),\quad w\in \Omega,\quad b\in\mathcal H_2^t
       \end{equation}
     depending on the chosen adjoint.
\end{proposition}

\begin{proof} We consider the case of the first adjoint; the case of $[*]$ is treated similarly. From \eqref{hermi-form-1} we can rewrite
  \eqref{point-eval-2} as

  \begin{equation}
    \label{point-eval-2-3}
      {\rm Tr}\, b^{\ct}[f(\cdot), K(\cdot, w)]_{\mathcal K}={\rm Tr}\,b^{\circledast}f(w)     ,\quad w\in \Omega,\quad b\in\mathcal H_2^t.
   \end{equation}
To conclude we use Corollary \ref{help-facto}.

  \end{proof}

  These definitions are readily extended to the case of matrix-valued functions with entries in $\mathcal H_2^t$. As already earlier in
  Section \ref{Hardy_q_t}, the matrix-valued case already has differences with the scalar case in the quaternionic setting; see \cite[p. 1767 (and in particular (62.38)]{acs-survey}.\\

    \bibliographystyle{plain}
  \def\lfhook#1{\setbox0=\hbox{#1}{\ooalign{\hidewidth
  \lower1.5ex\hbox{'}\hidewidth\crcr\unhbox0}}} \def\cprime{$'$}
  \def\cprime{$'$} \def\cprime{$'$} \def\cprime{$'$} \def\cprime{$'$}
  \def\cprime{$'$}

\bibliography{all}

\end{document}